\documentclass[a4paper,onecolumn, superscriptaddress,10pt, accepted=2022-03-14, issue=4, volume=4, shorttitle=papers/compositionality-4-4]{compositionalityarticle}

\pdfoutput=1

\usepackage[utf8]{inputenc}
\usepackage[english]{babel}
\usepackage[T1]{fontenc}

\newif\ifignore 
\ignorefalse

\newcommand{\auxproof}[1]{
\ifignore\mbox{}\newline
\textbf{PROOF:} \dotfill\newline
{\it #1}\mbox{}\newline
\textbf{ENDPROOF}\dotfill
\fi}

\usepackage{amsmath,amssymb,amstext,mathtools}
\usepackage{stmaryrd}
\usepackage{xspace}
\usepackage{hyperref}
\usepackage{url}
\usepackage{nccmath}
\usepackage{fancybox}
\usepackage{pdfrender}
\usepackage{multirow}

\usepackage{xypic}
\usepackage[all]{xy}
\xyoption{all}
\xyoption{2cell}
\UseTwocells
\xyoption{tips}
\SelectTips{cm}{}  

\newenvironment{myproof}[1][Proof]%
   { \begin{trivlist}%
     \item[\hskip \labelsep {\bfseries #1}.]%
   }%
   { \end{trivlist}%
   }

\newtheorem{theorem}{Theorem}
\newtheorem{lemma}[theorem]{Lemma}
\newtheorem{proposition}[theorem]{Proposition}
\newtheorem{corollary}[theorem]{Corollary}

\newtheorem{definition}{Definition}
\newtheorem{example}{Example}

\newtheorem{remark}{Remark}
\newtheorem{fact}{Fact}

\usepackage{wrapfig}
\newcommand*{\fatten}[1][.4pt]{%
  \textpdfrender{
    TextRenderingMode=FillStroke,
    LineWidth={\dimexpr(#1)\relax},
  }%
}

\ifdefined\mathsl\else
  \DeclareMathAlphabet{\mathsl}{\encodingdefault}{\rmdefault}{\mddefault}{\sldefault}
  \SetMathAlphabet{\mathsl}{bold}{\encodingdefault}{\rmdefault}{\bfdefault}{\sldefault}
\fi

\makeatletter
\newcommand{\mathoverlap}[2]{\mathpalette\mathoverlap@{{#1}{#2}}}
\newcommand{\mathoverlap@}[2]{\mathoverlap@@{#1}#2}
\newcommand{\mathoverlap@@}[3]{\ooalign{$\m@th#1#2$\crcr\hidewidth$\m@th#1#3$\hidewidth}}
\makeatother

\newcommand{\QEDbox}{\ensuremath{\square}}
\newcommand{\QED}{\hspace*{\fill}$\QEDbox$}

\newcommand{\klafter}{\mathrel{\bullet}}
\newcommand{\idmap}[1][]{\ensuremath{\mathrm{id}_{#1}}}
\newcommand{\after}{\mathrel{\circ}}
\newcommand{\supp}{\mathrm{supp}}

\newcommand{\NNO}{{\mathbb{N}}}
\newcommand{\set}[2]{\{#1\;|\;#2\}}
\newcommand{\setin}[3]{\{#1\in#2\;|\;#3\}}
\newcommand{\tuple}[1]{\langle#1\rangle}
\newcommand{\setsize}[1]{|{\kern.1em}#1{\kern.1em}|}
\newcommand{\bigsetsize}[1]{\big|{\kern.1em}#1{\kern.1em}\big|}
\newcommand{\no}[1]{#1^{\scriptscriptstyle \bot}} 

\newcommand{\mulnom}{\ensuremath{\mathsl{mn}}}
\newcommand{\multinomial}[1][]{\ensuremath{\mulnom[#1]}}
\newcommand{\negmultinomial}{\ensuremath{\mathsl{nmn}}}
\newcommand{\pol}{\ensuremath{\mathsl{pl}}}
\newcommand{\polya}[1][]{\ensuremath{\pol[#1]}}
\newcommand{\negpolya}{\ensuremath{\mathsl{npl}}}

\newcommand{\binomial}[1][]{\ensuremath{\mathsl{bn}[#1]}}
\newcommand{\hypgeom}{\ensuremath{\mathsl{hg}}}
\newcommand{\hypergeometric}[1][]{\ensuremath{\hypgeom[#1]}}
\newcommand{\mnff}{\ensuremath{\mathsl{mnff}}}
\newcommand{\hgff}{\ensuremath{\mathsl{hgff}}}
\newcommand{\plff}{\ensuremath{\mathsl{plff}}}
\newcommand{\mnclg}{\ensuremath{\mathsl{MN}}}
\newcommand{\hgclg}{\ensuremath{\mathsl{HG}}}
\newcommand{\plclg}{\ensuremath{\mathsl{PL}}}
\newcommand{\negmnclg}{\ensuremath{\mathsl{NMN}}}
\newcommand{\neghgclg}{\ensuremath{\mathsl{NHG}}}
\newcommand{\negplclg}{\ensuremath{\mathsl{NPL}}}
\newcommand{\neghypergeometric}{\ensuremath{\mathsl{nhg}}}
\newcommand{\negbinomial}[1][]{\ensuremath{\mathsl{nbn}[#1]}}
\newcommand{\mean}{\ensuremath{\mathsl{mean}}}

\newcommand{\facto}[1]{\ensuremath{#1{\kern-2.5pt}\raisebox{-2.5pt}{\includegraphics[width=0.9em]{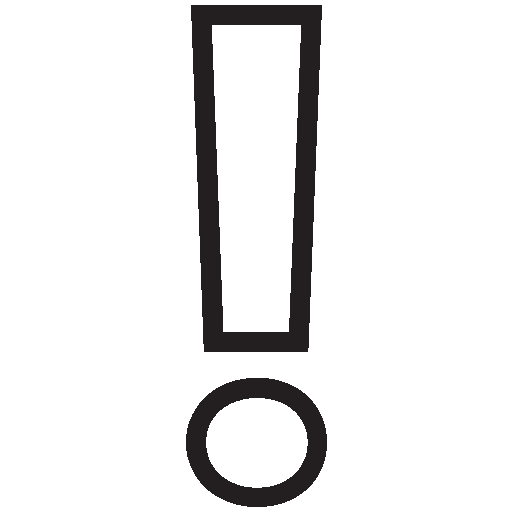}}}}
\newcommand{\sfacto}[1]{\ensuremath{#1{\kern-1.5pt}\raisebox{-1.5pt}{\includegraphics[width=0.6em]{exclamation}}}}
\newcommand{\coefm}[1]{\ensuremath{\fatten[0.6pt]{(}#1\fatten[0.6pt]{)}}}

\newcommand{\bibinom}[2]{\left({\kern-3pt}\binom{#1}{#2}{\kern-3pt}\right)}

\newcommand{\acc}{\ensuremath{\mathsl{acc}}}

\newcommand{\flrn}{\ensuremath{\mathsl{Flrn}}}

\newcommand{\distributionsymbol}{\mathcal{D}}

\newcommand{\multisetsymbol}{\mathcal{M}}

\newcommand{\Dst}{\distributionsymbol}
\newcommand{\infDst}{\ensuremath{\mathcal{D}_{\infty}}}

\newcommand{\Mlt}{\multisetsymbol}
\newcommand{\natMlt}{\multisetsymbol}
\newcommand{\nenatMlt}{\multisetsymbol_{*}}

\newcommand{\one}{\ensuremath{\mathbf{1}}}
\newcommand{\zero}{\ensuremath{\mathbf{0}}}

\newcommand{\Cat}[1]{\ensuremath{\mathbf{#1}}}

\newcommand{\Kl}{\mathcal{K}{\kern-.4ex}\ell}
\newcommand{\EM}{\mathcal{E}{\kern-.4ex}\mathcal{M}}

\newcommand{\Sets}{\Cat{Sets}\xspace}

\newcommand{\intd}{{\kern.2em}\mathrm{d}{\kern.03em}}

\newcommand{\ket}[1]{\ensuremath{|{\kern.1em}#1{\kern.1em}\rangle}}
\newcommand{\ketstrut}{\vrule height 8.5pt depth 4.5pt width 0pt}
\newcommand{\bigket}[1]{\ensuremath{\big|{\kern.1em}#1{\kern.1em}\big\rangle}}
\newcommand{\Bigket}[1]{\ensuremath{\left|\ketstrut{\kern.1em}#1{\kern.005em}\right>}}

\DeclareSymbolFont{T1op}{T1}{cmr}{m}{n}
\SetSymbolFont{T1op}{bold}{T1}{cmr}{bx}{n}
\DeclareMathSymbol{\mathguilsinglleft}{\mathopen}{T1op}{'016}
\DeclareMathSymbol{\mathguilsinglright}{\mathclose}{T1op}{'017}

\newcommand{\eg}{\textit{e.g.}\xspace}

\begin{document}
\title{Urns \& Tubes}

\author{Bart Jacobs}

\affiliation{Institute for Computing and Information Sciences,
Radboud University Nijmegen, P.O. Box 9010, 6500 GL Nijmegen, The Netherlands}
\email{b.jacobs@cs.ru.nl}

\begin{abstract}
Urn models play an important role to express various basic ideas in
probability theory. Here we extend this urn model with tubes. An urn
contains coloured balls, which can be drawn with probabilities
proportional to the numbers of balls of each colour. For each colour a
tube is assumed. These tubes have different sizes (lengths). The idea
is that after drawing a ball from the urn it is dropped in the tube of
the corresponding colour. We consider two associated probability
distributions. The \emph{first-full} distribution on colours gives for
each colour the probability that the corresponding tube is full first,
before any of the other tubes. The \emph{negative} distribution on
natural numbers captures for a number $k$ the probability that all
tubes are full for the first time after $k$ draws.

This paper uses multisets to systematically describe these first-full
and negative distributions in the urns \& tubes setting, in fully
multivariate form, for all three standard drawing modes (multinomial,
hypergeometric, and P\'olya).
\end{abstract}


\section{Introduction}\label{IntroSec}

Consider the situation sketched below~\eqref{UrnTubesPic}, with an urn
filled with coloured balls (on the left) and tubes of different
lengths (on the right), with one tube for each colour. Below there are
three colours: red (R), blue (B), and green (G), but in general there
can be $N\geq 2$ many colours --- and then also $N$ tubes. We consider
the following action: when a ball is drawn from the urn, it is dropped
in the tube of the corresponding colour. This action is repeated. 
\begin{equation}
\label{UrnTubesPic}
\vcenter{\hbox{
\begin{picture}(250,100)
\thicklines
\put(68, 60){\oval(20, 20)[tl]}
\put(0, 60){\oval(20, 20)[tr]}
\put(20, 10){\oval(20, 20)[bl]}
\put(48, 10){\oval(20, 20)[br]}
\put(20, 0){\line(1, 0){28}}
\put(10, 10){\line(0, 1){50}}
\put(58, 10){\line(0, 1){50}}
\color{green}
\put(20,8){\circle{12}}
\color{black}
\put(16,5){G}
\color{red}
\put(34,8){\circle{12}}
\color{black}
\put(31,5){R}
\color{blue}
\put(48,8){\circle{12}}
\color{black}
\put(45,5){B}
\color{red}
\put(20,22){\circle{12}}
\color{black}
\put(17,19){R}
\color{blue}
\put(34,22){\circle{12}}
\color{black}
\put(31,19){B}
\color{red}
\put(48,22){\circle{12}}
\color{black}
\put(45,19){R}
\color{red}
\put(20,36){\circle{12}}
\color{black}
\put(17,33){R}
\color{blue}
\put(34,36){\circle{12}}
\color{black}
\put(31,33){B}
\color{green}
\put(48,36){\circle{12}}
\color{black}
\put(44,33){G}
\put(128,0){Red}
\put(120, 10){\line(1, 0){30}}
\put(127, 10){\line(0, 1){41}}
\put(143, 10){\line(0, 1){41}}
\color{red}
\put(135,17){\circle{12}}
\color{black}
\put(132,14){R}
\put(176,0){Blue}
\put(170, 10){\line(1, 0){30}}
\put(177, 10){\line(0, 1){83}}
\put(193, 10){\line(0, 1){83}}
\color{blue}
\put(185,17){\circle{12}}
\color{black}
\put(182,14){B}
\color{blue}
\put(185,31){\circle{12}}
\color{black}
\put(182,28){B}
\color{blue}
\put(185,45){\circle{12}}
\color{black}
\put(182,42){B}
\put(223,0){Green}
\put(220, 10){\line(1, 0){30}}
\put(227, 10){\line(0, 1){70}}
\put(243, 10){\line(0, 1){70}}
\color{green}
\put(235,17){\circle{12}}
\color{black}
\put(231,14){G}
\color{green}
\put(235,31){\circle{12}}
\color{black}
\put(231,28){G}
\color{green}
\put(235,45){\circle{12}}
\color{black}
\put(231,42){G}
\color{green}
\put(235,60){\circle{12}}
\color{black}
\put(231,57){G}
\put(80, 90){\vector(1, 0){40}}
\put(65, 95){draws from the urn}
\end{picture}}}
\end{equation}

In this paper we consider this urns \& tubes setting in two
scenarios, involving either \emph{some tube} or \emph{all tubes} being
full for the first time. They both start from empty tubes.
\begin{enumerate}
\item The first scenario looks at the probability that some tube is
  completely full first, before any of the other tubes is full. As
  will be shown, this yields a distribution on colours, which we call
  the \emph{first-full} distribution.

\item In the second scenario we consider a distribution on natural
  numbers, where the probability assigned to number $k$ is the
  probability that all tubes are full for the first time after $k$
  draws. This means that there is some tube getting full at stage $k$,
  while all other tubes are already full --- possibly with overflows.
  Such distributions are known in the literature as \emph{negative}
  distributions. More on this at the end of this section.
\end{enumerate}

\noindent One can translate this abstract urns \& tubes setting to
more practical scenarios where the filling of the tubes may represent
something good or bad, like hospital beds of various types becoming
fully used up. Both the first-full and the negative scenarios may be
relevant in risk modeling, where the fullness probabilities of tubes
correspond to risks of reaching thresholds.

Intuitively, the first-full probability for a colour $C$ decreases
with the length of the $C$-coloured tube, and increases with the
proportion of $C$-coloured balls in the urn. It thus involves complex
dependencies. The main technical challenge is to prove that first-full
is actually a distribution, with first-full probabilities for each
colour adding up to one. For this purpose we reason compositionally,
via certain probabilistic automata, namely Markov models with output
(MOO), which terminate at some stage, after some number of
compositions, producing the relevant first-full distribution. The same
type of automaton can be used for negative distributions. The
categorical details behind this composition are explained in the
appendix

Commonly three modes of drawing balls from an urn are distinguished,
see \textit{e.g.}~\cite{JohnsonK77,Mahmoud08,PishroNik14,Ross18},
called \emph{multinomial}, \emph{hypergeometric} and \emph{P\'olya};
the last mode is also called \emph{P\'olya--Eggenberg} or
\emph{Dirichlet-multinomial}, see~\cite{JohnsonKB97}. We use ``0'',
``-1'' and ``+1'' as short-hand descriptions for these different
modes.  Explicitly, we use:
\begin{itemize}
\item[``0''] for the multinomial mode, in which the drawn ball is
  returned to the urn;

\item[``-1''] for the hypergeometric mode, where the drawn ball is not
  returned to the urn --- so that the urn is diminishing;

\item[``+1''] for the P\'olya mode, where the drawn ball is returned
  to the urn, together with an additional copy, of the same colour
  (called a reinforcement).
\end{itemize}

\noindent The distinction between multinomial and hypergeometric modes
is most familiar and is often expressed in terms of: \emph{with} or
\emph{without} replacement. The P\'olya mode is less well known.  The
additional ball that is added to the urn after drawing has a
strengthening effect that can capture situations with a cluster
dynamics, like in the spread of contagious diseases~\cite{HayhoeAG17}
or the flow of tourists~\cite{LauKW20}.

In a physical explanation of the first-full distributions we need for
the multinomial and P\'olya modes an auxiliary box of balls on the
side (with sufficiently many balls). In multinomial (resp.\ P\'olya)
mode, the ball drawn from the urn is dropped in the right tube, but
one (resp.\ two) ball(s) of the same colour are taken from the box and
added to the urn, before the next draw. In the P\'olya case the urn
grows in size with each draw. In the multinomial case the urn remains
the same and is best described as a probability distribution (over the
set of colours). In the hypergeometric mode the urn decreases in size;
we thus have to assume that initially the urn contains sufficiently
many balls of each colour: more than the length of each tube.

The urn \& tubes set-up as introduced here generalises the famous
`problem of points', studied in the 17th century by Pierre Fermat and
Blaise Pascal, that played an important role in the development of
modern probability theory --- in particular for the notion of
expectation. See~\cite{Edwards82} for a historical and~\cite{Ma16} for
a popular account. The problem of points involves a game between two
players that is terminated prematurely and where the stakes so far
have to be divided between the two. The solution there is to look at
the remaining number of steps for winning (and associated
probabilities) for each of the players. These remaining steps
translate directly into lengths of two tubes, one for each player,
with multinomial draws, see Subsection~\ref{PointsSubsec} below for
further details.

The current urn \& tubes set-up `inverts' the problem of points, and
also generalises it in two ways: (1)~urns \& tubes are analysed is in
fully multivariate form, and~(2) the analysis covers the three drawing
modes described above. 


Although the same urn \& tubes setting is used both for first-full and
for negative distributions, these distributions are really
different. First of all, negative distributions have the natural
numbers $\NNO$ as sample space, with infinite support --- in
multinomial and P\'olya mode. To $k\in\NNO$ the probability is
assigned that \emph{all} tubes are full, for the first time, after $k$
draws --- where some tubes may overflow, while others are not full
yet. Such negative distributions are studied in the literature, see
\textit{e.g.}~\cite{Panaretos81,SchusterS87,SibuyaYS64} and the
textbooks~\cite{JohnsonKK05,JohnsonKB97}, but typically with one tube
only. Still, these negative distributions are not mainstream, and are
even called `forgotten', see~\cite{MillerF07}. Here we describe
negative distributions, in the general urn \& tubes setting, in fully
multivariate form, for all three drawing modes (``0'', ``-1'',
``+1'').

In this paper we make extensive use of multisets, like in other recent
publications~\cite{Jacobs19b,Jacobs21b,Jacobs21a,JacobsS20}. A
multiset is like a subset, except that elements may occur multiple
times. For instance, an urn is a multiset, over the set of colours. A
draw of multiple balls from such an urn is a multiset. Also, the tubes
of different colours are represented as a multiset. Multisets form the
proper formalism for multivariate probabilities, see
Section~\ref{PrelimSec} below. Sending an arbitrary set $X$ to the set
of multisets over $X$ has the structure of a monad. Similarly, taking
distributions over a set forms a monad. These monad structures play an
important role in the various ways that the multinomial and
hypergeometric (and P\'olya) operations, as Kleisli maps, can be
composed, see~\cite{Jacobs21b} and~\cite{Jacobs21d}. In this paper the
underlying categorical structure is kept in the background. This is a
deliberate choice, in order not to limit the potential audience. For
instance, in the beginning of Section~\ref{FirstFullDstSec} the
composition steps for Markov models with output are spelled-out
concretely; their abstract categorical form as Kleisli composition is
elaborated in the appendix.

This paper is organised as follows. It starts with a concrete
description of first-full distributions, for all three modes (``0'',
``-1'', ``+1''), in Section~\ref{ExSec}. Subsequently,
Section~\ref{PrelimSec} introduces relevant notation and terminology
for multisets and distributions, and
Section~\ref{MulnomHypgeomPolyaSec} formulates multivariate versions
of the multinomial, hypergeometric and P\'olya distributions, as
distributions on multisets of a fixed size.  For snappy formulation of
the hypergeometric and P\'oly distributions we use binomial
coefficients with multisets instead of numbers, both for ordinary
binomial coefficients and for the multichoose version.  We introduce
suitable generalisations of Vandermonde's formula, for multisets, also
in multichoose form.

The next two sections are devoted to first-full distributions.  They
are defined in Section~\ref{FirstFullDefSec} in a pointwise manner, as
sums over multisets. These probabilities are illustrated in several
bar plots. Next, Section~\ref{FirstFullDstSec} introduces three
probabilistic automata, in the form of Markov models with output,
which are used to show that we actually get three distributions, with
probabilities adding up to one. The heart of the argument is that
composition (iteration) of the steps of these automata preserves
distributions.

Section~\ref{NegativeDrawSec} introduces and illustrates negative
distributions in the urns \& tubes setting. We pay special attention
to the bivariate case, with one tube only, which is the form in which
they occur in the literature. We illustrate how the probabilities add
up to one, and thus yield actual distributions, vie the same
compositional argument, by sketching the relevant Markov models with
output.

These descriptions of the first-full distributions set the scene for
two additional topics. It is known that the hypergeometric and P\'olya
distributions can be obtained via conditioning from binomial and from
negative binomial distributions. Section~\ref{BinomialSec} recalls
these results in the current setting, in uniform descriptions.
Finally, Section~\ref{CorollarySec} concentrates on the bivariate
first-full and negative distributions, with two tubes. It translates
the fact that probabilities in these distributions add up to one into
number-theoretic corollaries. These results seem to be new. It is left
as a challenge to prove them directly.

\section*{Acknowledgments} The urn \& tubes set-up in this
paper was developed without awareness of the problem of points. Thanks
are due to Onno Boxma for pointing out the connection.

\section{Examples of first-full distributions}\label{ExSec}

This section illustrates how first-full distributions come about, in
the three drawing modes. We keep things simple at this stage and use
only two colours, written $R$ for red and $B$ for blue. We assume
length $2$ for the red tube, and length $3$ for the blue tube.  We
briefly describe the draw probabilities in the three modes, in this
illustration.
\begin{itemize}
\item In the multinomial mode ``0'', we assume that there are three
  balls in the urn, one red and two blue. The probability of drawing
  red is thus $\frac{1}{3}$ and for blue it is $\frac{2}{3}$. These
  probabilities remain the same, since drawn balls are returned to the
  urn. This ``0'' case is elaborated in
  Example~\ref{FirstFullDstExMulnom} below.

\item In the hypergeometric ``-1'' mode we assume that the urn
  initially has three red and six blue balls. The initial probability
  of drawing red from this urn is thus $\frac{3}{9} = \frac{1}{3}$. It
  leaves an urn with two red and six blue balls. So the probability of
  drawing a second red ball from the resulting urn is $\frac{2}{8} =
  \frac{1}{4}$. Example~\ref{FirstFullDstExHypgeom} gives the
  first-full details in this mode.

\item In our illustration of the P\'olya ``+1'' mode we assume
  initially just one red and one blue ball in the urn. The probability
  of drawing red is thus initially $\frac{1}{2}$. Upon drawing red,
  not only the drawn red ball, but also an additional red ball, is
  added to the urn, so that it subsequently contains three balls, two
  red and one blue. The probability of drawing red is then
  $\frac{2}{3}$. The resulting dynamics is described in
  Example~\ref{FirstFullDstExPolya} below.
\end{itemize}

\noindent These three modes are elaborated below, in three separate
examples.




\begin{example}
\label{FirstFullDstExMulnom} 
As described above, in the multinomial case we assume a $\frac{1}{3}$
probability for $R$ = red, and $\frac{2}{3}$ for $B$ = blue. What are
the possible draws for getting one tube filled first? We list the
possible draws to fill the red tube (of length 2) first, on the left
below, with corresponding (multinomial) probabilities, and the draws
to fill the blue tube (of length 3) first on the right.
\begin{ceqn}
\[ \begin{array}{rcrcl}
R,R
& \quad &
\frac{1}{3}\cdot\frac{1}{3}
& = &
\frac{1}{9}
\\
B,R,R
& &
\frac{2}{3}\cdot\frac{1}{3}\cdot\frac{1}{3}
& = &
\frac{2}{27}
\\
R,B,R
& &
\frac{1}{3}\cdot\frac{2}{3}\cdot\frac{1}{3}
& = &
\frac{2}{27}
\\
B,B,R,R
& &
\frac{2}{3}\cdot\frac{2}{3}\cdot\frac{1}{3}\cdot\frac{1}{3}
& = &
\frac{4}{81}
\\
B,R,B,R
& &
\frac{2}{3}\cdot\frac{1}{3}\cdot\frac{2}{3}\cdot\frac{1}{3}
& = &
\frac{4}{81}
\\
R,B,B,R
& &
\frac{1}{3}\cdot\frac{2}{3}\cdot\frac{2}{3}\cdot\frac{1}{3}
& = &
\frac{4}{81}
\\[-1em]
& &
\llap{$\underline{\hspace*{5.5em}}$} & &
\\
& & \mbox{total for $R$}
& = &
\frac{11}{27}
\end{array}
\hspace*{8em}
\begin{array}{rcrcl}
B,B,B
& \quad &
\frac{2}{3}\cdot\frac{2}{3}\cdot\frac{2}{3}
& = &
\frac{8}{27}
\\
R,B,B,B
& &
\frac{1}{3}\cdot\frac{2}{3}\cdot\frac{2}{3}\cdot\frac{2}{3}
& = &
\frac{8}{81}
\\
B,R,B,B
& &
\frac{2}{3}\cdot\frac{1}{3}\cdot\frac{2}{3}\cdot\frac{2}{3}
& = &
\frac{8}{81}
\\
B,B,R,B
& &
\frac{2}{3}\cdot\frac{2}{3}\cdot\frac{1}{3}\cdot\frac{2}{3}
& = &
\frac{8}{81}
\\[-1em]
& &
\llap{$\underline{\hspace*{5em}}$} & &
\\
& & \mbox{total for $B$}
& = &
\frac{16}{27}.
\end{array} \]
\end{ceqn}

\noindent We see on the left that the two required $R$-draws can be
mixed with at most two $B$-draws, but the last draw must of course be
$R$, in order to completely fill the red tube first. Similarly, on the
right, the three required $B$-draws can be mixed with at most one
$R$-draw.

The probability that the blue tube is full first is the highest one.
The blue tube is longer than the red one (3 versus 2), but the
probability of drawing blue is higher (2 versus 1).\\
\end{example}

\begin{example}
\label{FirstFullDstExHypgeom} We turn to the hypergeometric mode
and start with an urn with three red and six blue balls. Now each
drawn ball is removed from the urn, which affects subsequent
probabilities. This gives different probabilities for the same draws
as in the previous example.
\begin{ceqn}
\[ \begin{array}{rcrcl}
R,R
& \quad &
\frac{3}{9}\cdot\frac{2}{8}
& = &
\frac{1}{12}
\\
B,R,R
& &
\frac{6}{9}\cdot\frac{3}{8}\cdot\frac{2}{7}
& = &
\frac{1}{14}
\\
R,B,R
& &
\frac{3}{9}\cdot\frac{6}{8}\cdot\frac{2}{7}
& = &
\frac{1}{14}
\\
B,B,R,R
& &
\frac{6}{9}\cdot\frac{5}{8}\cdot\frac{3}{7}\cdot\frac{2}{6}
& = &
\frac{5}{84}
\\
B,R,B,R
& &
\frac{6}{9}\cdot\frac{3}{8}\cdot\frac{5}{7}\cdot\frac{2}{6}
& = &
\frac{5}{84}
\\
R,B,B,R
& &
\frac{3}{9}\cdot\frac{6}{8}\cdot\frac{5}{7}\cdot\frac{2}{6}
& = &
\frac{5}{84}
\\[-1em]
& &
\llap{$\underline{\hspace*{5.5em}}$} & &
\\
& & \mbox{total for $R$}
& = &
\frac{17}{42}
\end{array}
\hspace*{8em}
\begin{array}{rcrcl}
B,B,B
& \quad &
\frac{6}{9}\cdot\frac{5}{8}\cdot\frac{4}{7}
& = &
\frac{5}{21}
\\
R,B,B,B
& &
\frac{3}{9}\cdot\frac{6}{8}\cdot\frac{5}{7}\cdot\frac{4}{6}
& = &
\frac{5}{42}
\\
B,R,B,B
& &
\frac{6}{9}\cdot\frac{3}{8}\cdot\frac{5}{7}\cdot\frac{4}{6}
& = &
\frac{5}{42}
\\
B,B,R,B
& &
\frac{6}{9}\cdot\frac{5}{8}\cdot\frac{3}{7}\cdot\frac{4}{6}
& = &
\frac{5}{42}
\\[-1em]
& &
\llap{$\underline{\hspace*{5em}}$} & &
\\
& & \mbox{total for $B$}
& = &
\frac{25}{42}.
\end{array} \]
\end{ceqn}

\noindent In this hypergeometric mode the two first-full probabilities
are different from the ones in Example~\ref{FirstFullDstExMulnom}, but
they still add up to one. Again, blue `wins'.\\
\end{example}

\begin{example}
\label{FirstFullDstExPolya} Finally we consider first-full in
P\'olya mode, with (initial) urn containing one red and one blue
ball. The probabilities for the various draws are then as follows.

\begin{ceqn}
\[ \begin{array}{rcrcl}
R,R
& \quad &
\frac{1}{2}\cdot\frac{2}{3}
& = &
\frac{1}{3}
\\
B,R,R
& &
\frac{1}{2}\cdot\frac{1}{3}\cdot\frac{2}{4}
& = &
\frac{1}{12}
\\
R,B,R
& &
\frac{1}{2}\cdot\frac{1}{3}\cdot\frac{2}{4}
& = &
\frac{1}{12}
\\
B,B,R,R
& &
\frac{1}{2}\cdot\frac{2}{3}\cdot\frac{1}{4}\cdot\frac{2}{5}
& = &
\frac{1}{30}
\\
B,R,B,R
& &
\frac{1}{2}\cdot\frac{1}{3}\cdot\frac{2}{4}\cdot\frac{2}{5}
& = &
\frac{1}{30}
\\
R,B,B,R
& &
\frac{1}{2}\cdot\frac{1}{3}\cdot\frac{2}{4}\cdot\frac{2}{5}
& = &
\frac{1}{30}
\\[-1em]
& &
\llap{$\underline{\hspace*{5.5em}}$} & &
\\
& & \mbox{total for $R$}
& = &
\frac{3}{5}
\end{array}
\hspace*{8em}
\begin{array}{rcrcl}
B,B,B
& \quad &
\frac{1}{2}\cdot\frac{2}{3}\cdot\frac{3}{4}
& = &
\frac{1}{4}
\\
R,B,B,B
& &
\frac{1}{2}\cdot\frac{1}{3}\cdot\frac{2}{4}\cdot\frac{3}{5}
& = &
\frac{1}{20}
\\
B,R,B,B
& &
\frac{1}{2}\cdot\frac{1}{3}\cdot\frac{2}{4}\cdot\frac{3}{5}
& = &
\frac{1}{20}
\\
B,B,R,B
& &
\frac{1}{2}\cdot\frac{2}{3}\cdot\frac{1}{4}\cdot\frac{3}{5}
& = &
\frac{1}{20}
\\[-1em]
& &
\llap{$\underline{\hspace*{5em}}$} & &
\\
& & \mbox{total for $B$}
& = &
\frac{2}{5}.
\end{array} \]
\end{ceqn}

\noindent We obtain a third first-full distribution, now with higher
probability for red.
\end{example}

\section{Preliminaries, on multisets and distributions}\label{PrelimSec}

We briefly describe the notation and terminology for multisets
and distributions, in two separate subsections.

\subsection{Multisets}\label{MltSubsec}

As mentioned in the introduction, a multiset (or bag) is a finite
`subset' in which elements may occur multiple times. We use a `ket'
notation $\ket{-}$ borrowed from quantum theory, as convenient way of
writing such multisets. For instance, the initial urn with tree red
and six blue balls in Example~\ref{FirstFullDstExHypgeom} forms a
multiset $3\ket{R} + 6\ket{B}$. And the three tubes
in~\eqref{UrnTubesPic} form a multiset $3\ket{R} + 6\ket{B} +
5\ket{G}$. In general, a multiset over a set $X$ is a finite formal
combination $\sum_{i}n_{i}\ket{x_i}$ with $n_{i}\in\NNO$ and $x_{i}
\in X$. Alternatively, a multiset is a function $\varphi \colon X
\rightarrow \NNO$ with finite support $\supp(\varphi) =
\setin{x}{X}{\varphi(x) > 0}$. The number $\varphi(x)\in\NNO$ tells
how many times the element $x$ occurs in the multiset $\varphi$. We
freely switch between the formal sum and the function notation.

We shall write $\Mlt(X)$ for the set of multisets over $X$. Notice
that each multiset is finite, in our description, but the underlying
set $X$ itself need not be finite. Via pointwise addition, multisets
form a commutative monoid, and in fact, $\Mlt(X)$ is the free
commutative monoid on $X$, via the unit map $\eta\colon X \rightarrow
\Mlt(X)$, given by $\eta(x) = 1\ket{x}$. We shall write
$\zero\in\Mlt(X)$ for the empty multiset, with $\zero(x) = 0$ for all
$x\in X$.



We associate several numbers with a multiset.

\begin{definition}
\label{MltNumberDef}
For a multiset $\varphi\in\Mlt(X)$, write:
\begin{enumerate}
\item $\|\varphi\| \coloneqq \sum_{x}\varphi(x)$ for the size of
  $\varphi$, taking multiplicities into account;

\item $\facto{\varphi} \coloneqq \prod_{x} \varphi(x)!$ for the
  multiset factorial;

\item $\coefm{\varphi} \coloneqq
  \displaystyle\frac{\|\varphi\|!}{\facto{\varphi}}$ for the
  multinomial coefficient.
\end{enumerate}
\end{definition}

We are often interested in multisets of a particular size $K\in\NNO$,
so we define a subset:
\[ \begin{array}{rcl}
\Mlt[K](X)
& \coloneqq &
\setin{\varphi}{\Mlt(X)}{\|\varphi\|=K}.
\end{array} \]

\noindent This $\natMlt[K]$ is a functor, but not a monad.

Sequences can be turned into multisets, via `accumulator' functions
$\acc \colon X^{K} \rightarrow \Mlt[K](X)$, given by $\acc(x_{1},
\ldots, x_{K}) \coloneqq 1\ket{x_1} + \cdots + 1\ket{x_K}$.  Thus, for
instance, $\acc(a,a,b,a) = 3\ket{a} + 1\ket{b}$. The multinomial
coefficient $\coefm{\varphi}$ is used in this paper in the following
two ways.

\begin{fact}
\label{CoefficientFact}
\begin{enumerate}
\item \label{CoefficientFactCount} For a multiset $\varphi$ there are
  $\coefm{\varphi}$ lists that accumulate to $\varphi$, that is,
  $\big|\,\acc^{-1}(\varphi)\,\big| = \coefm{\varphi}$.

\item \label{CoefficientFactMulnomThm} For real numbers $a_{1},
  \ldots, a_{n}$, the multinomial theorem says:
\[ \begin{array}[b]{rcl}
\big(a_{1} + \cdots + a_{n}\big)^{K}
& = &
\displaystyle\sum_{\varphi\in\Mlt[K](\{1,\ldots n\})} \coefm{\varphi} \cdot
   \textstyle{\displaystyle\prod}_{i} \, a_{i}^{\varphi(i)}.
\end{array} \eqno{\QEDbox} \]
\end{enumerate}
\end{fact}

Multisets can be ordered pointwise, giving rise to some subtle distinctions.

\begin{definition}
\label{MltOrderDef}
Let $\varphi,\psi\in\Mlt(X)$ be given. We write:
\begin{enumerate}
\item $\varphi\leq\psi$ if $\varphi(x)\leq\psi(x)$ for all $x\in X$;
  in that case we define the multiset difference $\psi-\varphi$ via
  pointwise subtraction, as: $(\psi-\varphi)(x) = \psi(x) -
  \varphi(x)$;

\item $\varphi\leq_{K}\psi$ if $\|\varphi\| = K$ and $\varphi\leq\psi$;

\item $\varphi < \psi$ if $\varphi\leq\psi$ but $\varphi\neq\psi$;

\item $\varphi \prec \psi$ if $\varphi(x) < \psi(x)$ for all $x\in X$.
\end{enumerate}
\end{definition}

\noindent The relation $\prec$ will be called \emph{fully below}. It
is different from $<$, \textit{e.g.} in:
\[ \begin{array}{rclcrcl}
2\ket{a} + 3\ket{b}
& < &
3\ket{a} + 3\ket{b}
& \mbox{\qquad and \qquad} &
2\ket{a} + 2\ket{b}
& \prec &
3\ket{a} + 3\ket{b}.
\end{array} \]

For multisets $\varphi,\psi\in\Mlt(X)$ with $\varphi\leq\psi$,
we define the \emph{multiset binomial} as:
\begin{ceqn}
\begin{equation}
\label{MltBinomEqn}
\begin{array}{rcl}
\displaystyle\binom{\psi}{\varphi}
\hspace*{\arraycolsep}\coloneqq\hspace*{\arraycolsep}
\displaystyle\frac{\facto{\psi}}{\facto{\varphi}\cdot \facto{(\psi-\varphi)}}
& = &
\displaystyle\frac{\prod_{x}\psi(x)!}{\big(\prod_{x}\varphi(x)!\big)\cdot
   \big(\prod_{x}(\psi(x)-\varphi(x)\big)!)}
\\[+1em]
& = &
\displaystyle\prod_{x\in X}\, \frac{\psi(x)!}
   {\varphi(x)!\cdot(\psi(x)-\varphi(x))!}
\\[+1em]
& = &
\displaystyle\prod_{x\in X}\, \binom{\psi(x)}{\varphi(x)}.
\end{array}
\end{equation}
\end{ceqn}

\noindent Intuitively, this is the number of ways $\varphi$ can sit
inside $\psi$.

The next result guarantees that hypergeometric draws form a
distribution. It is well known, but not in this form given below, with
binomials for multisets.

\begin{proposition}
\label{HypgeomCountProp}
For a multiset $\psi\in\Mlt(X)$ of size $L = \|\psi\|$ and for
a number $K\leq L$,
\[ \begin{array}{rcl}
\displaystyle\sum_{\varphi\leq_{K}\psi}\, \binom{\psi}{\varphi}
& = &
\displaystyle\binom{L}{K}.
\end{array} \]
\end{proposition}

The binary case, when the set $X$ has two elements, is known as
Vandermonde's formula, see~\eqref{VandermondeBinaryEqn} below. The
above generalisation can be obtained from it by induction on the
number of elements in the support of $\psi$. For completeness, we
include the proof. It uses Pascal's rule, which says:
\begin{equation}
\label{PascalRuleEqn}
\begin{array}{rcl}
\displaystyle\binom{n}{m} + \binom{n}{m+1}
& = &
\displaystyle\binom{n+1}{m}.
\end{array}
\end{equation}

\begin{myproof}
We use induction on the number of elements in the support
$\supp(\psi)$ of the multiset $\psi$. We go through some initial
values explicitly. If the number of elements is $0$, then $\psi =
\zero$ and so $L = 0 = K$ and $\varphi\leq_{K} \psi$ means $\varphi =
\zero$, so that the result holds. Similarly, if $\supp(\psi)$ is a
singleton, say $\{x\}$, then $L = \psi(x)$.  For $K \leq L$ and
$\varphi\leq_{K}\psi$ we get $\supp(\varphi) = \{x\}$ and $K =
\varphi(x)$. The result then obviously holds.

The case where $\supp(\psi) = \{x,y\}$ captures the ordinary form of
Vandermonde's formula. We reformulate it for numbers $B,G\in\NNO$ and
$K\leq B+G$. Then:
\begin{equation}
\label{VandermondeBinaryEqn}
\begin{array}{rcl}
\displaystyle\binom{B+G}{K}
& = &
\displaystyle\sum_{b\leq B, \, g\leq G, \, b+g = K}\, 
   \binom{B}{b}\cdot\binom{G}{g}.
\end{array}
\end{equation}

\noindent Intuitively: if you select $K$ children out of $B$ boys and
$G$ girls, the number of options is given by the sum over the options
for $b\leq B$ boys times the options for $g\leq G$ girls, with $b+g =
K$. The equation~\eqref{VandermondeBinaryEqn} is standard, so a proof
(\textit{e.g.}~by induction on $G$) is skipped.

\auxproof{
When $G=0$ both sides amount to $\binom{B}{K}$ so we quickly proceed
to the induction step. The case $K=0$ is trivial, so we may assume $K
> 0$.
\[ \begin{array}{rcl}
\lefteqn{\sum_{b\leq B, \, g\leq G+1, \, b+g = K}\, 
   \binom{B}{b}\cdot\binom{G\!+\!1}{g}}
\\[+1.4em]
& = &
\displaystyle\binom{B}{K}\cdot\binom{G\!+\!1}{0} + 
   \binom{B}{K\!-\!1}\cdot\binom{G\!+\!1}{1} + \cdots +
   \binom{B}{K\!-\!G}\cdot\binom{G\!+\!1}{G} + 
   \binom{B}{K\!-\!G\!-\!1}\cdot\binom{G\!+\!1}{G\!+\!1}
\\[+1em]
& \smash{\stackrel{\eqref{PascalRuleEqn}}{=}} &
\displaystyle\binom{B}{K}\cdot\binom{G}{0} + 
   \binom{B}{K\!-\!1}\cdot\binom{G}{1} + 
   \binom{B}{K\!-\!1}\cdot\binom{G}{0} 
\\[+0.8em]
& & \qquad +\, \cdots +\, \displaystyle
   \binom{B}{K\!-\!G}\cdot\binom{G}{G} + 
   \binom{B}{K\!-\!G}\cdot\binom{G}{G\!-\!1} + 
   \binom{B}{K\!-\!G\!-\!1}\cdot\binom{G}{G}
\\[+1em]
& = &
\displaystyle\sum_{b\leq B, \, g\leq G, \, b+g = K}\, 
   \binom{B}{b}\cdot\binom{G}{g} +
   \sum_{b\leq B, \, g\leq G, \, b+g = K-1}\,
   \binom{B}{b}\cdot\binom{G}{g}
\\[+1.4em]
& \smash{\stackrel{\text{(IH)}}{=}} &
\displaystyle\binom{B\!+\!G}{K} + \binom{B\!+\!G}{K\!-\!1}
\\[+1em]
& \smash{\stackrel{\eqref{PascalRuleEqn}}{=}} &
\displaystyle\binom{B\!+\!G\!+\!1}{K}.
\end{array} \]
}
For the induction step, let $\supp(\psi) = \{x_{1}, \ldots, x_{n},
y\}$, for $n\geq 2$.  Writing $\ell = \psi(y)$, $L' = L - \ell$ and
$\psi' = \psi - \ell\ket{y}\in\natMlt[L'](X)$ gives:
\[ \begin{array}[b]{rcl}
\displaystyle\sum_{\varphi \leq_{K} \psi} \, \binom{\psi}{\varphi}
\hspace*{\arraycolsep}=\hspace*{\arraycolsep}
\displaystyle\sum_{\varphi \leq_{K} \psi} \, 
   {\displaystyle\prod}_{x}\, \binom{\psi(x)}{\varphi(x)}
& = &
\displaystyle\sum_{n\leq \ell} \, 
   \sum_{\varphi \leq_{K-n} \psi'} \, \binom{\ell}{n} \cdot 
   {\displaystyle\prod}_{i}\, \displaystyle\binom{\psi'(x_{i})}{\varphi(x_{i})}
\\[+1.4em]
& \smash{\stackrel{\text{(IH)}}{=}} &
\displaystyle\sum_{n\leq \ell, \, K-n\leq L-\ell} \, 
   \binom{\ell}{n} \cdot \binom{L\!-\!\ell}{K\!-\!n}
\hspace*{\arraycolsep}\smash{\stackrel{\eqref{VandermondeBinaryEqn}}{=}}\hspace*{\arraycolsep}
\binom{L}{K}.
\end{array} \eqno{\QEDbox} \]
\end{myproof}

We recall that for $n>0$ and $m\geq 0$ there is the \emph{multichoose}
coefficient, defined for $n \geq 1$ and $m\geq 0$ as:
\[ \begin{array}{rcccl}
\displaystyle\bibinom{n}{m}
& \coloneqq &
\displaystyle\binom{n+m-1}{m}
& = &
\displaystyle\frac{(n+m-1)!}{m!\cdot (n-1)!}.
\end{array} \]

\noindent Interestingly, where $\binom{n}{m}$ is the number of
\emph{subsets} of size $m$ of an $n$-element set, $\bibinom{n}{m}$ is
the number of \emph{multisets} of size $m$ over an $n$-element set.
It is easy to see that:
\begin{equation}
\label{BibinomSuccEqn}
\begin{array}{rcl}
\displaystyle\bibinom{n+1}{m+1}
& = &
\displaystyle\bibinom{n+1}{m} + \bibinom{n}{m+1}.
\end{array}
\end{equation}

\noindent We extend multichoose from numbers to multisets, in line
with~\eqref{MltBinomEqn}:
\[ \begin{array}{rcl}
\displaystyle\bibinom{\psi}{\varphi}
& \coloneqq &
\displaystyle\prod_{x\in\supp(\psi)} \, \bibinom{\psi(x)}{\varphi(x)}.
\end{array} \]

\noindent There is the following multichoose analogue of
Proposition~\ref{HypgeomCountProp}.

\begin{proposition}
\label{PolyaCountProp}
For a multiset $\psi\in\Mlt(X)$ of size $L = \|\psi\| > 0$ and for
any number $K\geq 0$,
\[ \begin{array}{rcl}
\displaystyle\sum_{\varphi\in\Mlt[K](\supp(\psi))}\, \bibinom{\psi}{\varphi}
& = &
\displaystyle\bibinom{L}{K}.
\end{array} \]
\end{proposition}

\begin{myproof}
We start with a double-bracket analogue
of~\eqref{VandermondeBinaryEqn}.  Fix $B\geq 1$ and $G \geq 1$. For
all $K$ one has:
\begin{equation}
\label{VandermondeBibinaryEqn}
\begin{array}{rcl}
\displaystyle\bibinom{B+G}{K}
& = &
\displaystyle\sum_{0\leq k\leq K}\, \bibinom{B}{k} \cdot \bibinom{G}{K-k}.
\end{array}
\end{equation}

\noindent We first prove this equation by induction on $B\geq $1. In
both the base case $B=1$ and the induction step we shall use induction
on $K$. We shall try to keep the structure clear by using nested
bullets.
\begin{itemize}
\item We first prove Equation~\eqref{VandermondeBibinaryEqn} for
  $B=1$, by induction on $K$.
\begin{itemize}
\item When $K=0$ both sides in~\eqref{VandermondeBibinaryEqn} are
  equal to $1$.

\item Assume Equation~\eqref{VandermondeBibinaryEqn} holds for $K$
  (and $B=1$). 
\[ \!\!\begin{array}{rcl}
\displaystyle\sum_{0\leq k\leq K+1} 
   \bibinom{1}{k} \cdot \bibinom{G}{(K\!+\!1)\!-\!k}
& = &
\displaystyle\sum_{0\leq k\leq K+1} \bibinom{G}{K\!-\!(k\!-\!1)}
\\[+1.4em]
& = &
\displaystyle\bibinom{G}{K\!+\!1} + \sum_{0\leq \ell\leq K} 
    \bibinom{1}{\ell} \cdot \bibinom{G}{K\!-\!\ell}
\\[+1.4em]
& \smash{\stackrel{\text{(IH)}}{=}} &
\displaystyle\bibinom{G}{K\!+\!1} + \bibinom{G\!+\!1}{K}
\\[+1.2em]
& \smash{\stackrel{\eqref{BibinomSuccEqn}}{=}} &
\displaystyle\bibinom{G\!+\!1}{K\!+\!1}.
\end{array} \]
\end{itemize}

\item Now assume Equation~\eqref{VandermondeBibinaryEqn} holds for $B$
  (for all $G,K$). In order to show that it then also holds for $B+1$
  we use induction on $K$.
\begin{itemize}
\item When $K=0$ both sides in~\eqref{VandermondeBibinaryEqn} are
  equal to $1$.

\item Now assume that Equation~\eqref{VandermondeBibinaryEqn} holds
  for $K$, and for $B$. Then:

 \allowdisplaybreaks{
 \begin{alignat}{2}
& \notag \lefteqn{\sum_{0\leq k\leq K+1} 
   \bibinom{B\!+\!1}{k} \cdot \bibinom{G}{(K\!+\!1)\!-\!k}}
\\
& \notag = 
\displaystyle\bibinom{G}{K\!+\!1} + \sum_{0\leq k\leq K} 
   \bibinom{B\!+\!1}{k\!+\!1} \cdot \bibinom{G}{K\!-\!k}
\\
& \notag \smash{\stackrel{\eqref{BibinomSuccEqn}}{=}} 
\displaystyle\bibinom{G}{K\!+\!1} + \sum_{0\leq k\leq K}
   \left[\,\bibinom{B}{k\!+\!1} + 
   \bibinom{B\!+\!1}{k}\,\right] \cdot \bibinom{G}{K\!-\!k}
\\
& \notag = 
\displaystyle\bibinom{G}{K\!+\!1} + \sum_{0\leq k\leq K} 
   \bibinom{B}{k\!+\!1} \cdot \bibinom{G}{K\!-\!k} \;+
   \sum_{0\leq k\leq K} \bibinom{B\!+\!1}{k} \cdot \bibinom{G}{K\!-\!k}
\\
& \notag \smash{\stackrel{\text{(IH, $K$)}}{=}} 
\displaystyle\sum_{0\leq k\leq K+1} 
   \bibinom{B}{k} \cdot \bibinom{G}{(K\!+\!1)\!-\!k} \;+\; 
   \bibinom{(B\!+\!1)\!+\!G}{K}
\\
& \notag \smash{\stackrel{\text{(IH, $B$)}}{=}} 
\displaystyle\bibinom{B\!+\!G}{K\!+\!1} \,+\, \bibinom{(B\!+\!1)\!+\!G}{K}
\\
& \notag \smash{\stackrel{\eqref{BibinomSuccEqn}}{=}} 
\displaystyle\bibinom{(B\!+\!1)\!+\!G}{K\!+\!1}.
\end{alignat} 
}

\end{itemize}
\end{itemize}

\noindent This completes the proof of~\eqref{VandermondeBibinaryEqn}.
We proceed with the equation in the proposition, via induction on the
number of elements in the support of $\psi$.  By assumption the
support cannot be empty, so the induction starts when the support is a
singleton, say $\supp(\psi) = \{x\}$. But then $\psi(x) = \|\psi\| =
L$ and $\varphi(x) = \|\varphi\| = K$, so the result obviously holds.

Now let $\supp(\psi) = S \cup \{y\}$ where $y\not\in S$ and $S$ is not
empty. Write:
\[ \begin{array}{rclcrccclcrclcrcccl}
L
& = &
\|\psi\|
& \quad &
\ell
& = &
\psi(y)
& > &
0
& \quad &
\psi'
& = &
\psi - \ell\ket{y}
& \quad &
L'
& = &
L - \ell 
& > &
0.
\end{array} \]

\noindent By construction $S = \supp(\psi')$ and $L' = \|\psi'\|$. Now:
\[ \begin{array}[b]{rcl}
\displaystyle\sum_{\varphi\in\natMlt[K](S\cup \{y\})} \, 
   \bibinom{\psi}{\varphi}
& = &
\displaystyle\sum_{\varphi\in\natMlt[K](S\cup \{y\})} \, 
   \prod_{x\in S\cup\{y\}}\, \bibinom{\psi(x)}{\varphi(x)}
\\[1.2em]
& = &
\displaystyle\sum_{0\leq k\leq K} \, \sum_{\varphi\in\natMlt[K-k](S)} \, 
   \bibinom{\psi(y)}{k} \cdot 
   \prod_{x\in S}\, \bibinom{\psi(x)}{\varphi(x)}
\\[+1.4em]
& = &
\displaystyle\sum_{0\leq k\leq K} \, \bibinom{\ell}{k}
   \cdot \displaystyle\sum_{\varphi\in\natMlt[K-k](S)} \, 
   \bibinom{\psi'}{\varphi}
\\[+1.4em]
& \smash{\stackrel{\text{(IH)}}{=}} &
\displaystyle\sum_{0\leq k\leq K} \, 
  \bibinom{\ell}{k} \cdot \bibinom{L'}{K\!-\!k}
\hspace*{\arraycolsep}\smash{\stackrel{\eqref{VandermondeBibinaryEqn}}{=}}\hspace*{\arraycolsep}
\displaystyle\bibinom{\ell\!+\!L'}{K}
\hspace*{\arraycolsep}=\hspace*{\arraycolsep}
\displaystyle\bibinom{L}{K}.
\end{array} \eqno{\QEDbox} \]
\end{myproof}

\subsection{Probability distributions}\label{DstSubsec}

In this paper we concentrate on finite, discrete probability
distributions. Such a distribution, over a set $X$, is a finite formal
convex combination $\sum_{i}r_{i}\ket{x_i}$ with $r_{i}\in[0,1]$
satisfying $\sum_{i}r_{i} = 1$ and with $x_{i}\in X$. Alternatively,
it may be described as a function $\omega \colon X \rightarrow [0,1]$
with finite support $\supp(\omega) \coloneqq \setin{x}{X}{\omega(x) >
  0}$ and with $\sum_{x} \omega(x) = 1$.  We shall write $\Dst(X)$ for
the set of distributions on a set $X$.  This $\Dst$ forms a monad,
just like $\Mlt$.

Distributions on a product set $X \times Y$ are often called
\emph{joint} distributions. One way to obtains such distributions is
to put $\omega\in\Dst(X)$ and $\rho\in\Dst(Y)$ in parallel as
$\omega\otimes\rho\in \Dst(X\times Y)$, where:
\[ \begin{array}{rcl}
\omega\otimes\rho
& = &
\displaystyle\sum_{x\in X, y\in Y} \, \omega(x)\cdot\rho(y)
   \,\bigket{x,y}.
\end{array} \]

\noindent We then write $\omega^{K} = \omega\otimes \cdots
\otimes\omega \in \Dst(X^{K})$, for numbers $K\geq 1$.

Each non-empty multiset can be turned into a distribution, via
normalisation. We shall call this operation \emph{frequentist
  learning}, written as $\flrn$, since it involves learning by
counting. Explicitly:
\[ \begin{array}{rclcrcl}
\flrn\left({\displaystyle\sum}_{i}\, n_{i}\ket{x_i}\right)
& \coloneqq &
{\displaystyle\sum}_{i}\, \displaystyle\frac{n_i}{n}\ket{x_i}
& \mbox{\quad where \quad} &
n
& = &
\sum_{i} n_{i}.
\end{array} \]

\noindent Alternatively, $\flrn(\varphi)(x) = \displaystyle
\frac{\varphi(x)}{\|\varphi\|}$, or simply, $\flrn(\varphi) =
\displaystyle \frac{1}{\|\varphi\|}\cdot\varphi$.


Multisets and distributions as defined above have finite support.  We
shall also need (discrete) distributions with infinite support. Therefore
we define, for an arbitrary set $X$,
\[ \begin{array}{rcl}
\infDst(X)
& \coloneqq &
\set{\omega\colon X \rightarrow [0,1]}{\sum_{x}\omega(x) = 1}.
\end{array} \]

\noindent It can be shown that the support of $\omega\in\infDst(X)$ is
necessarily countable or finite. In practice one typically encounters
$X = \NNO$. For instance, the Poisson distribution can be described as
an element of $\infDst(\NNO)$. Later on we shall describe negative
distributions that also live in $\infDst(\NNO)$.

\section{Multinomial, hypergeometric, and P\'olya 
   distributions}\label{MulnomHypgeomPolyaSec}

This section introduces the multinomial, hypergeometric and P\'olya
distributions, in multivariate form. This is most conveniently done
via (binomial / multichoose) coefficients for multisets, which is
non-standard. The formulations that are used below can be derived in a
compositional manner via iterated drawing of single elements, using a
suitable form of Kleisli composition, see~\cite{Jacobs21a, Jacobs19d}
for details.

\subsection{Multinomial distributions}\label{MulnomSubsec}

Since the urn remains unchanged for multinomial draws, it is most
appropriate to describe it as a distribution $\omega\in\Dst(X)$,
for a set of colours $X$. The multinomial distribution
$\multinomial[K](\omega)$ is a distribution on draws of size $K$, and
thus an element of the set $\Dst\big(\Mlt[K](X)\big)$. Explicitly,
\begin{ceqn}
\begin{equation}
\label{MulnomEqn}
\begin{array}{rcl}
\multinomial[K](\omega)
& \coloneqq &
\displaystyle\sum_{\varphi\in\Mlt[K](X)}\, \coefm{\varphi}\cdot
   \prod_{x\in X} \omega(x)^{\varphi(x)}\,\bigket{\varphi}.
\end{array}
\end{equation}
\end{ceqn}

\noindent The probabilities in this multinomial distribution add up to
one by the Multinomial Theorem, see
Fact~\ref{CoefficientFact}~\eqref{CoefficientFactMulnomThm}. For instance,
\[ \begin{array}{rcl}
\textstyle\multinomial[3]
   \big(\frac{1}{3}\ket{a} + \frac{1}{2}\ket{b} + \frac{1}{6}\ket{c}\big)
& = &
\frac{1}{27}\Bigket{3\ket{a}} + \frac{1}{6}\Bigket{2\ket{a} + 1\ket{b}}
+ \frac{1}{4}\Bigket{1\ket{a} + 2\ket{b}} + \frac{1}{8}\Bigket{3\ket{b}}
\\[+0.5em]
& & \;+\, 
\frac{1}{18}\Bigket{2\ket{a} + 1\ket{c}}
+ \frac{1}{6}\Bigket{1\ket{a} + 1\ket{b} + 1\ket{c}}
+ \frac{1}{8}\Bigket{2\ket{b} + 1\ket{c}}
\\[+0.5em]
& & \;+\, 
\frac{1}{36}\Bigket{1\ket{a} + 2\ket{c}}
+ \frac{1}{24}\Bigket{1\ket{b} + 2\ket{c}}
+ \frac{1}{216}\Bigket{3\ket{c}}
\end{array} \]

\noindent Notice that the right-hand-side is a distribution over
multisets. The multisets are written inside the `big' kets
$\bigket{-}$ using `small' ket $\ket{-}$ for the individual colours
$a,b,c$. The probabilities of these multisets, as draws, are written
before the big kets. This may require some parsing if you see this
notation style for the first time.

The next result expresses multinomial probabilities in terms of
sequences (of drawn balls).

\begin{lemma}
\label{MultinomialLem}
For $\omega\in\Dst(X)$ and $\varphi\in\natMlt[K](X)$ one has:
\[ \begin{array}{rcccl}
\multinomial[K](\omega)(\varphi)
& = &
\displaystyle\sum_{\vec{x}\in\acc^{-1}(\varphi)} \omega^{K}(\vec{x})
& = &
\displaystyle\sum_{\vec{x}\in\acc^{-1}(\varphi)} \textstyle {\displaystyle\prod}_{i}\,
   \omega(x_{i}).
\end{array} \eqno{\QEDbox} \]
\end{lemma}
The bivariate (or binary) form of these multinomial distributions
involves a map $\Dst(2) \rightarrow \Dst\big(\natMlt[K](2)\big)$,
where $2 = \{0,1\}$. Via the isomorphisms $\Dst(2) \cong [0,1]$ and
$\natMlt[K](2) \cong \{0,1,2,\ldots, K\}$ this map is often described
as a \emph{binomial} $\binomial[K] \colon [0,1] \rightarrow
\Dst\big(\{0,1,\ldots,K\}\big)$, given on $r\in[0,1]$ as:
\begin{ceqn}
\begin{equation}
\label{BinomEqn}
\begin{array}{rcl}
\binomial[K](r)
& \coloneqq &
\displaystyle\sum_{0\leq k\leq K}\, \binom{K}{k}\cdot r^{k}\cdot
   (1-r)^{K-k}\,\bigket{k}
\\[+1.3em]
& = &
\displaystyle\sum_{0\leq k\leq K}\, 
   \multinomial[K]\Big(r\ket{0} + (1\!-\!r)\ket{1}\Big)
      \Big(k\ket{0} + (K\!-\!k)\bigket{1}\Big)\,\bigket{k}
\end{array}
\end{equation}
\end{ceqn}

\subsection{Hypergeometric distributions}\label{HypgeomSubsec}

Proposition~\ref{HypgeomCountProp} guarantees that the
probabilities add up to one in the following multivariate definition
of the hypergeometric distribution, again on multisets of size $K$. It
assumes an urn $\upsilon$ of size $L = \|\upsilon\| \geq K$.
\begin{ceqn}
\begin{equation}
\label{HypgeomEqn}
\begin{array}{rcl}
\hypergeometric[K](\upsilon)
& \coloneqq &
\displaystyle\sum_{\varphi\leq_{K}\upsilon}\, 
   \frac{\binom{\upsilon}{\varphi}}{\binom{L}{K}}\,\bigket{\varphi}.
\end{array}
\end{equation}
\end{ceqn}

\noindent For instance,
\[ \begin{array}{rcl}
\hypergeometric[3]\big(4\ket{a} + 6\ket{b}\big)
& = &
\frac{1}{30}\Bigket{3\ket{a}} + \frac{3}{10}\Bigket{2\ket{a} + 1\ket{b}} + 
\frac{1}{2}\Bigket{1\ket{a} + 2\ket{b}} + \frac{1}{6}\Bigket{3\ket{b}}.
\end{array} \]


\begin{lemma}
\label{HypgeomLem}
For an urn $u\in\Mlt(X)$ and a draw $\varphi\leq_{K}\upsilon$,
\[ \begin{array}{rcl}
\hypergeometric[K](\upsilon)(\varphi)
& = &
\displaystyle\sum_{\vec{x} \in \acc^{-1}(\varphi)} \,
   \prod_{0\leq i<K} \flrn\Big(\upsilon-\acc(x_{1}, \ldots, x_{i})\Big)(x_{i+1}).
\end{array} \eqno{\QEDbox} \]
\end{lemma}

\subsection{P\'olya distributions}\label{PolyaSubsec}

The P\'olya distribution can be described in a similar way, using the
multichoose binomial coefficients. It yields a distribution on
multisets of size $K$, for a non-empty urn $\upsilon$, via:
\begin{ceqn}
\begin{equation}
\label{PolyaEqn}
\begin{array}{rcl}
\polya[K](\upsilon)
& \coloneqq &
\displaystyle\sum_{\varphi\in\Mlt[K](\supp(\upsilon))}\, 
   \frac{\big(\!\binom{\upsilon}{\varphi}\!\big)}
   {\big(\!\binom{L}{K}\!\big)}\,\bigket{\varphi}.
\end{array}
\end{equation}
\end{ceqn}

\noindent This is well-defined by Proposition~\ref{PolyaCountProp}. A
subtle point is that the draws $\varphi$ must be restricted to
elements that occur in the urn $\upsilon$. That's achieved by summing over
$\varphi\in\Mlt[K](\supp(\upsilon))$, so that $\supp(\varphi) \subseteq
\supp(\upsilon)$.

The P\'olya distribution is known~\cite{Mahmoud08}, sometimes as
Dirichlet-multinomial. Its formulation in terms of multichoose
multinomial coefficients of multisets~\eqref{PolyaEqn}, in analogy
with the multinomial coefficients of multisets in the hypergeometric
distribution~\eqref{HypgeomEqn}, seems new. Formulation that come
close are~\cite[Eqn.~(A.1)]{SibuyaYS64}
or~\cite[Eqn.~(40.7)]{JohnsonKB97}. The details that this captures the
P\'olya urn --- where a drawn ball is returned together with an
additional copy of the same colour --- are elaborated
in~\cite{Jacobs19d}.

Here is an example of a P\'olya distribution, for the same urn as
above, in the hypergeometric illustration.
\[ \begin{array}{rcl}
\polya[3]\big(4\ket{a} + 6\ket{b}\big)
& = &
\frac{1}{11}\Bigket{3\ket{a}} + \frac{3}{11}\Bigket{2\ket{a} + 1\ket{b}} + 
\frac{21}{55}\Bigket{1\ket{a} + 2\ket{b}} + \frac{14}{55}\Bigket{3\ket{b}}.
\end{array} \]


\begin{lemma}
\label{PolyaLem}
For an urn $\upsilon\in\Mlt(X)$ and a draw $\varphi\in\Mlt[K](X)$ with
$\supp(\varphi) \subseteq \supp(\upsilon)$,
\[ \begin{array}{rcl}
\polya[K](\upsilon)(\varphi)
& = &
\displaystyle\sum_{\vec{x} \in \acc^{-1}(\varphi)} \,
   \prod_{0\leq i<K} \flrn\Big(\upsilon+\acc(x_{1}, \ldots, x_{i})\Big)(x_{i+1}).
\end{array} \eqno{\QEDbox} \]
\end{lemma}

\section{First-full definitions}\label{FirstFullDefSec}

From the illustrations in Section~\ref{ExSec} we can extract the
general formulations for the first-full probabilities. They use the
fully-below relation $\prec$ between multisets from
Definition~\ref{MltOrderDef}, given by $\varphi \prec \psi$ iff
$\varphi(x) < \psi(x)$ for all $x$.  At this stage we only give the
probabilities pointwise. Proving that they add up to one, and thus
form a probability distribution, is achieved later, in
Theorems~\ref{MnclgThm}, \ref{HgclgThm} and~\ref{PoclgThm}.

We write $\one = \sum_{x\in X}1\ket{x}$ for the multiset of singletons
on a finite set $X$ of colours. The tubes in our urns \& tubes setting
are represented as a multiset $\tau\in\Mlt(X)$. We require
$\tau\geq\one$, so that each tube has at least length $1$. Empty tubes
are irrelevant and can be ignored.

The definitions below involve draws $\varphi \prec \tau$, so that none
of the tubes is full yet. For colour $x$ we take those draws $\varphi$
with $\varphi(x) = \tau(x) - 1$, so that only one ball is missing in
tube $x$. The probability of this last ball is included in the three
formulations below, resp.\ as $\omega(x)$, as $\flrn(\psi-\varphi)(x)$
and as $\flrn(\psi+\varphi)(x)$.

\begin{definition}
\label{FirstFullDstDef}
Let $X$ be a finite set of colours with a multiset of tubes
$\tau\geq\one$ over $X$, and let $x\in X$ be an arbitrary element.
\begin{enumerate}
\item \label{FirstFullDstDefMulnom} Let $\omega\in\Dst(X)$ be a
  distribution with full support. The \emph{multinomial first-full}
  probabilities are given via the function $\mnff(\omega,\tau)\colon X
  \rightarrow [0,1]$ determined by:
\[ \begin{array}{rcl}
\mnff(\omega,\tau)(x)
& \coloneqq &
\displaystyle \sum_{\smash{\begin{array}{c} \scriptstyle \\[-1em]
   \scriptstyle \varphi\prec\tau, \\[-0.8em]
   \scriptstyle \varphi(x) = \tau(x)-1 \end{array}}}
     \mulnom(\omega)(\varphi)\cdot\omega(x).
\end{array} \]

\smallskip

\item \label{FirstFullDstDefHypgeom} Let $\upsilon\in\natMlt(X)$ be an urn
  / multiset with $\upsilon\geq\tau$. The \emph{hypergeometric first-full}
  probabilities $\hgff(\upsilon,\tau) \colon X \rightarrow [0,1]$ are
  defined as:
\[ \begin{array}{rcl}
\hgff(\upsilon,\tau)(x)
& \coloneqq &
\displaystyle \sum_{\smash{\begin{array}{c} \scriptstyle \\[-1em]
   \scriptstyle \varphi\prec\tau, \\[-0.8em]
   \scriptstyle \varphi(x) = \tau(x)-1 \end{array}}}
     \hypgeom(\upsilon)(\varphi)\cdot\flrn(\upsilon-\varphi)(x).
\end{array} \]

\smallskip

\item \label{FirstFullDstDefPolya} Let $\upsilon\in\natMlt(X)$ be an urn
  with $\upsilon\geq\one$. The \emph{P\'olya first-full} probabilities
  $\plff(\upsilon,\tau) \colon X \rightarrow [0,1]$ are:
\[ \begin{array}{rcl}
\plff(\upsilon,\tau)(x)
& \coloneqq &
\displaystyle \sum_{\smash{\begin{array}{c} \scriptstyle \\[-1em]
   \scriptstyle \varphi\prec\tau, \\[-0.8em]
   \scriptstyle \varphi(x) = \tau(x)-1 \end{array}}}
     \pol(\upsilon)(\varphi)\cdot\flrn(\upsilon+\varphi)(x).
\end{array} \]

\smallskip

\end{enumerate}
\end{definition}

\begin{figure}[h]
\begin{center}
\begin{tabular}{c||c|c}
$\vcenter{\xymatrix@C-1.5pc@R-1.5pc{
\mbox{}\ar@{-}[dr] & \textbf{tubes}
\\
\textbf{urns} & \mbox{}
}}$
   & $\begin{array}{c}
   \vcenter{\hbox{\includegraphics[scale=0.14]{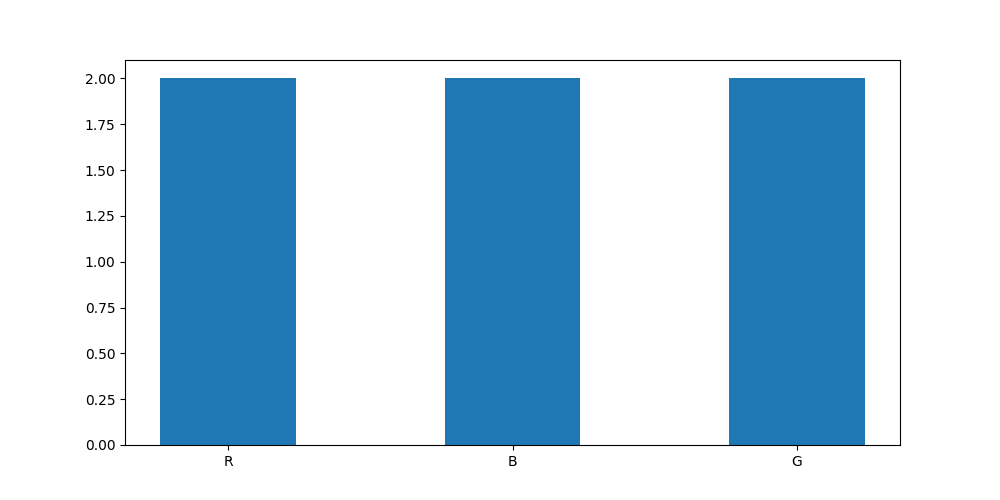}}}
\\[-0.2em]
2\ket{a} + 2\ket{b} + 2\ket{c}
\end{array}$
   & $\begin{array}{c}
   \vcenter{\hbox{\includegraphics[scale=0.14]{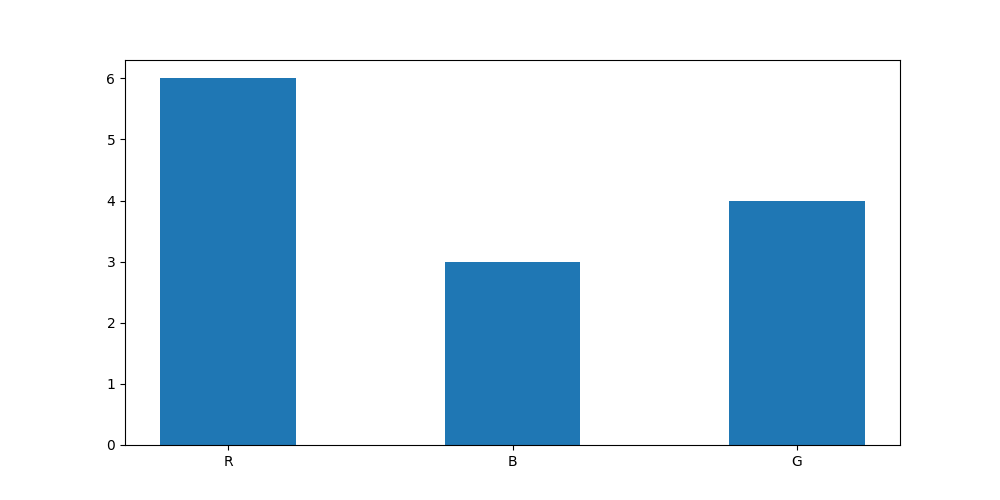}}}
\\[-0.2em]
6\ket{a} + 3\ket{b} + 4\ket{c}
\end{array}$
\\
\hline\hline
\\[-0.5em]
$\hspace*{-1em}\begin{array}{c}
\mbox{multinomial}
\\[-0.4em]
\mbox{first-full}
\\[-0.1em]
\vcenter{\hbox{\includegraphics[scale=0.14]{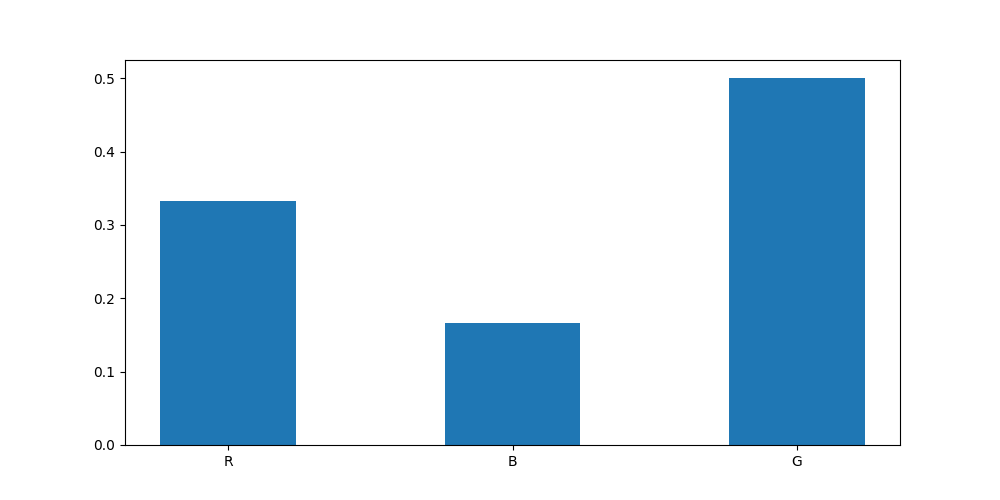}}}
\\[-0.2em]
\frac{1}{3}\ket{a} + \frac{1}{6}\ket{b} + \frac{1}{2}\ket{c}
\end{array}\hspace*{-1em}$
   & $\hspace*{-0.7em}\begin{array}{c}
     \mbox{} \\[-0.2em]
     \vcenter{\hbox{\includegraphics[scale=0.14]{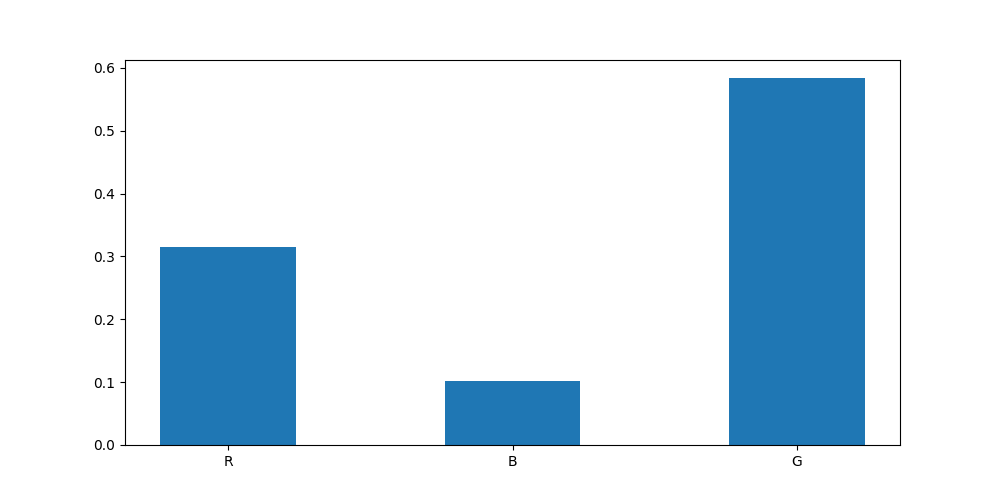}}}
\\
     \frac{17}{54}\ket{a} + \frac{11}{108}\ket{b} + \frac{7}{12}\ket{c}
\end{array}\hspace*{-0.3em}$
   & $\hspace*{-0.3em}\begin{array}{c}
     \mbox{} \\[-0.2em]
     \vcenter{\hbox{\includegraphics[scale=0.14]{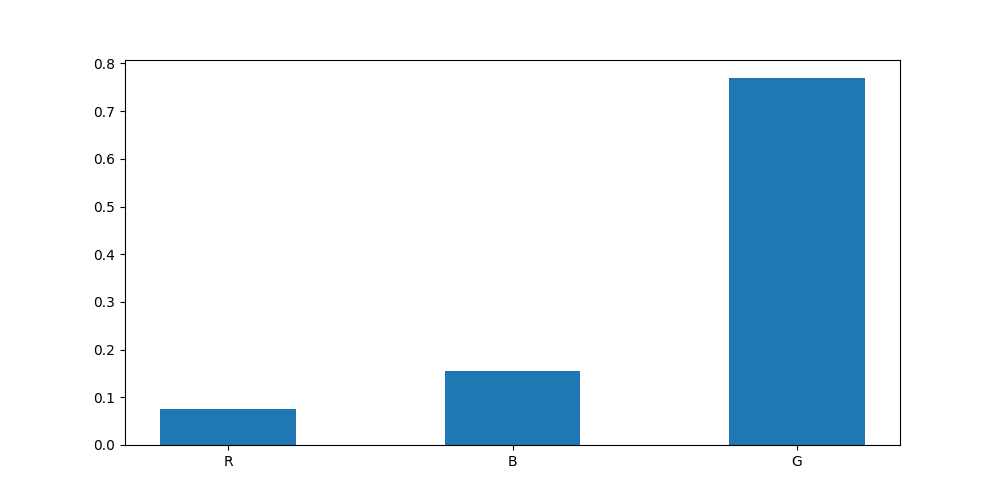}}}
\\
\frac{331}{4374}\ket{a} + \frac{5443}{34992}\ket{b} + \frac{2989}{3888}\ket{c}
\end{array}\hspace*{-1em}$
\\
\\[-0.5em]
\hline
\\[-0.5em]
$\hspace*{-1em}\begin{array}{c}
\mbox{hypergeometric}
\\[-0.4em]
\mbox{first-full}
\\[-0.1em]
\vcenter{\hbox{\includegraphics[scale=0.14]{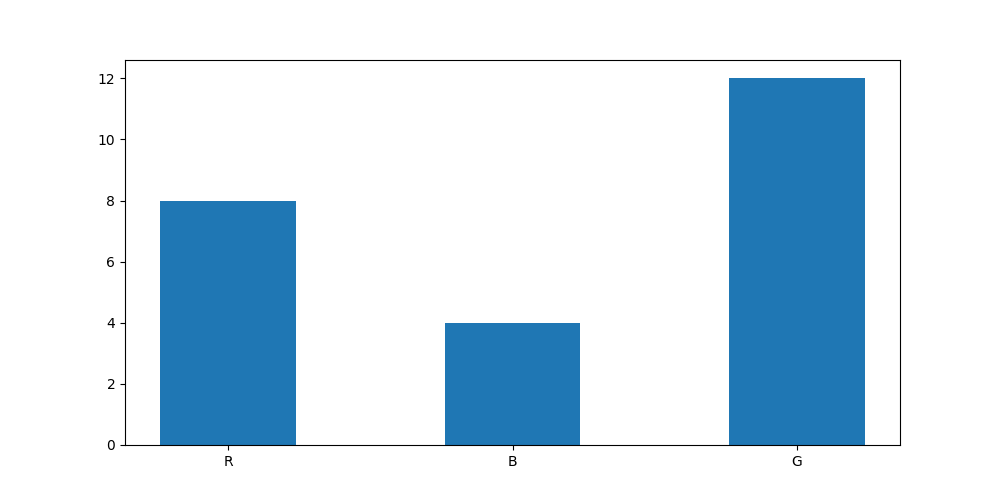}}}
\\[-0.2em]
8\ket{a} + 4\ket{b} + 12\ket{c}
\end{array}\hspace*{-1em}$
   & $\hspace*{-0.5em}\begin{array}{c}
     \mbox{} \\[-0.2em]
     \vcenter{\hbox{\includegraphics[scale=0.14]{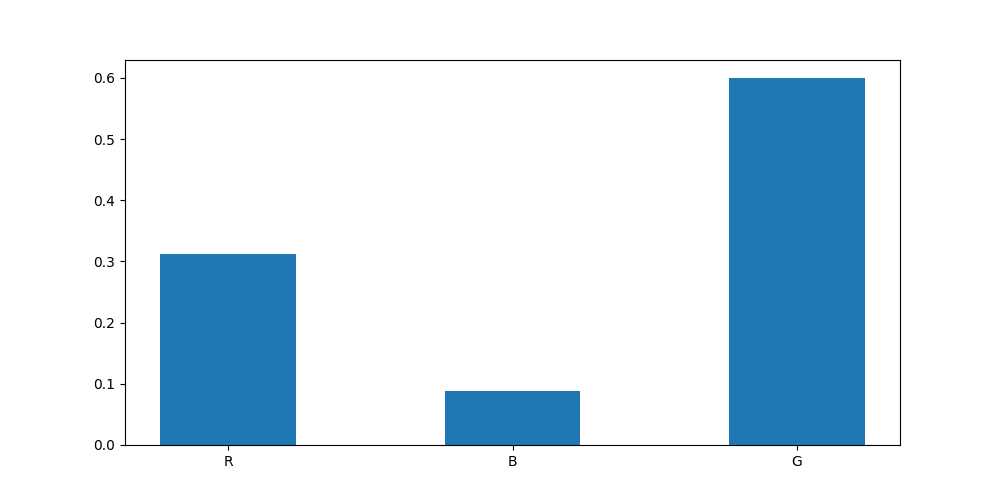}}}
\\
\frac{79}{253}\ket{a} + \frac{313}{3542}\ket{b} + \frac{193}{322}\ket{c}
\end{array}\hspace*{-0.3em}$
   & $\hspace*{-0.3em}\begin{array}{c}
     \mbox{} \\[-0.2em]
     \vcenter{\hbox{\includegraphics[scale=0.14]{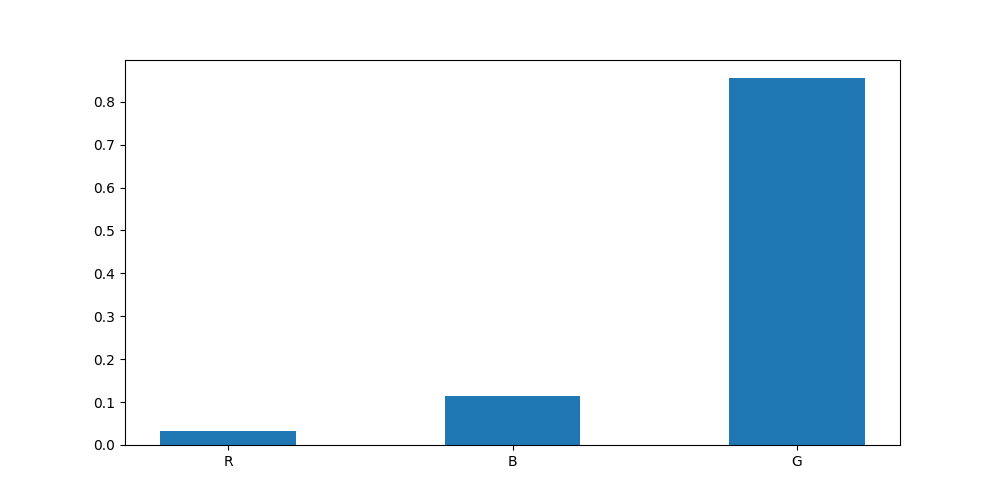}}}
\\
\frac{38843}{1225785}\ket{a} + \frac{1952813}{17160990}\ket{b} 
   + \frac{88875}{104006}\ket{c}
\end{array}\hspace*{-1em}$
\\
\\[-0.5em]
\hline
\\[-0.5em]
$\hspace*{-1em}\begin{array}{c}
\mbox{P\'olya first-full}
\\[-0.1em]
\vcenter{\hbox{\includegraphics[scale=0.14]{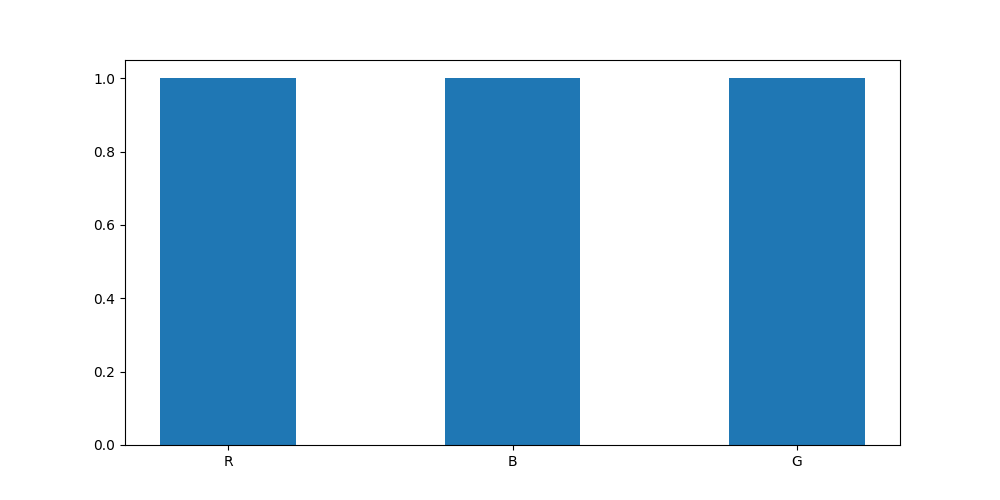}}}
\\[-0.2em]
1\ket{a} + 1\ket{b} + 1\ket{c}
\end{array}\hspace*{-1em}$
   & $\hspace*{-0.7em}\begin{array}{c}
     \mbox{} \\[-0.2em]
     \vcenter{\hbox{\includegraphics[scale=0.14]{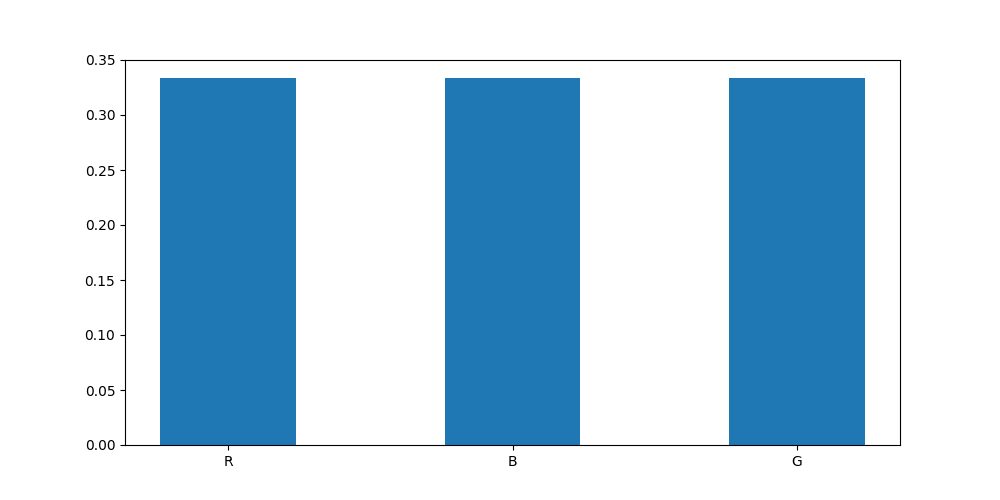}}}
\\
\frac{1}{3}\ket{a} + \frac{1}{3}\ket{b} + \frac{1}{3}\ket{c}
\end{array}\hspace*{-0.3em}$
   & $\hspace*{-0.3em}\begin{array}{c}
     \mbox{} \\[-0.2em]
     \vcenter{\hbox{\includegraphics[scale=0.14]{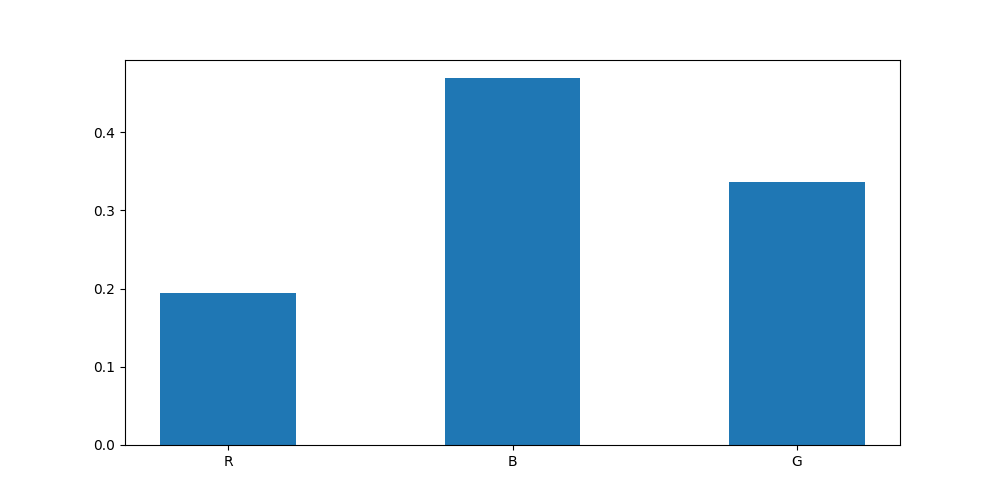}}}
\\
\frac{38}{195}\ket{a} + \frac{128}{273}\ket{b} + \frac{153}{455}\ket{c}
\end{array}\hspace*{-1em}$
\end{tabular}
\end{center}
\caption{First-full distributions on space of colours $\{a,b,c\}$
  arising from the two tube configurations in the top row, and from
  the three urns in the left column. The inner distributions are the
  results, for the multinomial, hypergeometric and P\'olya modes.}

\label{FirstFullDstFig}
\end{figure}

Earlier we have written multinomial, hypergeometric and P\'olya
distributions as $\multinomial[K]$, $\hypergeometric[K]$ and
$\polya[K]$, with explicit parameter $K\in\NNO$ for the size of the
draw. For convenience we have omitted this $K$ in the above
formulations. It may be added as $K=\|\varphi\|$, but that makes the
notation unnecessarily heavy.

In the above definition we require full support of the
urn/distribution $\omega$, for convenience.  We could have been more
relaxed and required only $\supp(\omega) \subseteq \supp(\tau)$ and
$\supp(\upsilon) \subseteq \supp(\tau)$. When these are proper
inclusions, there are tubes that will never receive any balls. Then we
might as well exclude them altogether.

Figure~\ref{FirstFullDstFig} presents illustrations of these different
first-full probabilities, for two different multisets of tubes, at the
top of the second and third column. In the second column the three
tubes have the same length; the corresponding first-full plots then
resemble the urns.  In the third column the tubes differ; the highest
first-full probabilities are determined not only by the lowest tubes,
but also by the highest numbers in the urns. The bar plots are based
on distributions that are computed via the formulations in
Definition~\ref{FirstFullDstDef}.

\subsection{The problem of points}\label{PointsSubsec}

We briefly elaborate the connection between (multinomial) first-full
distributions and the ancient problem of points, as discussed in the
introduction. We do so via an example in~\cite{Ma16}, with two
players, called $A$ and $B$, playing a game that ends when one of the
players has won 4 times. The winner then gets 64 coins. Each time, the
probability of winning for $A$ is $\frac{6}{10}$ and is $\frac{4}{10}$
for $B$.  

A particular game is terminated abruptly at a stage where $A$ has won
$a<4$ times and $B$ has won $b<4$ times. The question that has
occupied Fermat and Pascal is how to fairly divide the stake of 64
coins at such an unfinished stage. Their solution is to look at the
chances of $A$ and $B$ to win, if they were to continue from where the
game was terminated. One then looks at the number of times $4-a$ and
$4-b$ that $A$ and $B$ still need to win. This can be reformulated in
terms of tubes to be filled.

Thus, the distribution capturing the chances for $A$ and $B$ to still
win in this (aborted) situation of the game --- if the game would be
continued --- is a multinomial first-full:
\begin{equation}
\label{ProblemOfPointsEqn}
\begin{array}{rcl}
\rho(a,b)
& \coloneqq &
\mnff\big(\frac{6}{10}\ket{A} + \frac{4}{10}\ket{B},
   \, (4-a)\ket{A} + (4-b)\ket{B}\big). 
\end{array}
\end{equation}


\noindent For instance, $\rho(1,2) = \frac{297}{625}\ket{A} +
\frac{328}{625}\ket{B}$. The division of stakes from the problem of
points can now be formulated in terms of such first-full
distributions.  Figure~3 in~\cite{Ma16}, reconstructed here in
Figure~\ref{OriginalProblemOfPointsFig}, contains, for $a = 1$ and
$b=2$, as fair share for $A$:
\[ \begin{array}{rcccl}
\rho(1,2)(A) \cdot 64
& = &
\frac{297}{625} \cdot 64
& \approx &
30.4128 \quad \mbox{coins.}
\end{array} \]

\noindent In this way all numbers in Figure~3 of~\cite{Ma16} can be
reconstructed, for all numbers $0 \leq a < 4$ and $0 \leq b < 4$,
see Figure~\ref{ReconstructedProblemOfPointsFig}.

\begin{figure}[h]
\begin{center}
{\small\begin{tabular}{c|c||c|c|c|c|c|l}
\hspace*{-1em}\begin{tabular}{c}
{\#}points \\[-0.4em]
remaining \\[-0.4em]
for $B$
\\
\end{tabular}\hspace*{-1em} 
& 
\hspace*{-1em}\begin{tabular}{c}
{\#}points \\[-0.4em]
already \\[-0.4em]
won by $B$
\\
\end{tabular}\hspace*{-1em} 
& 
\multicolumn{4}{c}{} 
\\
\cline{1-7}
\cline{1-7}
0 & 4 & \textcolor{red}{0} & \textcolor{red}{0} & 
   \textcolor{red}{0} & \textcolor{red}{0} & 
\\
\cline{1-7}
1 & 3 & 8.2944 & 13.824 & 23.04 & 38.4 & \textcolor{red}{64}
\\
\cline{1-7}
2 & 2 & 21.56544 & 30.4128 & 41.472 & 53.76 & \textcolor{red}{64}
\\
\cline{1-7}
3 & 1 & 34.83648 & 43.68384 & 52.5312 & 59.904 & \textcolor{red}{64}
\\
\cline{1-7}
4 & 0 & 45.453312 & 52.5312 & 58.42944 & 62.3616 & \textcolor{red}{64}
\\
\hline
\hline
\multicolumn{1}{c}{} & & 0 & 1 & 2 & 3 & 4 & 
   {\hspace*{-0.5em}{\#}points already won by $A$\hspace*{-0.5em}}
\\
\cline{3-8}
\multicolumn{1}{c}{} & & 4 & 3 & 2 & 1 & 0 & 
  {\hspace*{-0.5em}{\#}points remaining for $A$\hspace*{-0.5em}}
\end{tabular}}
\end{center}
\caption{The stakes for player $A$ in the problem of points
  illustration copied from~\cite{Ma16}, see
  Subsection~\ref{PointsSubsec}.}
\label{OriginalProblemOfPointsFig}
\end{figure}

Fermat and Pascal solved the problem of points in binary form, for two
players. Our multivariate formulation of the first-full multinomial
distribution in
Definition~\ref{FirstFullDstDef}~\eqref{FirstFullDstDefMulnom} can be
used when there are multiple (finitely many) players.



\begin{figure}[]
\begin{center}
{\small\begin{tabular}{c|c|c|c|c|c}
$\!\!\!b$
\\
\cline{1-5}
$\!\!\!3$ & $\begin{array}{rcl}
\lefteqn{\textstyle\frac{5184}{625}\ket{A} + \frac{34816}{625}\ket{B}}
\\
& = &
8.2944\ket{A} 
\\[-0.2em]
& & \;\; +\; 55.7056\ket{B}
\end{array}$
&
$\begin{array}{rcl}
\lefteqn{\textstyle\frac{1728}{125}\ket{A} + \frac{6272}{125}\ket{B}}
\\
& = &
13.824\ket{A}
\\[-0.2em]
& & \;\; +\; 50.176\ket{B}
\end{array}$
&
$\begin{array}{rcl}
\lefteqn{\textstyle\frac{576}{25}\ket{A} + \frac{1024}{25}\ket{B}}
\\
& = &
23.04\ket{A}
\\[-0.2em]
& & \;\; +\; 40.96\ket{B}
\end{array}$
&
$\begin{array}{rcl}
\lefteqn{\textstyle\frac{192}{5}\ket{A} + \frac{128}{5}\ket{B}}
\\
& = &
38.4\ket{A}
\\[-0.2em]
& & \;\; +\; 25.6\ket{B}
\end{array}$
\\
\cline{1-5}
$\!\!\!2$ & $\begin{array}{rcl}
\lefteqn{\textstyle\frac{67392}{3125}\ket{A} + \frac{132608}{3125}\ket{B}}
\\
& = &
21.56544\ket{A}
\\[-0.2em]
& & \;\; +\; 42.43456\ket{B}
\end{array}$
&
$\begin{array}{rcl}
\lefteqn{\textstyle\frac{19008}{625}\ket{A} + \frac{20992}{625}\ket{B}}
\\
& = &
30.4128\ket{A}
\\[-0.2em]
& & \;\; +\; 33.5872\ket{B}
\end{array}$
&
$\begin{array}{rcl}
\lefteqn{\textstyle\frac{5184}{125}\ket{A} + \frac{2816}{125}\ket{B}}
\\
& = &
41.472\ket{A}
\\[-0.2em]
& & \;\; +\; 22.528\ket{B}
\end{array}$
&
$\begin{array}{rcl}
\lefteqn{\textstyle\frac{1344}{25}\ket{A} + \frac{256}{25}\ket{B}}
\\
& = &
53.76\ket{A}
\\[-0.2em]
& & \;\; +\; 10.24\ket{B}
\end{array}$
\\
\cline{1-5}
$\!\!\!1$ & $\begin{array}{rcl}
\lefteqn{\textstyle\frac{108864}{3125}\ket{A} + \frac{91136}{3125}\ket{B}}
\\
& = &
34.83648\ket{A}
\\[-0.2em]
& & \;\; +\; 29.16352\ket{B}
\end{array}$
&
$\begin{array}{rcl}
\lefteqn{\textstyle\frac{136512}{3125}\ket{A} + \frac{63488}{3125}\ket{B}}
\\
& = &
43.68384\ket{A}
\\[-0.2em]
& & \;\; +\; 20.31616\ket{B}
\end{array}$
&
$\begin{array}{rcl}
\lefteqn{\textstyle\frac{32832}{625}\ket{A} + \frac{7168}{625}\ket{B}}
\\
& = &
52.5312\ket{A}
\\[-0.2em]
& & \;\; +\; 11.4688\ket{B}
\end{array}$
&
$\begin{array}{rcl}
\lefteqn{\textstyle\frac{7488}{125}\ket{A} + \frac{512}{125}\ket{B}}
\\
& = &
59.904\ket{A}
\\[-0.2em]
& & \;\; +\; 4.096\ket{B}
\end{array}$
\\
\cline{1-5}
$\!\!\!0$ & $\begin{array}{rcl}
\lefteqn{\textstyle\frac{710208}{15625}\ket{A} + \frac{289792}{15625}\ket{B}}
\\
& = &
45.453312\ket{A}
\\[-0.2em]
& & \;\; +\; 18.546688\ket{B} \;\;
\end{array}$
&
$\begin{array}{rcl}
\lefteqn{\textstyle\frac{32832}{625}\ket{A} + \frac{7168}{625}\ket{B}}
\\
& = &
52.5312\ket{A}
\\[-0.2em]
& & \;\; +\; 11.4688\ket{B}
\end{array}$
&
$\begin{array}{rcl}
\lefteqn{\textstyle\frac{182592}{3125}\ket{A} + \frac{17408}{3125}\ket{B}}
\\
& = &
58.42944\ket{A}
\\[-0.2em]
& & \;\; +\; 5.57056\ket{B}\quad
\end{array}$
&
$\begin{array}{rcl}
\lefteqn{\textstyle\frac{38976}{625}\ket{A} + \frac{1024}{625}\ket{B}}
\\
& = &
62.3616\ket{A}
\\[-0.2em]
& & \;\; +\; 1.6384\ket{B}\quad
\end{array}$
\\
\hline
& 0 & 1 & 2 & 3 & $a\!\!\!$
\\
\end{tabular}}
\end{center}
\caption{The stakes for both players $A$ and $B$ in the problem of
  points illustration from~\cite{Ma16}, reconstructed as $64 \cdot
  \rho(a,b)$ via first-full multinomials~\eqref{ProblemOfPointsEqn},
  where $a = \#$points already won by $A$ and $b = \#$points already
  won by $B$. This table corresponds to the numbers in the central
  $4\times 4$ part of Figure~\ref{OriginalProblemOfPointsFig}.}
\label{ReconstructedProblemOfPointsFig}
\end{figure}

\section{First-full yields distributions}\label{FirstFullDstSec}

Our goal in this section is to prove that the probabilities in the
three pointwise first-full formulations in
Definition~\ref{FirstFullDstDef} all add up to one, and thus form
proper probability distributions. The trick is to use the multiset of
tubes as a position in a probabilistic automaton that changes with
every draw-and-drop action. The automaton precisely records the
relevant probabilities and terminates after a finite number of
iterations, with a first-full distribution on colours as result.  This
works because in each composition step distributions are
preserved. Hence, if we start with a distribution, we will also end up
with a distribution, namely a first-full one.

We use Markov models with output (MOOs) as probabilistic
automata. They can be described as functions (coalgebras) of the form:
\begin{equation}
\label{MMODiag}
\xymatrix{
Y\ar[r]^-{c} & \Dst\big(Y+X\big)
}
\end{equation}

\noindent where $Y$ is a set of positions and $X$ is a set of outputs.
The $+$ is a coproduct (disjoint union). We will not use separate
`coprojection' functions for these coproducts; the types of the
elements will make it clear whether they live in the left or right
component of a coproduct $Y+X$. What's important is that a function
$c$ as above can be composed with itself, giving the required
iterations (or transitions) of the automaton: with a successor
position in $Y$ the automaton can continue, and with an output in $X$,
the automaton halts.

A compositional argument underlying the next
iterations~\eqref{IterationEqn} of an automaton~\eqref{MMODiag} is
provided in the appendix. At this stage we use such iterations $c^{n}
\colon Y \rightarrow \Dst(Y+X)$, for $n\in\NNO$, via the concrete
formulations given below, where $y\in Y$ is a start position.
\begin{equation}
\label{IterationEqn}
\begin{array}{rcl}
c^{0}(y)
& \coloneqq &
1\bigket{y}
\\[+0.3em]
c^{n+1}(y)
& \coloneqq &
\displaystyle \sum_{z\in Y+X} \textstyle
   \left(\sum_{y'\in Y} c(y)(y')\cdot c^{n}(y')(z)\right)\bigket{z} \;+\;
   \displaystyle \sum_{x\in X} c(y)(x)\bigket{x}.
\end{array}
\end{equation}

\noindent The first sum defines the transitions and the second sum the
outputs. Notice that these are defined as proper distributions. Hence
via iterated composition only distributions arise. This fact will be
crucial.

In the next three subsections we define three appropriate Markov models
with output~\eqref{MMODiag} with transitions that incorporate the
first-full steps. 
\subsection{Multinomial first-full distributions}\label{MnffCoalgSubsec}
In multinomial mode an urn is represented as a distribution
$\omega\in\Dst(X)$, with full support. We shall write as set of tubes:
\[ \begin{array}{rclcrcl}
\natMlt_{\geq 1}(X)
& \coloneqq &
\setin{\tau}{\natMlt(X)}{\tau \geq \one}
& \quad\mbox{where}\quad &
\one
& = &
\displaystyle\sum_{x\in X}\, 1\ket{x}.
\end{array} \]

\noindent Associated with urn/distribution $\omega$ we define the
following multinomial Markov model with output $\mnclg(\omega)$, with
tubes in $\natMlt_{\geq 1}(X)$ as positions.
\[ \vcenter{\xymatrix@R-2.5pc{
\natMlt_{\geq 1}(X)\ar[r]^-{\mnclg(\omega)} & 
   \Dst\Big(\natMlt_{\geq 1}(X) + X\Big)
\\
\tau\ar@{|->}[r] & 
   \displaystyle\sum_{x,\,\tau(x)>1}\!\omega(x)\bigket{\tau\!-\!1\ket{x}}
   \;+
   \displaystyle\sum_{x,\,\tau(x)=1}\!\omega(x)\bigket{x}.
}} \]

\noindent The aim is to iterate this Markov model with output
$\mnclg(\omega)$, using composition for such models, as described
in~\eqref{IterationEqn}. We illustrate the resulting dynamics by
redoing Example~\ref{FirstFullDstExMulnom}, with state $\omega =
\frac{1}{3}\ket{a} + \frac{2}{3}\ket{b}$ and tubes $\tau = 2\ket{a} +
3\ket{b}$. Then:
\[ \begin{array}{rcl}
\mnclg(\omega)(\tau)
& = &
\frac{1}{3}\Bigket{1\ket{a} + 3\ket{b}} + 
   \frac{2}{3}\Bigket{2\ket{a} + 2\ket{b}}
\\[+0.5em]
\mnclg(\omega)^{2}(\tau)
& = &
\frac{1}{3}\cdot\frac{1}{3}\bigket{a} + 
   \frac{1}{3}\cdot\frac{2}{3}\Bigket{1\ket{a} + 2\ket{b}} +
   \frac{2}{3}\cdot\frac{1}{3}\Bigket{1\ket{a} + 2\ket{b}} +
   \frac{2}{3}\cdot\frac{2}{3}\Bigket{2\ket{a} + 1\ket{b}}
\\[+0.4em]
& = &
\frac{1}{9}\bigket{a} + 
   \frac{4}{9}\Bigket{1\ket{a} + 2\ket{b}} +
   \frac{4}{9}\Bigket{2\ket{a} + 1\ket{b}}
\\[+0.5em]
\mnclg(\omega)^{3}(\tau)
& = &
\frac{1}{9}\bigket{a} + 
   \frac{4}{9}\cdot\frac{1}{3}\bigket{a} +
   \frac{4}{9}\cdot\frac{2}{3}\Bigket{1\ket{a} + 1\ket{b}} +
   \frac{4}{9}\cdot\frac{1}{3}\Bigket{1\ket{a} + 1\ket{b}} +
   \frac{4}{9}\cdot\frac{2}{3}\bigket{b}
\\[+0.4em]
& = &
\frac{7}{27}\bigket{a} + 
   \frac{12}{27}\Bigket{1\ket{a} + 1\ket{b}} +
   \frac{8}{27}\bigket{b}
\\[+0.5em]
\mnclg(\omega)^{4}(\tau)
& = &
\frac{7}{27}\bigket{a} + 
   \frac{12}{27}\cdot\frac{1}{3}\bigket{a} +
   \frac{12}{27}\cdot\frac{2}{3}\bigket{b} +
   \frac{8}{27}\bigket{b}
\\[+0.4em]
& = &
\frac{11}{27}\bigket{a} + \frac{16}{27}\bigket{b}.
\end{array} \]

\noindent This is precisely the outcome that we obtained in
Example~\ref{FirstFullDstExMulnom} by manually checking all options.

We formulate at a more general level what's happening via iteration.

\begin{lemma}
\label{MnclgLem}
Consider the above Markov model with output $\mnclg(\omega)$ for a
distribution $\omega\in\Dst(X)$.
\begin{enumerate}
\item \label{MnclgLemLst} For a multiset $\varphi\in\natMlt_{\geq 1}(X)$,
\[ \begin{array}{rcl}
\mnclg(\omega)^{n}(\tau)(\varphi)
& = &
\displaystyle\sum\,\Big\{\,\omega^{n}(\ell)\,\Big|\,
   \ell=\tuple{x_{1}, \ldots, x_{n}} \in X^{n} 
    \mbox{ with }\varphi = \tau\!-\! 1\ket{x_{1}}\!- \cdots 
   -\! 1\ket{x_n}\,\Big\}
\\[+0.5em]
& = &
\displaystyle\sum\,\Big\{\,\multinomial[n](\omega)(\chi)\,\Big|\;
   \chi\leq\tau\!-\!\one \mbox{ with } \|\chi\| = n
    \mbox{ and }\varphi = \tau \!-\! \chi\,\Big\}.
\end{array} \]

\item \label{MnclgLemMlt} For an element $x\in X$,
\[ \hspace*{-0.5em}\begin{array}{rcl}
\mnclg(\omega)^{n+1}(\tau)(x)
& = &
\displaystyle\sum\,\Big\{\,\multinomial[K](\omega)(\chi) \cdot \omega(x)
   \,\Big|\, K\leq n  \mbox{ and } \chi \leq \tau - \one
\\
& & \hspace*{5em} \mbox{ with } \|\chi\| = K  \mbox{ and } 
   (\tau\!-\!\chi)(x) = 1 \,\Big\}.
\end{array} \]
\end{enumerate}
\end{lemma}

\begin{myproof}
\begin{enumerate}
\item We first prove the first equation, by induction on $n$. The case
  $n=0$ is trivial since we have on the left-hand-side
  $\mnclg(\omega)^{0}(\tau) = 1\ket{\tau}$, and on the right-hand-side a
  sum over the empty sequence $\ell = \tuple{}$ for which by
  definition $\omega^{0}(\ell) = 1$. Next,
\[ \begin{array}{rcl}
\mnclg(\omega)^{n+1}(\tau)(\varphi)
& = &
\displaystyle\sum_{x\in X} \omega(x) \cdot 
   \mnclg(\omega)^{n}(\tau\!-\!1\ket{x}))(\varphi)
\\[+1.2em]
& \smash{\stackrel{\text{(IH)}}{=}} &
\displaystyle\sum_{x\in X} \omega(x) \cdot \sum\,\Big\{\,\omega^{n}(\ell)\,\Big|\,
   \ell=\tuple{x_{1}, \ldots, x_{n}} \in X^{n} 
\\[-0.3em]
& & \hspace*{5em}
    \mbox{with }\varphi = 
   \tau- 1\ket{x} \!-\! 1\ket{x_{1}} \!- \cdots -\! 1\ket{x_n}\,\Big\}
\\[+0.3em]
& = &
\displaystyle\sum\,\Big\{\,\omega^{n+1}(\ell)\,\Big|\,
   \ell=\tuple{x_{1}, \ldots, x_{n},x_{n+1}} \in X^{n+1} 
\\
& & \hspace*{4em}
    \mbox{with }\varphi = 
   \tau \!-\! 1\ket{x_{1}} \!- \cdots -\! 1\ket{x_n} \!-\! 1\ket{x_{n+1}}\,\Big\}.
\end{array} \]

\noindent The second equation in point~\eqref{MnclgLemLst} 
follows from Lemma~\ref{MultinomialLem}.

\item Using the previous point:
\[ \hspace*{-3em}\begin{array}[b]{rcl}
\mnclg(\omega)^{n+1}(\tau)(x)
& = &
\displaystyle\sum\,\Big\{\,\mnclg(\omega)^{K}(\tau)(\varphi) \cdot 
   \mnclg(\omega)(\varphi)(x) \,\Big|\, K\leq n \mbox{ and }
   \varphi\in\natMlt_{\geq 1} \rlap{$(X)\,\Big\}$}
\\[+0.8em]
& = &
\displaystyle\sum\,\Big\{\,\mnclg(\omega)^{K}(\tau)(\varphi) \cdot \omega(x)
   \,\Big|\, K\leq n, \,
   \varphi\in\natMlt_{\geq 1}(X) \mbox{ with } \varphi(x)
   \rlap{$\, = 1 \,\Big\}$}
\\[+0.8em]
& = &
\displaystyle\sum\,\Big\{\,\multinomial[K](\omega)(\chi) \cdot \omega(x)
   \,\Big|\, K\leq n  \mbox{ and } \chi \leq \tau - \one
\\
& & \hspace*{11em} \mbox{ with } \|\chi\| = K  \mbox{ and } 
   (\tau-\chi)(x) = 1 \,\Big\}.
\end{array} \eqno{\QEDbox} \]
\end{enumerate}
\end{myproof}

We now show that after suitably many iterations of the Markov model
with output $\mnclg(\omega)$ a multinomial first-full distribution
remains, see Definition~\ref{FirstFullDstDef}.

\begin{theorem}
\label{MnclgThm}
Let set $X$ have $N$ elements and let tubes $\tau\in\natMlt_{\geq
  1}(X)$ have size (combined length) $L = \|\tau\| \geq N$. For
$\omega\in\Dst(X)$ one has:
\[ \begin{array}{rclcrcl}
\supp\Big(\mnclg(\omega)^{L-N+1}(\tau)\Big)
& \subseteq &
X,
& \mbox{\qquad and then\qquad} &
\mnclg(\omega)^{L-N+1}(\tau)
& = &
\mnff(\omega,\tau).
\end{array} \]

\noindent In particular, this shows that multinomial first-full
$\mnff(\omega,\tau)$ is a probability distribution, with probabilities
adding up to one.
\end{theorem}

\begin{myproof}
Since $\tau \geq \one$, we have $L = \|\tau\| \geq \|\one\| = N$. With
each transition of the Markov model $\mnclg(\omega)$, say going from
multiset $\varphi$ to $\varphi'$, one has $\|\varphi'\| =
\|\varphi\|-1$. Hence after $L-N$ steps, starting from $\tau$, at most
a multiset of singletons remains. It transitions to single elements in
one step. Hence after at most $L-N+1$ steps, $\mnclg(\omega)(\tau)$
stabilises as distribution over elements $x\in X$, in the
$X$-component of $\natMlt_{\geq 1}(X)+X$. By
Lemma~\ref{MnclgLem}~\eqref{MnclgLemMlt} we then get:
\[ \begin{array}{rcl}
\mnclg(\omega)^{L-N+1}(\tau)(x)
& = &
\displaystyle\sum\,\Big\{\,\multinomial[K](\omega)(\chi) \cdot \omega(x)
   \,\Big|\, K\leq L-N  \mbox{ and } \chi\in\natMlt[K](X)
\\
& & \hspace*{6em} \mbox{ with } \chi \leq \tau - \one \mbox{ and } 
   (\tau-\chi)(x) = 1 \,\Big\}
\\[+0.5em]
& = &
\displaystyle\sum\,\Big\{\,\mulnom(\omega)(\chi) \cdot \omega(x)
   \,\Big|\; \chi \prec \tau \mbox{ with } \chi(x) = \tau(x) - 1 \,\Big\}
\\[+0.5em]
& = &
\mnff(\omega, \tau)(x), \qquad 
   \mbox{see Definition~\ref{FirstFullDstDef}~\eqref{FirstFullDstDefMulnom}.}
\end{array} \]

\noindent We use that $\chi \leq \tau - \one$ iff $\chi \prec \tau$,
where, recall, $\prec$ is the fully-below order. \QED
\end{myproof}

\subsection{Hypergeometric first-full distributions}\label{HgffCoalgSubsec}

Recall that for a hypergeometric first-full distribution we use an
urn as a multiset $\upsilon$, from which each drawn ball is actually
removed, and then dropped in the right tube. Thus, the probability of
drawing a particularly coloured ball changes throughout the filling of
the tubes.  Hence, if we wish to turn the situation into a Markov model
with output, we have to carry the urn along. This leads to the
following set-up.

Let's use the \emph{ad-hoc} notation:
\[ \begin{array}{rcl}
\natMlt_{\geq\geq 1}(X)
& \coloneqq &
\setin{(\upsilon,\tau)}{\natMlt(X)\times\natMlt(X)}{
   \upsilon \geq \tau \geq \one}.
\end{array} \]

\noindent It is the set of positions in the following hypergeometric
MMO.
\[ \vcenter{\xymatrix@R-2.5pc@C-0.6pc{
\natMlt_{\geq\geq 1}(X)\ar[r]^-{\hgclg} & 
   \Dst\Big(\natMlt_{\geq\geq 1}(X) + X\Big)
\\
(\upsilon,\tau)\ar@{|->}[r] & 
   \displaystyle\!\sum_{x,\,\tau(x)>1}\!\flrn(\upsilon)(x)
   \bigket{\upsilon\!-\!1\ket{x},\, \tau\!-\!1\ket{x}} \,+
   \displaystyle\sum_{x,\,\tau(x)=1}\!\flrn(\upsilon)(x)\bigket{x}.
}} \]

\noindent We now redo Example~\ref{FirstFullDstExHypgeom}, with urn
$\upsilon = 3\ket{a} + 6\ket{b}$ and tubes $\tau = 2\ket{a} +
3\ket{b}$. Then:

%
%
%

\[ \begin{array}{rcl}
\hgclg(\upsilon,\tau)
& = &
\frac{3}{9}\Bigket{2\ket{a} + 6\ket{b}, 1\ket{a} + 3\ket{b}} + 
   \frac{6}{9}\Bigket{3\ket{a} + 5\ket{b}, 2\ket{a} + 2\ket{b}}
\\[+0.5em]
\hgclg^{2}(\upsilon,\tau)
& = &
\frac{1}{3}\cdot\frac{2}{8}\bigket{a} + 
   \frac{1}{3}\cdot\frac{6}{8}\Bigket{2\ket{a} + 5\ket{b}, 1\ket{a} + 2\ket{b}}
\\[+0.4em]
& & \quad +\,
   \frac{2}{3}\cdot\frac{3}{8}\Bigket{2\ket{a} + 5\ket{b}, 1\ket{a} + 2\ket{b}}
  + \frac{2}{3}\cdot\frac{5}{8}\Bigket{3\ket{a} + 4\ket{b}, 2\ket{a} + 1\ket{b}}
\\[+0.5em]
\hgclg^{3}(\upsilon,\tau)
& = &
\frac{1}{12}\bigket{a} + 
   \frac{1}{2}\cdot\frac{2}{7}\bigket{a} +
   \frac{1}{2}\cdot\frac{5}{7}\Bigket{2\ket{a} + 4\ket{b}, 1\ket{a} + 1\ket{b}}
\\[+0.4em]
& & \quad +\,
   \frac{5}{12}\cdot\frac{3}{7}\Bigket{2\ket{a} + 4\ket{b}, 1\ket{a} + 1\ket{b}}
   + \frac{5}{12}\cdot\frac{4}{7}\bigket{b}
\\[+0.5em]
\hgclg^{4}(\upsilon,\tau)
& = &
\frac{19}{84}\bigket{a} + 
   \frac{15}{28}\cdot\frac{2}{6}\bigket{a} +
   \frac{15}{28}\cdot\frac{4}{6}\bigket{b} + 
   \frac{5}{21}\bigket{b}
\\[+0.4em]
& = &
\frac{17}{42}\bigket{a} + \frac{25}{42}\bigket{b}.
\end{array} \]

The next lemma makes explicit what's going on.

\begin{lemma}
\label{HgclgLem}
Let $\upsilon,\tau\in\natMlt(X)$ be multisets with
$\upsilon\geq\tau\geq\one$.
\begin{enumerate}
\item \label{HgclgLemLst} For $\upsilon',\tau'\in\natMlt(X)$ with
  $\upsilon'\geq\tau'\geq\one$,
\[ \hspace*{-2em}\begin{array}{rcl}
\hgclg^{n}(\upsilon,\tau)(\upsilon',\tau')
& = &
\displaystyle \sum\,\Big\{\, \prod_{0\leq i<n} 
   \flrn\Big(\upsilon-\acc(x_{1}, \ldots, x_{i})\Big)(x_{i+1}) \,\Big|\, 
   \ell = \tuple{x_{1}, \ldots, x_{n}} \in X^{n} 
\\[+0.2em]
& & \hspace*{3em}\mbox{with }
   \upsilon' = \upsilon - \acc(\ell) \mbox{ and } 
   \tau' = \tau - \acc(\ell) \,\Big\}
\\
& = &
\displaystyle \sum\,\Big\{\, \hypergeometric[n](\upsilon)(\chi) \,\Big|\; 
   \chi\leq\tau\!-\!\one \mbox{ with } \|\chi\| = n, \, 
   \upsilon'\! = \upsilon \!-\! \chi,
   \, \tau'\! = \tau \!-\! \chi \,\Big\}.
\end{array} \]

\item \label{HgclgLemMlt} For an element $x\in X$,
\[ \begin{array}[b]{rcl}
\hgclg^{n+1}(\upsilon,\tau)(x)
& = &
\displaystyle\sum\,\Big\{\,\hypergeometric[K](\upsilon)(\chi) \cdot 
   \flrn(\upsilon\!-\!\chi)(x)
   \,\Big|\, K\leq n \mbox{ and } \chi \leq \tau \!-\! \one 
\\
& & \hspace*{6em} \mbox{ with }  \|\chi\| = K \mbox{ and } 
   (\tau\!-\!\chi)(x) = 1 \,\Big\}.
\end{array} \]
\end{enumerate}
\end{lemma}

\begin{myproof}
\begin{enumerate}
\item We first prove the first equation by induction on $n$. The case
  $n=0$ is trivial, so we proceed with the induction step:
  
  \allowdisplaybreaks{
\begin{alignat}{2}
\lefteqn{\hgclg^{n+1}(\upsilon,\tau)(\upsilon',\tau')}
\\
&\notag =
\displaystyle\sum_{x\in X}\, \flrn(\upsilon)(x) \cdot 
   \hgclg^{n}(\upsilon-1\ket{x},\tau-1\ket{x})(\upsilon',\tau')
\\
&\notag \smash{\stackrel{\text{(IH)}}{=}} 
\displaystyle\sum_{x\in X}\, \flrn(\upsilon)(x) \cdot 
   \sum\,\Big\{\, \prod_{0\leq i<n} 
   \flrn\Big(\upsilon-1\ket{x}-\acc(x_{1}, \ldots, x_{i})\Big)(x_{i+1})
   \,\Big|\, 
\\
&\notag \hspace*{8em}
   \ell = \tuple{x_{1}, \ldots, x_{n}} \in X^{n} \mbox{ with }
   \upsilon'\! = \upsilon - 1\ket{x} - \acc(\ell)
\\
& \notag \hspace*{11em}
  \mbox{and } 
   \tau'\! = \tau - 1\ket{x} - \acc(\ell) \,\Big\}
\\
&\notag = 
\displaystyle \sum\,\Big\{\, \prod_{0\leq i<n+1} 
   \flrn\Big(\upsilon-\acc(x_{1}, \ldots, x_{i}))(x_{i+1})\Big) \,\Big|\, 
\\
&\notag \hspace*{6em} \ell = \tuple{x_{1}, \ldots, x_{n},x_{n+1}} \in X^{n+1} 
   \mbox{ with }
   \upsilon'\! = \upsilon - \acc(\ell)
\\
&\notag \hspace*{9em} \mbox{and } 
   \tau'\! = \tau - \acc(\ell) \,\Big\}.
\end{alignat}}

\noindent The second equation follows from Lemma~\ref{HypgeomLem}.

\item Via the previous point:
\[ \hspace*{-1em}\begin{array}[b]{rcl}
\hgclg^{n+1}(\upsilon,\tau)(x)
& = &
\displaystyle\sum\,\Big\{\,\hgclg^{K}(\upsilon,\tau)(\upsilon',\tau') \cdot
   \hgclg(\upsilon',\tau')(x) \,\Big|\, K\leq n \mbox{ and } \tau'(x) = 
   \rlap{$1 \,\Big\}$}
\\[+0.5em]
& = &
\displaystyle\sum\,\Big\{\,\hypergeometric[K](\upsilon)(\chi) \cdot 
   \flrn(\upsilon')(x) \,\Big|\, K\leq n, \, \|\chi\| = K, \tau'(x) = 1,
\\
& & \hspace*{6em} \chi \leq \tau \!-\! \one \mbox{ and }
   \upsilon'\! = \upsilon-\chi \mbox{ and }
   \tau'\! = \tau - \chi \,\Big\}
\\[+0.5em]
& = &
\displaystyle\sum\,\Big\{\,\hypergeometric[K](\upsilon)(\chi) \cdot 
   \flrn(\upsilon\!-\!\chi)(x)
   \,\Big|\, K\leq n \mbox{ and } \chi \leq \tau \!-\! \one 
\\
& & \hspace*{6em} \mbox{ with }  \|\chi\| = K \mbox{ and } 
   (\tau\!-\!\chi)(x) = 1 \,\Big\}.
\end{array}
\eqno{\QEDbox} \]
\end{enumerate}
\end{myproof}

We now obtain that the hypergeometric first-full probabilities form a
distribution, in the same way as in Theorem~\ref{MnclgThm} for the
multinomial mode.

\begin{theorem}
\label{HgclgThm}
Let set $X$ have $N$ elements and let tubes $\tau\in\natMlt_{\geq
  1}(X)$ have size $L = \|\tau\|$. For urn $\upsilon\geq\tau$ one
has:
\[ \begin{array}{rclcrcl}
\supp\Big(\hgclg^{L-N+1}(\upsilon,\tau)\Big)
& \subseteq &
X
& \mbox{\qquad and \qquad} &
\hgclg^{L-N+1}(\upsilon,\tau)
& = &
\hgff(\upsilon,\tau).
\end{array} \]

\noindent As a result, hypergeometric first-full $\hgff(\upsilon,\tau)$ is
a probability distribution. \QED
\end{theorem}

\auxproof{
\begin{myproof}
Termination in $L-N+1$ steps works as in (the proof of)
Theorem~\ref{MnclgThm}. Then, by Lemma~\ref{HgclgLem}~\eqref{HgclgLemMlt},
\[ \begin{array}[b]{rcl}
\lefteqn{\hgclg^{L-N+1}(\upsilon,\tau)}
\\
& = &
\displaystyle\sum\,\Big\{\,\hypgeom(\upsilon)(\chi) \cdot 
   \flrn(\upsilon\!-\!\chi)(x)
   \,\Big|\; \chi \leq \tau \!-\! \one 
   \mbox{ with }  (\tau\!-\!\chi)(x) = 1 \,\Big\}
\\[+0.5em]
& = &
\displaystyle\sum\,\Big\{\,\hypgeom(\upsilon)(\chi) \cdot 
   \flrn(\upsilon\!-\!\chi)(x)
   \,\Big|\; \chi \prec \tau
   \mbox{ with }  \chi(x) = \tau(x) - 1 \,\Big\}
\\[+0.2em]
& = &
\hgff(\upsilon,\tau), \qquad \mbox{see Definition~\ref{FirstFullDstDef}~\eqref{FirstFullDstDefHypgeom}}
\end{array} \eqno{\QEDbox} \]
\end{myproof}
}

\subsection{P\'olya first-full distributions}\label{PlffCoalgSubsec}

The P\'olya first-full mode is very similar to the hypergeometric
first-full mode, except that the drawn ball is not removed from the
urn (``-1''), but it is returned together with another ball of the
same colour (``+1''). In this case the urn is a multiset $\upsilon$ with
as only requirement $\upsilon\geq\one$ so that at least one ball of each
colour is present. We thus use a P\'olya MMO of the following form.
\[ \vcenter{\xymatrix@R-2.5pc@C-1pc{
\natMlt_{\geq 1}(X)\times \natMlt_{\geq 1}(X)\ar[r]^-{\plclg} & 
   \Dst\Big(\natMlt_{\geq 1}(X)\times\natMlt_{\geq 1}(X) + X\Big)
\\
(\upsilon,\tau)\ar@{|->}[r] & 
   \displaystyle\sum_{x,\,\tau(x)>1}\!\flrn(\upsilon)(x)
   \bigket{\upsilon\!+\!1\ket{x},\, \tau\!-\!1\ket{x}} 
   \;+\; \sum_{x,\,\tau(x)=1}\!\flrn(\upsilon)(x)\bigket{x}.
}} \]

\noindent We recalculate the outcome of
Example~\ref{FirstFullDstExPolya} as illustration, with urn $\upsilon =
1\ket{a} + 1\ket{b}$ and tubes $\tau = 2\ket{a} + 3\ket{b}$. Then:

\allowdisplaybreaks{
\begin{alignat}{2}
\plclg(\upsilon,\tau)
&\notag =
\frac{1}{2}\Bigket{2\ket{a} + 1\ket{b}, 1\ket{a} + 3\ket{b}} + 
   \frac{1}{2}\Bigket{1\ket{a} + 2\ket{b}, 2\ket{a} + 2\ket{b}}
\\
\plclg^{2}(\upsilon,\tau)
&\notag = 
\frac{1}{2}\cdot\frac{2}{3}\bigket{a} + 
   \frac{1}{2}\cdot\frac{1}{3}\Bigket{2\ket{a} + 2\ket{b}, 1\ket{a} + 2\ket{b}}
\\
&\notag  \quad +\,
   \frac{1}{2}\cdot\frac{1}{3}\Bigket{2\ket{a} + 2\ket{b}, 1\ket{a} + 2\ket{b}}
  + \frac{1}{2}\cdot\frac{2}{3}\Bigket{1\ket{a} + 3\ket{b}, 2\ket{a} + 1\ket{b}}
\\[+0.5em]
\plclg^{3}(\upsilon,\tau)
&\notag = 
\frac{1}{3}\bigket{a} + 
   \frac{1}{3}\cdot\frac{1}{2}\bigket{a} +
   \frac{1}{3}\cdot\frac{1}{2}\Bigket{2\ket{a} + 3\ket{b}, 1\ket{a} + 1\ket{b}}
\\
&\notag  \quad +\,
   \frac{1}{3}\cdot\frac{1}{4}\Bigket{2\ket{a} + 3\ket{b}, 1\ket{a} + 1\ket{b}}
   + \frac{1}{3}\cdot\frac{3}{4}\bigket{b}
\\
\plclg^{4}(\upsilon,\tau)
&\notag =
\frac{1}{2}\bigket{a} + 
   \frac{1}{4}\cdot\frac{2}{5}\bigket{a} +
   \frac{1}{4}\cdot\frac{3}{5}\bigket{b} + 
   \frac{1}{4}\bigket{b}
\\
&\notag =
\frac{3}{5}\bigket{a} + \frac{2}{5}\bigket{b}.
\end{alignat}}

We proceed with a pattern that is by now familiar. That's why we
only state the results and leave the proofs to the interested reader.

\begin{lemma}
\label{PoclgLem}
Let $\upsilon,\tau\in\natMlt(X)$ be multisets with
$\upsilon\geq\tau\geq\one$.
\begin{enumerate}
\item \label{PoclgLemLst} For $\upsilon',\tau'\in\natMlt(X)$ with
  $\upsilon',\tau'\geq\one$,
\[ \begin{array}{rcl}
\plclg^{n}(\upsilon,\tau)(\upsilon',\tau')
& = &
\displaystyle \sum\,\Big\{\, \prod_{0\leq i<n} 
   \flrn\Big((\upsilon+\acc(x_{1}, \ldots, x_{i})\Big)(x_{i+1}) \,\Big|\, 
   \ell = \tuple{x_{1}, \ldots, x_{n}} \in X^{n} 
\\[+0.2em]
& & \hspace*{3em}\mbox{with }
   \upsilon' = \upsilon + \acc(\ell) \mbox{ and } 
   \tau' = \tau - \acc(\ell) \,\Big\}
\\
& = &
\displaystyle \sum\,\Big\{\, \polya[n](\upsilon)(\chi) \,\Big|\; 
   \chi\leq\tau\!-\!\one \mbox{ with } \|\chi\| = n, \, 
   \upsilon'\! = \upsilon \!+\! \chi,
   \, \tau'\! = \tau \!-\! \chi \,\Big\}.
\end{array} \]

\item \label{PoclgLemMlt} For an element $x\in X$,
\[ \begin{array}[b]{rcl}
\plclg^{n+1}(\upsilon,\tau)(x)
& = &
\displaystyle\sum\,\Big\{\,\polya[K](\upsilon)(\chi) \cdot 
   \flrn(\upsilon\!+\!\chi)(x)
   \,\Big|\, K\leq n \mbox{ and } \chi \leq \tau \!-\! \one 
\\
& & \hspace*{6em} \mbox{ with }  \|\chi\| = K \mbox{ and } 
   (\tau\!-\!\chi)(x) = 1 \,\Big\}.
\end{array} \eqno{\QEDbox} \]
\end{enumerate}
\end{lemma}

\auxproof{
\begin{myproof}
\begin{enumerate}
\item We first prove the first equation by induction on $n$. The case
  $n=0$ is trivial, so we proceed with the induction step:
\[ \begin{array}{rcl}
\lefteqn{\plclg^{n+1}(\upsilon,\tau)(\upsilon',\tau')}
\\[+0.3em]
& = &
\displaystyle\sum_{x\in X}\, \flrn(\upsilon)(x) \cdot 
   \plclg^{n}(\upsilon+1\ket{x},\tau-1\ket{x})(\upsilon',\tau')
\\[+1.2em]
& \smash{\stackrel{\text{(IH)}}{=}} &
\displaystyle\sum_{x\in X}\, \frac{\upsilon(x)}{\|\upsilon\|} \cdot 
   \sum\,\Big\{\, \prod_{0\leq i<n} 
   \frac{(\upsilon+1\ket{x}+\acc(x_{1}, \ldots, x_{i}))(x_{i+1})}
   {\|\upsilon\|+1+i} \,\Big|\, 
\\[+0.2em]
& & \hspace*{8em}
   \ell = \tuple{x_{1}, \ldots, x_{n}} \in X^{n} \mbox{ with }
   \upsilon'\! = \upsilon + 1\ket{x} + \acc(\ell)
\\[+0.2em]
& & \hspace*{11em}
  \mbox{and } 
   \tau'\! = \tau - 1\ket{x} - \acc(\ell) \,\Big\}
\\[+0.8em]
& = &
\displaystyle \sum\,\Big\{\, \prod_{0\leq i<n+1} 
   \frac{(\upsilon+\acc(x_{1}, \ldots, x_{i}))(x_{i+1})}
   {\|\upsilon\|+i} \,\Big|\, 
   \ell = \tuple{x_{1}, \ldots, x_{n},x_{n+1}} \in X^{n+1} 
\\[+0.2em]
& & \hspace*{3em}\mbox{with }
   \upsilon'\! = \upsilon + \acc(\ell) \mbox{ and } 
   \tau'\! = \tau - \acc(\ell) \,\Big\}.
\end{array} \]

\noindent For the second equation we reason as in Lemma~\ref{OtALem}:
\[ \begin{array}{rcl}
\lefteqn{\displaystyle \sum\,\Big\{\, \prod_{0\leq i<n} 
   \frac{(\upsilon+\acc(x_{1}, \ldots, x_{i}))(x_{i+1})}
   {\|\upsilon\|+i} \,\Big|\, 
   \ell = \tuple{x_{1}, \ldots, x_{n}} \in X^{n}}
\\[+0.2em]
& & \hspace*{3em}\mbox{with }
   \upsilon'\! = \upsilon \!+\! \acc(\ell) \mbox{ and } 
   \tau'\! = \tau \!-\! \acc(\ell) \,\Big\}
\\[+0.5em]
& = &
\displaystyle \sum\,\Big\{\, \coefm{\chi} \cdot 
   \frac{(\|\upsilon\|-1)!}{(\|\upsilon\|+n-1)!}
   \cdot \prod_{y\in\supp(\upsilon)} \frac{(\upsilon(y)+\chi(y)-1)!}{(\upsilon(y)-1)!} 
   \,\Big|\; \chi\in\natMlt[n](X) 
\\[+0.2em]
& & \hspace*{3em}\mbox{with } \chi\leq\tau\!-\!\one \mbox{ and }
   \upsilon'\! = \upsilon \!+\! \chi \mbox{ and } 
   \tau'\! = \tau \!-\! \chi \,\Big\},
   \quad\mbox{by Lemma~\ref{OtALem}}
\\[+0.5em]
& = &
\displaystyle \sum\,\Big\{\, \frac{n!\cdot (\|\upsilon\|-1)!}{(\|\upsilon\|+n-1)!} \cdot
   \prod_{y\in\supp(\upsilon)} \frac{(\upsilon(y)+\chi(y)-1)!}{\chi(y)!\cdot (\upsilon(y)-1)!} 
   \,\Big|\; \chi\in\natMlt[n](X) 
\\[+0.2em]
& & \hspace*{3em}\mbox{with } \chi\leq\tau\!-\!\one \mbox{ and }
   \upsilon'\! = \upsilon \!+\! \chi \mbox{ and } 
   \tau'\! = \tau \!-\! \chi \,\Big\},
\\[+0.5em]
& = &
\displaystyle \sum\,\Big\{\, \frac{\bibinom{\upsilon}{\chi}}
   {\bibinom{\|\upsilon\|}{\|\chi\|}} \,\Big|\; 
   \chi\leq\tau\!-\!\one \mbox{ with } \|\chi\| = n,
   \, \upsilon'\! = \upsilon \!+\! \chi, \, 
   \tau'\! = \tau \!-\! \chi \,\Big\},
\\[+0.5em]
& = &
\displaystyle \sum\,\Big\{\, \polya[n](\upsilon)(\chi) \,\Big|\; 
   \chi\leq\tau\!-\!\one \mbox{ with } \|\chi\| = n, \, 
   \upsilon'\! = \upsilon \!+\! \chi,
   \, \tau'\! = \tau \!-\! \chi \,\Big\}.
\end{array} \]

\item Via the previous point:
\[ \hspace*{-1em}\begin{array}[b]{rcl}
\lefteqn{\plclg^{n+1}(\upsilon,\tau)(x)}
\\[+0.3em]
& = &
\displaystyle\sum\,\Big\{\,\plclg^{K}(\upsilon,\tau)(\upsilon',\tau') \cdot
   \plclg(\upsilon',\tau')(x) \,\Big|\, K\leq n \mbox{ and } \tau'(x) = 
   \rlap{$1 \,\Big\}$}
\\[+0.5em]
& = &
\displaystyle\sum\,\Big\{\,\polya[K](\upsilon)(\chi) \cdot 
   \flrn(\upsilon')(x) \,\Big|\, K\leq n, \, \|\chi\| = K, \tau'(x) = 1,
\\
& & \hspace*{6em} \chi \leq \tau \!-\! \one \mbox{ and }
   \upsilon'\! = \upsilon+\chi \mbox{ and }
   \tau'\! = \tau - \chi \,\Big\}
\\[+0.5em]
& = &
\displaystyle\sum\,\Big\{\,\polya[K](\upsilon)(\chi) \cdot 
   \flrn(\upsilon\!+\!\chi)(x)
   \,\Big|\, K\leq n \mbox{ and } \chi \leq \tau \!-\! \one 
\\
& & \hspace*{6em} \mbox{ with }  \|\chi\| = K \mbox{ and } 
   (\tau\!-\!\chi)(x) = 1 \,\Big\}.
\end{array} \eqno{\QEDbox} \]
\end{enumerate}
\end{myproof}
}

\begin{theorem}
\label{PoclgThm}
Let set $X$ have $N$ elements and let tubes $\tau\in\natMlt_{\geq
  1}(X)$ have size $L = \|\tau\| \geq N$. For urn $\upsilon\geq\one$ one
gets:
\[ \begin{array}{rclcrcl}
\supp\Big(\plclg^{L-N+1}(\upsilon,\tau)\Big)
& \subseteq &
X
& \qquad\mbox{and}\qquad &
\plclg^{L-N+1}(\upsilon,\tau)
& = &
\plff(\upsilon,\tau).
\end{array} \]

\noindent In particular, P\'olya first-full $\plff(\upsilon,\tau)$ is a
probability distribution. \QED
\end{theorem}

\auxproof{
\begin{myproof}
Termination in $L-N+1$ steps works as before. Then, by
Lemma~\ref{PoclgLem}~\eqref{PoclgLemMlt},
\[ \begin{array}[b]{rcl}
\lefteqn{\plclg^{L-N+1}(\upsilon,\tau)}
\\
& = &
\displaystyle\sum\,\Big\{\,\pol(\upsilon)(\chi) \cdot 
   \flrn(\upsilon\!+\!\chi)(x)
   \,\Big|\; \chi \leq \tau \!-\! \one 
   \mbox{ with }  (\tau\!-\!\chi)(x) = 1 \,\Big\}
\\[+0.5em]
& = &
\displaystyle\sum\,\Big\{\,\pol(\upsilon)(\chi) \cdot 
   \flrn(\upsilon\!+\!\chi)(x)
   \,\Big|\; \chi \prec \tau
   \mbox{ with }  \chi(x) = \tau(x) - 1 \,\Big\}
\\[+0.2em]
& = &
\plff(\upsilon,\tau), \qquad \mbox{see Definition~\ref{FirstFullDstDef}~\eqref{FirstFullDstDefPolya}}
\end{array} \eqno{\QEDbox} \]
\end{myproof}
}

\begin{remark}
\label{FirstFullImplementationRem}
\emph{In the end we have two (equivalent) ways to compute first-full
distributions, namely via their pointwise formulations (in
Definition~\ref{FirstFullDstDef}) and via their three MMO's, as
described above. The latter can easily be turned into recursive
definitions. Experiments with both implementations show that the
recursive approach is slower than the one based on the
definitions. This is not surprising since the recursive approach
computes probabilities for sequences, with much duplication, instead
of for multisets (as accumulations of those sequences). In contrast,
reasoning with sequences is easier than with multisets. In the end
that is the whole reason why we use the MMO-approach --- for proving
that combined first-full probabilities form a distribution.}
\end{remark}

\section{Negative distributions}\label{NegativeDrawSec}


The urns \& tubes set-up that we have used to introduce first-full
distributions can also be used to describe `negative' distributions.
The latter are known from the literature, in bivariate form, with one
tube only.  Here we use our multiset-based approach to describe them
systematically, in multivariate form, for all three modes
(multinomial, hypergeometric and P\'olya). We concentrate on the
definitions and on illustrations. The fact that these definitions lead
to actual probability distributions is addressed, by giving the
corresponding Markov models with output, but without all the
mathematical details. After all, these distributions are not new.

We start with an example.  Consider a group of people consisting of
five males ($M$) and four females ($F$). From this group we like to
form a committee with two male and two female members. Iteratively we
choose members from the group, at random, until the committee is
formed. How many choices are needed? More precisely, what is the
probability --- that the committee is first formed --- for each number
of choices. We might be done after four choices, if they immediately
involve two males and two females. But we may also first pick three
men, and then two females, which involves five choices. What is the
highest possible number of choices? It is seven, when all five men are
chosen, before the two females. Hence this situation involves a
distribution on the set $\{4,5,6,7\}$. It is:
\[ \textstyle\frac{10}{21}\bigket{4} \,+\, \frac{20}{63}\bigket{5} \,+\, 
   \frac{10}{63}\bigket{6} \,+\, \frac{1}{21}\bigket{7}
\quad\mbox{in}\quad
\vcenter{\hbox{\includegraphics[width=7cm]{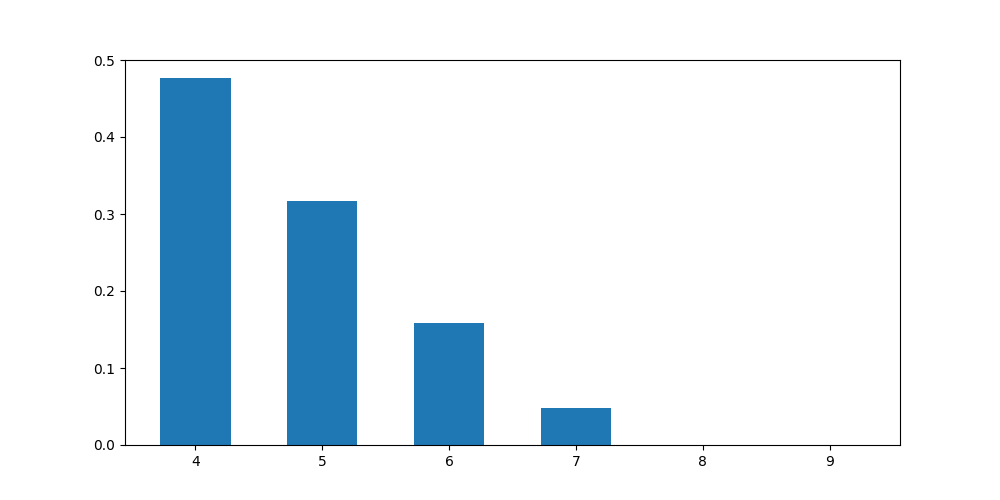}}} \]

\noindent How does it come about? We can see the group of people as an
urn $\upsilon = 5\ket{M} + 4\ket{F}$ from which we `draw' candidate
committee members, in hypergeometric mode: after drawing a member from
the urn/group, this person is either put in the committee, or is
skipped, when there are already two committee members with this
person's gender in the committee.
\begin{itemize}
\item The probability of being done in four steps is given by the
  hypergeometric distribution at multiset $2\ket{M} +
  2\ket{F}$. Indeed,
\[ \begin{array}{rcl}
\hypergeometric[4](\upsilon)\Big(2\ket{M} + 2\ket{F}\Big)
& = &
\frac{10}{21}.
\end{array} \]

\item We may need five steps in two cases: (1)~when we first choose three
  men and one woman, in any order, and finally a woman, or (2)~when we
  first choose three women and one man, in any order, and finally a
  man. This leads to the probability:
\[ \begin{array}{rcccl}
\hypergeometric[4](\upsilon)\Big(3\ket{M} + 1\ket{F}\Big) \cdot \frac{3}{5} +
\hypergeometric[4](\upsilon)\Big(1\ket{M} + 3\ket{F}\Big) \cdot \frac{4}{5}
& = &
\frac{20}{63} \cdot \frac{3}{5} + \frac{10}{63} \cdot \frac{4}{5}
& = &
\frac{20}{63}. 
\end{array} \]

\item We need six steps when we first choose four men and one woman, and
  then one woman, or when we first choose one man and four women. In
  the latter case we are done at the next selection, because we can
  only choose a male. The associated probability is thus obtained as:
\[ \begin{array}{rcccl}
\hypergeometric[5](\upsilon)\Big(4\ket{M} + 1\ket{F}\Big) \cdot \frac{3}{4} +
\hypergeometric[5](\upsilon)\Big(1\ket{M} + 4\ket{F}\Big)
& = &
\frac{10}{63} \cdot \frac{3}{4} + \frac{5}{126}
& = &
\frac{10}{63}.
\end{array} \]

\item Finally, there is only one possibility that requires the maximum
  number of seven steps, namely when we first choose five men and one woman.
  The associated probability is simply:
\[ \begin{array}{rcl}
\hypergeometric[6](\upsilon)\Big(5\ket{M} + 1\ket{F}\Big)
& = &
\frac{1}{21}.
\end{array} \]
\end{itemize}












In this example we may consider the committee that needs to be filled
with two males and two females as a pair of tubes, both of length two.
Thus, the urns \& tubes model can be used here as well, but with a
different question, namely what is the probability of filling \emph{all}
tubes in a certain number of steps.

Thus, abstractly, our starting point is the same as in the previous
section, see Picture~\ref{UrnTubesPic}: we have an urn filled with
coloured balls, together with coloured tubes. The question that we now
look at is as follows.
\begin{quote}
Suppose we draw $k$ balls from the urn, to fill the tubes, for
$k\in\NNO$. What is the probability that all tubes are full for the
first time after drawing these $k$ balls? This means that there is one
tube that becomes full with the $k$-th ball, while sufficiently many
balls --- typically more than needed --- have already been drawn to
fill all other tubes.
\end{quote}

\noindent The filling of all tubes can be seen as a desired condition,
or as a risk. The probability distribution that we are after gives for
each $k\in\NNO$ the probability of reaching this threshold condition
for the first time. 
\begin{itemize}
\item Historically the distributions that arise in this manner are
  called \emph{negative}. They are not very well known, and are even
  called `forgotten' in~\cite{MillerF07}. These negative distributions
  may occur in different forms, depending on the mode of drawing
  (``-1'', ``0'', or ``+1''). Accordingly, we shall speak of
  \emph{negative hypergeometric}, \emph{negative multinomial}, and
  \emph{negative P\'olya} distributions. Below we cover all three
  modes.

\item The negative hypergeometric distribution has finite support,
  since at some stage the urn is empty. In multinomial mode the urn
  does not change, and in P\'olya mode the urn grows in size. Hence in
  these last two cases the support of the negative distributions are
  infinite subsets of the natural numbers. Recall that we write
  $\infDst(\NNO)$ for set of discrete distributions on $\NNO$ with
  (possibly) infinite support.

\item In the literature (see
  \eg~\cite{JohnsonKK05,MillerF07,Panaretos81,SchusterS87,SibuyaYS64})
  negative distributions are studied only for the (simple) case with a
  single tube and usually with only two colours. Here we deal with the
  general, multivariate and multi-tube, scenario, where there are
  multiple colours and as many tubes as colours. Like before, we shall
  write $\tau$ for the tubes and $L \coloneqq \|\tau\|$ for the sum of
  the lengths of all tubes. We assume that $L > 0$ so that there is at
  least one non-empty tube that can be filled. We shall write
  $\nenatMlt(X) \hookrightarrow \natMlt(X)$ for the subset of
  \emph{non-empty} multisets; we thus require
  $\tau\in\nenatMlt(X)$. The negative distributions on $\NNO$ that we
  are after will `start at $L$': they are zero at $k<L$, since one
  needs to draw at least $L$ balls to fill all tubes.

\item As in the first-fill case, in Section~\ref{FirstFullDefSec},
  there is a challenge to show that negative probabilities add up to
  one, and thus form a proper distribution. Again we use Markov models
  with output, like in the previous section, but without elaborating
  all details. Previously, we only had finitely many possible
  transitions. Here, in the negative setting, there may be infinitely
  many transitions, leading to infinite supports.
\end{itemize}

\begin{definition}
\label{NegativeDstDef}
Let $X$ be a finite set (of colours) with $\setsize{X} \geq 2$ and let
$\tau\in\nenatMlt(X)$ be an $X$-indexed collection of tubes.
\begin{enumerate}
\item \label{NegativeDstDefMulnom} For $\omega\in\Dst(X)$, we define
  at $k>0$ the \emph{negative multinomial} probability as:
\[ \begin{array}{rcl}
\negmultinomial(\omega, \tau)(k)
& \coloneqq &
\displaystyle\sum_{x\in \supp(\tau)} 
   \sum_{\smash{\begin{array}{c} \scriptstyle \\[-0.5em]
   \scriptstyle \varphi\in\natMlt[k-1](X), \\[-0.6em]
   \scriptstyle \tau-1\ket{x} \,\leq\, \varphi, \, \varphi(x) = \tau(x) - 1 
   \end{array}}}  \!\! \multinomial[k\!-\!1](\omega)(\varphi)\cdot\omega(x).
\end{array} \]

\smallskip

\item \label{NegativeDstDefHypgeom} For an urn $\upsilon\in\natMlt(X)$
  with $\upsilon\geq\tau$, we define the \emph{negative hypergeometric}
  probability at $k>0$ as:
\[ \begin{array}{rcl}
\neghypergeometric(\upsilon, \tau)(k)
& \coloneqq &
\displaystyle\sum_{x\in \supp(\tau)} 
   \sum_{\smash{\begin{array}{c} \scriptstyle \\[-0.5em]
   \scriptstyle \varphi\,\leq_{k-1}\,\upsilon, \\[-0.6em]
   \scriptstyle \tau-1\ket{x} \,\leq\, \varphi, 
   \, \varphi(x) = \tau(x) - 1 
   \end{array}}}  \!\! \hypgeom[k\!-\!1](\upsilon)(\varphi)\cdot
   \flrn(\upsilon-\varphi)(x).
\end{array} \]

\smallskip

\item \label{NegativeDstDefPolya} Finally, for an urn
  $\upsilon\in\natMlt(X)$ with $\upsilon\geq\one$, and for $k>0$ we define the
  \emph{negative P\'olya} probability as:
\[ \begin{array}{rcl}
\negpolya(\upsilon, \tau)(k)
& \coloneqq &
\displaystyle\sum_{x\in \supp(\tau)} 
   \sum_{\smash{\begin{array}{c} \scriptstyle \\[-0.5em]
   \scriptstyle \varphi\in\natMlt[k-1](X), \\[-0.6em]
   \scriptstyle \tau-1\ket{x} \,\leq\, \varphi, 
   \, \varphi(x) = \tau(x) - 1 
   \end{array}}}  \!\! \polya[k\!-\!1](\upsilon)(\varphi)\cdot
   \flrn(\upsilon+\varphi)(x).
\end{array} \]

\smallskip

\end{enumerate}

\noindent These distributions may be extended to $k=0$ by setting them
to zero there.
\end{definition}

In each of the above three cases we sum over draws $\varphi$
satisfying $\tau-1\ket{x} \leq \varphi$ and $\varphi(x) = \tau(x) -
1$. The inequality $\leq$ implies that $\tau(y) \leq \varphi(y)$ for
all $y\neq x$, so that all tubes are full (possibly with overflow)
after drawing $\varphi$, except for colour $x$. The equality
$\varphi(x) = \tau(x) - 1$ says that there is precisely one ball of
colour $x$ missing to ensure that all tubes are full. The probability
of additionally drawing this missing ball of colour $x$ is multiplied
in each of the above three cases with the probability of the draw
$\varphi$, as $\omega(x)$ in item~\eqref{NegativeDstDefMulnom}, as
$\flrn(\upsilon-\varphi)(x)$ in item~\eqref{NegativeDstDefHypgeom},
and as $\flrn(\upsilon+\varphi)(x)$ in
item~\eqref{NegativeDstDefPolya}. In the first-fill probabilities in
Definition~\ref{FirstFullDstDef} this is done analogously.

Computing negative distributions by hand is laborious because it
involves summing over all colours and over all draws, of a certain
size. However, this can be automated without too much effort.

\begin{example}
\label{NegativeDstEx}
Take $X = \{a,b,c\}$ with tubes $\tau = 2\ket{a} + 4\ket{b} +
3\ket{c}$, having total length $L = \|\tau\| = 9$. For a state $\omega
= \frac{1}{6}\ket{a} + \frac{1}{2}\ket{b} + \frac{1}{3}\ket{c}$. A
first part of the resulting negative multinomial distribution
$\negmultinomial(\omega,\tau)$ on $\NNO$ looks as follows.
\begin{center}
\includegraphics[width=11cm]{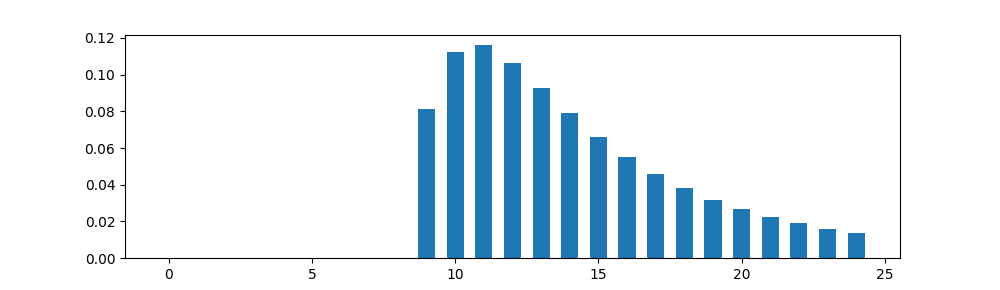}
\end{center}


\noindent The sum of the probabilities in this picture is
approximately $0.92$. The remaining $0.08$ is in the long tail. An
exact description of the first four probabilities is:
\begin{equation}
\label{NegativeMulnomFracEx}
\begin{array}{rcl}
\negmultinomial(\omega, \tau)
& = &
\frac{35}{432}\ket{9} + \frac{875}{7776}\ket{10} + 
   \frac{3605}{31104}\ket{11} + \frac{1243}{11664}\ket{12} + \cdots
\end{array}
\end{equation}


With urn $\upsilon = 10\ket{a} + 6\ket{b} + 8\ket{c}$ the negative
hypergeometric distribution $\neghypergeometric(\upsilon,\tau)$ runs from
$L=\|\tau\|=9$ to $\|\upsilon\| - 1 = 23$ and looks in its entirety as follows.
\begin{center}
\includegraphics[width=11cm]{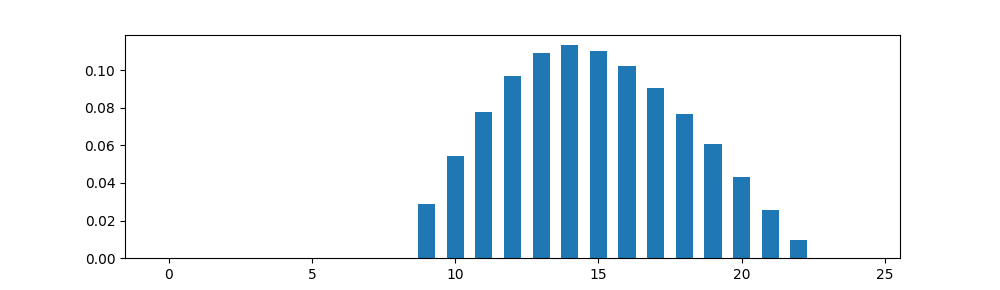}
\end{center}


Using urn $\upsilon = 3\ket{a} + 2\ket{b} + 1\ket{c}$ the negative P\'olya
distribution $\negpolya(\upsilon, \tau)$ on $\NNO$ starts as described
below.
\begin{center}
\includegraphics[width=11cm]{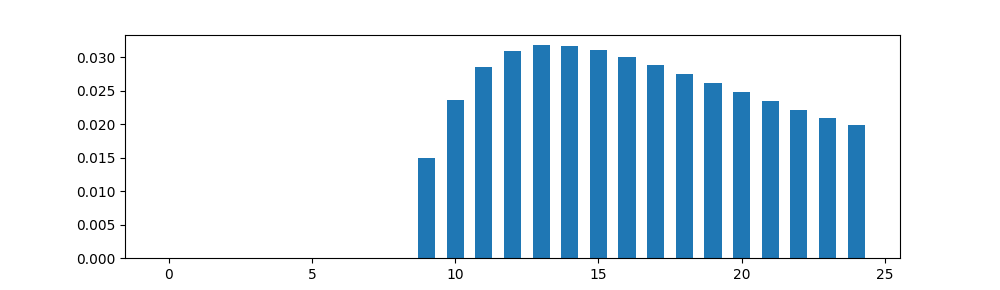}
\end{center}

\noindent This pictures only contains about $0.42$ of all
probabilities. This negative P\'olya is thus heavy-tailed and its
probabilities are less concentrated at the beginning than in the
negative multinomial.




\end{example}

\subsection{The common one-tube situation}\label{OneTubeSubsec}

In the introduction to this section we mentioned that the negative
distributions that are commonly considered in the literature involve
one tube only.  We describe what happens then, as special case of the
above formulations in Definition~\ref{NegativeDstDef}, and recover
familiar formulations. It turns out, in all three drawing modes, that
the relevant probabilities can also be described via the
`non-negative' bivariate distribution.

\begin{theorem}
\label{SingleTubeThm}
Let $X$ be a set of colours, with a special fixed element $y\in X$,
and with single-tube multiset $\tau = m\ket{y}$ for $m > 0$.
\begin{enumerate}
\item \label{SingleTubeThmMulnom} Let $\omega\in\Dst(X)$ satisfy $0 <
  \omega(y) < 1$. For $k\geq 0$,
\[ \begin{array}{rcl}
\negmultinomial\big(\omega, \, m\ket{y}\big)(m\!+\!k)
& = &
\displaystyle \frac{m}{m\!+\!k}\cdot 
   \multinomial[m\!+\!k]\Big(\omega(y)\ket{0} + (1\!-\!\omega(y))\ket{1}\Big)
    \Big(m\ket{0} + k\ket{1}\Big)
\\[+0.5em]
& = &
\negbinomial[m]\big(\omega(y)\big)(m\!+\!k).
\end{array} \]

\noindent The latter expression involves the negative binomial
distribution, of the form:
\begin{ceqn}
\begin{equation}
\label{NegBinomEqn}
\begin{array}{rcl}
\negbinomial[m](s)
& \coloneqq &
\displaystyle \sum_{i \geq 0} \bibinom{m}{i} 
   \cdot s^{m}\cdot (1-s)^{i} \,\bigket{m+i} \,\in\, \infDst(\NNO).
\end{array}
\end{equation}
\end{ceqn}

\item \label{SingleTubeThmHypgeom} Similarly, for an urn
  $\upsilon\in\natMlt[L](X)$ with $\upsilon(y) \geq m$, one has, for $k\geq
  0$,
\[ \begin{array}{rcl}
\neghypergeometric\big(\upsilon, \, m\ket{y}\big)(m\!+\!k)
& = &
\displaystyle \frac{m}{m\!+\!k} \cdot
   \hypergeometric[m\!+\!k]\big(\upsilon(y)\ket{0} + (L\!-\!\upsilon(y))\ket{1}\big)
    \big(m\ket{0} + k\ket{1}\big)
\\[+0.5em]
& = &
\displaystyle \frac{m}{m\!+\!k} \cdot
   \frac{\binom{\upsilon(y)}{m}\cdot\binom{L-\upsilon(y)}{k}}
   {\binom{L}{m+k}}.
\end{array} \]

\item \label{SingleTubeThmPolya} Also negative P\'olya with one tube
  reduces to bivariate non-negative form. For an urn
$\upsilon\in\natMlt[K](X)$ with $\upsilon(y) > 0$ one has:
\[ \begin{array}{rcl}
\negpolya\big(\upsilon, \, m\ket{y}\big)(m\!+\!k)
& = &
\displaystyle \frac{m}{m\!+\!k} \cdot
   \polya[m\!+\!k]\big(\upsilon(y)\ket{0} + (L\!-\!\upsilon(y))\ket{1}\big)
    \big(m\ket{0} + k\ket{1}\big)
\\[+0.5em]
& = &
\displaystyle \frac{m}{m\!+\!k} \cdot
   \frac{\big(\!\binom{\upsilon(y)}{m}\!\big)\cdot
      \big(\!\binom{L-\upsilon(y)}{k}\!\big)}
   {\big(\!\binom{L}{m+k}\!\big)}.
\end{array} \]
\end{enumerate}
\end{theorem}




















\begin{myproof}
\begin{enumerate}
\item In presence of a single tube the negative multinomial becomes a
  single sum over multisets:

\[ \begin{array}{rcl}
\negmultinomial\big(\omega, m\ket{y}\big)(m\!+\!k)
& = &
\displaystyle \sum_{\varphi\in\natMlt[m+k-1](X), \, \varphi(y) = m-1}
   \multinomial[m\!+\!k\!-\!1](\omega)(\varphi)\cdot\omega(y)
\\[+1.3em]
& = &
\displaystyle \sum_{\varphi\in\natMlt[k](X-y)} \, 
   \multinomial[m\!+\!k\!-\!1](\omega)\big((m\!-\!1)\ket{y} + \varphi\big)
   \cdot\omega(y)
\\[+1em]
& = &
\displaystyle \sum_{\varphi\in\natMlt[k-m](X-y)} \, 
   \frac{(m\!+\!k\!-\!1)!}{(m\!-\!1)! \cdot \facto{\varphi}} \cdot
   \omega(y)^{m} \cdot \prod_{x\neq y} \omega(x)^{\varphi(x)}
\\[+1.5em]
& = &
\displaystyle \binom{m\!+\!k\!-\!1}{m\!-\!1} \cdot \omega(y)^{m} \cdot 
   \sum_{\varphi\in\natMlt[k](X-y)} \, \frac{k!}{\facto{\varphi}} \cdot
   \prod_{x\neq y} \omega(x)^{\varphi(x)}
\\[+1em]
& = &
\displaystyle \bibinom{m}{k} \cdot \omega(y)^{m} \cdot 
   \textstyle\big(\sum_{x\neq y}\omega(x)\big)^{k} \qquad
   \mbox{by Fact~\ref{CoefficientFact}~\eqref{CoefficientFactMulnomThm}}
\\[+1em]
& = &
\displaystyle \bibinom{m}{k} \cdot \omega(y)^{m} \cdot 
   (1-\omega(y))^{k}
\\[+0.8em]
& = &
\negbinomial[m]\big(\omega(y)\big)(m\!+\!k).
\end{array} \]

\noindent At the same time we can write:
\[ \hspace*{-1em}\begin{array}{rcl}
\displaystyle\bibinom{m}{k} \cdot \omega(y)^{m} \cdot (1-\omega(y))^{k}
& = &
\displaystyle \frac{m}{m\!+\!k}\cdot 
   \binom{m\!+\!k}{k} \cdot \omega(y)^{m} \cdot (1-\omega(y))^{k}
\\[+1.0em]
& = &
\displaystyle \frac{m}{m\!+\!k}\cdot 
   \multinomial[m\!+\!k]\Big(\omega(y)\ket{0} + (1\!-\!\omega(y))\ket{1}\Big)
    \Big(m\ket{0} + k\ket{1}\Big).
\end{array} \]

\auxproof{
We use:
\[\begin{array}{rcl}
\displaystyle\bibinom{n}{m}
\hspace*{\arraycolsep}=\hspace*{\arraycolsep}
\displaystyle\binom{n+m-1}{m}
& = &
\displaystyle\frac{(n+m-1)!}{m!\cdot (n-1)!}
\\[+1em]
& = &
\displaystyle\frac{n}{n+m} \cdot \frac{(n+m-1)!}{m!\cdot n!}
\hspace*{\arraycolsep}=\hspace*{\arraycolsep}
\displaystyle\frac{n}{n+m} \cdot \binom{n+m}{m}.
\end{array} \]
}

\item We write $\upsilon' = \upsilon - \upsilon(y)\ket{y}$ for the urn from
which all balls of colour $y$ have been removed.
\[ \begin{array}{rcl}
\lefteqn{\neghypergeometric\big(\upsilon, m\ket{y}\big)(m\!+\!k)}
\\[+1em]
& = &
\displaystyle \sum_{\varphi\leq_{m+k-1} \upsilon, \, \varphi(y) = m-1}
   \hypergeometric[m\!+\!k\!-\!1](\upsilon)(\varphi)
   \cdot\flrn\big(\upsilon-\varphi\big)(y)
\\[+1.3em]
& = &
\displaystyle \sum_{\varphi\leq_{k} \upsilon'} \, 
   \hypergeometric[m\!+\!k\!-\!1](\upsilon)\big((m\!-\!1)\ket{y} + \varphi\big)
   \cdot\flrn\Big(\upsilon-(m\!-\!1)\ket{y} - \varphi\Big)(y)
\\[+1em]
& = &
\displaystyle \sum_{\varphi\leq_{k} \upsilon'} \, 
    \frac{\binom{\upsilon(y)}{m-1} \cdot \prod_{x\neq y} \binom{\upsilon(x)}{\varphi(x)}}
         {\binom{L}{m+k-1}} \cdot 
    \frac{\upsilon(y) - (m-1)}{L - (m+k-1)} 
    \qquad \mbox{where }L = \|\upsilon\|
\\[+1.3em]
& = &
\displaystyle \frac{m\cdot \binom{\upsilon(y)}{m} \cdot 
   \sum_{\varphi\leq_{k} \upsilon'} \binom{\upsilon'}{\varphi}}
         {(m\!+\!k)\cdot\binom{L}{m+k}}
   \qquad\mbox{since }\; \textstyle \binom{n}{m}\cdot (n\!-\!m) = 
   (m\!+\!1)\cdot\binom{n}{m+1}
\\[+1.3em]
& = &
\displaystyle \frac{m}{m\!+\!k} \cdot
   \frac{\binom{\upsilon(y)}{m} \cdot \binom{L-\upsilon(y)}{k}}{\binom{L}{m+k}}
   \qquad \mbox{by Proposition~\ref{HypgeomCountProp}}
\\[+1.2em]
& = &
\displaystyle \frac{m}{m\!+\!k} \cdot
   \hypergeometric[m\!+\!k]\big(\upsilon(y)\ket{0} + (L\!-\!\upsilon(y))\ket{1}\big)
    \big(m\ket{0} + k\ket{1}\big).
\end{array} \]

\item In the P\'olya case we proceed in a similar manner, for
$\upsilon\in\natMlt[L](X)$.

\allowdisplaybreaks{
\begin{alignat}{2}
&\notag\negpolya\big(\upsilon, m\ket{y}\big)(m\!+\!k)\\
&\notag =
\displaystyle \sum_{\varphi\in\natMlt[m+k-1](X), \, \varphi(y) = m-1}
   \polya[m\!+\!k\!-\!1](\upsilon)(\varphi)\cdot\flrn(\upsilon+\varphi)(y)
\\
&\notag = 
\displaystyle \sum_{\varphi\in\natMlt[k](X-y)} \, 
   \polya[m\!+\!k\!-\!1](\upsilon)\big((m\!-\!1)\ket{y} + \varphi\big)
\\
&\notag  \hspace*{10em}
   \cdot\,\flrn\Big(\upsilon+(m\!-\!1)\ket{y} + \varphi\Big)(y)
\\
&\notag = 
\displaystyle \sum_{\varphi\in\natMlt[k](X-y)} \, 
   \frac{\big(\!\binom{\upsilon(y)}{m-1}\!\big) \cdot 
      \prod_{x\neq y} \big(\!\binom{\upsilon(x)}{\varphi(x)}\!\big)}
         {\big(\!\binom{L}{m+k-1}\!\big)} \cdot 
    \frac{\upsilon(y) + (m-1)}{L + (m+k-1)}
\\
&\notag = 
\displaystyle \frac{m\cdot \big(\!\binom{\upsilon(y)}{m}\!\big) \cdot 
   \sum_{\varphi\in\natMlt[k](X-y)} \big(\!\binom{\upsilon'}{\varphi}\!\big)}
         {(m\!+\!k)\cdot\big(\!\binom{L}{m+k}\!\big)}
   \qquad \mbox{where } \upsilon' = \upsilon - \upsilon(y)\ket{y}
\\
&\notag = 
\displaystyle \frac{m}{m\!+\!k} \cdot
   \frac{\big(\!\binom{\upsilon(y)}{m}\!\big) \cdot 
    \big(\!\binom{L-\upsilon(y)}{k}\!\big)}{\big(\!\binom{L}{m+k}\!\big)}
   \qquad \mbox{by Proposition~\ref{PolyaCountProp}}
\\
& \notag= 
\displaystyle \frac{m}{m\!+\!k} \cdot
   \polya[m\!+\!k]\big(\upsilon(y)\ket{0} + (L\!-\!\upsilon(y))\ket{1}\big)
    \big(m\ket{0} + k\ket{1}\big).
\end{alignat} }
 \flushright \QEDbox

%

\end{enumerate}
\end{myproof}

\auxproof{
Let's abbreviate $\upsilon_{x} = \upsilon - \upsilon(x)\ket{x}$.
\[ \begin{array}{rcl}
\lefteqn{\mean\big(\neghypergeometric(\upsilon, \tau)\big)}
\\
& = &
\displaystyle\sum_{k > 0} \neghypergeometric(\upsilon, \tau)(k)\cdot k
\\[+1.2em]
& = &
\displaystyle\sum_{k > 0} \sum_{x\in\supp(\tau)} 
   \sum_{\varphi \leq_{k-\tau(x)} \upsilon_{X}} 
   \hypergeometric[k\!-\!1](\upsilon)\big((\tau(x)\!-\!1)\ket{x} + \varphi\big)
\\
& & \hspace*{12em}
   \cdot\, \flrn(\upsilon(x) - (\tau(x)\!-\!1)\ket{x} - \varphi)(x) \cdot k
\\[+1.2em]
& = &
\displaystyle\sum_{x\in\supp(\tau)}  \sum_{k\geq\tau(x)} \sum_{\varphi \leq_{k-\tau(x)} \upsilon_{x}}
   \frac{\binom{\upsilon(x)}{\tau(x)-1} \cdot 
      \prod_{y\neq x} \binom{\upsilon_{x}(y)}{\varphi(y)}}{\binom{\|\upsilon\|}{k-1}}
   \cdot \frac{\upsilon(x)\!-\!\tau(x)\!+\!1}{\|\upsilon\|\!-\!k\!+\!1} \cdot k
\\[+1.2em]
& = &
\displaystyle\sum_{x\in\supp(\tau)}  \sum_{k\geq\tau(x)} 
   \frac{\tau(x) \cdot \binom{\upsilon(x)}{\tau(x)} \cdot 
    \sum_{\varphi \leq_{k-\tau(x)} \upsilon_{x}} \binom{\upsilon_{x}}{\varphi}}
    {k \cdot \binom{\|\upsilon\|}{k}} \cdot k
\\[+1.2em]
& = &
\displaystyle\sum_{x\in\supp(\tau)}  \tau(x) \cdot \sum_{k\geq\tau(x)} 
   \frac{\binom{\upsilon(x)}{\tau(x)} \cdot 
    \binom{\|\upsilon\|-\upsilon(x)}{k-\tau(x)}}{\binom{\|\upsilon\|}{k}}
\\[+1.2em]
& = &
\displaystyle\sum_{x\in\supp(\tau)}  \tau(x) \cdot \sum_{k\geq\tau(x)} 
   \hypergeometric[k]\big(\upsilon(x)\ket{x} + (\|\upsilon\|-\upsilon(x))\ket{\no{x}}\big)
   \big(\tau(x)\ket{x} + (k-\tau(x))\ket{\no{x}}\big)
\\[+1.2em]
& = &
\displaystyle\sum_{x\in\supp(\tau)}  \tau(x) \cdot \sum_{0\leq k\leq\|\upsilon\|-\tau(x)} 
   \hypergeometric[\tau(x)\!+\!k]
   \big(\upsilon(x)\ket{x} + (\|\upsilon\|-\upsilon(x))\ket{\no{x}}\big)
   \big(\tau(x)\ket{x} + k\ket{\no{x}}\big)
\\[+1.2em]
& = &
\displaystyle\sum_{x\in\supp(\tau)}  \tau(x) \cdot 
   \frac{\|\upsilon\|+1}{\upsilon(x)+1} \qquad 
   \mbox{by Exericse~\ref{BinHypgeomSumExc}}.
\end{array} \]
}

\subsection{Negatives yield distributions}\label{NegativeDstSubsec}

We will illustrate that the probabilities in the `negative'
formulations in Definition~\ref{NegativeDstDef} yield actual
distributions. We shall proceed as in Section~\ref{FirstFullDstSec}
and introduce appropriate Markov models with output (MMO). The
positions in these MMOs are tuples involving a `stage' number
$i\in\NNO$. Each step involves one of the following three options.
\begin{enumerate}
\item For colour $x$ with already full tube, so $\tau(x)=0$, one can
  draw another ball of colour $x$ and move to a next position in the
  MMO. We then have an overflow situation for colour $x$ so this next
  position has the same tubes $\tau$ and an incremented stage, which
  change from $i$ to $i+1$.

\item In case colour $x$ is the last one whose tube needs to be
  filled, we have $\tau(x) > 0$ and $\|\tau\|=1$, or equivalently
  $\tau(x) = 1$ and $\tau(y) = 0$ for all $y\neq x$.  Then we can draw
  this last ball and move to the output $i+1$, from which no further
  transitions are possible.

\item When colour $x$'s tube is not full yet, and there are other
  non-full tubes as well --- so when $\tau(x) > 0$ and $\|\tau\| > 1$
  --- we can draw a ball of colour $x$ and drop it in the tube of
  colour $x$. The next position then involves new tubes $\tau -
  1\ket{x}$, where the number of missing balls of colour $x$ is
  reduced by $1$, with incremented stage $i+1$.
\end{enumerate}

\subsection*{Negative multinomial distributions} 

The MMO for negative multinomials has pairs $(\tau,i)$ as positions,
with tubes $\tau\in\nenatMlt(X)$ and stage $i\in\NNO$. The above three
options are captured as follows.
\begin{equation}
\label{NegMulnomDiag}
\vcenter{\xymatrix@R-3.5pc{
\nenatMlt(X)\times \NNO\ar[r]^-{\negmnclg(\omega)} & 
   \Dst\Big(\nenatMlt(X)\times\NNO \;+\; \NNO\Big)
\\
(\tau,i)\ar@{|->}[r] & 
   {\begin{array}[t]{l}
   \displaystyle\sum_{x,\,\tau(x)=0}\!\omega(x)
   \bigket{\tau, i\!+\!1} \,+\,
   \sum_{x,\,\tau(x)>0,\,\|\tau\|=1}\!\omega(x)\bigket{i\!+\!1}
   \\[+1em] \qquad + \displaystyle
   \sum_{x,\,\tau(x)>0,\,\|\tau\|>1}\!\omega(x)\bigket{\tau\!-\!\ket{x},i\!+\!1}.
   \end{array}}
}}
\end{equation}

\noindent Starting from initial position $(\tau,0)$ one eventually
ends up with an output in $\NNO$. We reason informally from the
contrapositive: an infinite sequences $(\rho,k) \rightarrow (\rho,k+1)
\rightarrow (\rho, k+2) \rightarrow \cdots$ exists only if the colours
$x\in\supp(\rho)$ are never drawn when the size of the draws
goes to infinity. This is impossible.

\subsection*{Negative hypergeometric distributions}

In the hypergeometric (and P\'olya) mode the urn changes with every
draw, so we have to incorporate not only the tubes $\tau$ but also the
urn $\upsilon$ in the positions of our MMO. For convenience, we
introduce the following special notation.
\[ \begin{array}{rcl}
\natMlt_{\geq *}(X)
& \coloneqq &
\setin{(\upsilon,\tau)}{\nenatMlt(X)\times\nenatMlt(X)}{\upsilon \geq \tau}.
\end{array} \]

\noindent The inequality $\upsilon\geq\tau$ expresses that the urn
contains sufficiently many balls of each colour to fill the
tubes. This inequality acts as an invariant for the following negative
hypergeometric MMO.
\begin{equation}
\label{NegHypgeomDiag}
\vcenter{\xymatrix@R-4.5pc@C-1pc{
\natMlt_{\geq *}(X)\times \NNO\ar[r]^-{\neghgclg} & 
   \Dst\Big(\natMlt_{\geq *}(X)\times\NNO + \NNO\Big)
\\
(\upsilon,\tau,i)\ar@{|->}[r] & 
   {\begin{array}[t]{l}
   \displaystyle\sum_{x,\,\tau(x)=0, \, \upsilon(x) > 0}\!\flrn(\upsilon)(x)
   \bigket{\upsilon\!-\!1\ket{x}, \tau, i\!+\!1}
   \\[+1.2em] 
   \qquad + \displaystyle
   \sum_{x,\,\tau(x)>0,\,\|\tau\|=1}\!\flrn(\upsilon)(x)\bigket{i\!+\!1}
   \\[+1.3em] 
   \qquad + \displaystyle
   \sum_{x,\,\tau(x)>0,\,\|\tau\|>1}\!\flrn(\upsilon)(x)
     \bigket{\upsilon\!-\!1\ket{x},\tau\!-\!\ket{x},i\!+\!1}.
   \end{array}}
}}
\end{equation}

\noindent It is obvious that there are no infinite transitions
starting from $(\upsilon, \tau, 0)$ since in each non-output step
the urn $\upsilon$ decreases in size.
\subsection*{Negative P\'olya} 

We now require that initially, the urn $\upsilon$ contains for all
colours of the tubes $\tau$ at least one ball. This can expressed as
inclusion of supports. Hence we define:
\[ \begin{array}{rcl}
\natMlt_{\supseteq *}(X)
& \coloneqq &
\setin{(\upsilon,\tau)}{\nenatMlt(X)\times\nenatMlt(X)}
   {\supp(\upsilon) \supseteq \supp(\tau)}.
\end{array} \]

\noindent This is used in the following negative P\'olya MMO.  It
looks very much like the hypergeometric one in~\eqref{NegHypgeomDiag},
with removal of balls from the urn replaced by addition.

\begin{equation}
\label{NegPolyaDiag}
\vcenter{\xymatrix@R-4.5pc@C-1.2pc{
\natMlt_{\supseteq *}(X)\times \NNO\ar[r]^-{\negplclg} & 
   \Dst\Big(\natMlt_{\supseteq *}(X)\times\NNO + \NNO\Big)
\\
(\upsilon,\tau,i)\ar@{|->}[r] & 
   {\begin{array}[t]{l}
   \displaystyle\sum_{x,\,\tau(x)=0}\!\flrn(\upsilon)(x)
   \bigket{\upsilon\!+\!1\ket{x}, \tau, i\!+\!1}
   \\[+1.2em] 
   \qquad + \displaystyle
   \sum_{x,\,\tau(x)>0,\,\|\tau\|=1}\!\flrn(\upsilon)(x)\bigket{i\!+\!1}
   \\[+1.3em] 
   \qquad + \displaystyle
   \sum_{x,\,\tau(x)>0,\,\|\tau\|>1}\!\flrn(\upsilon)(x)
     \bigket{\upsilon\!+\!1\ket{x},\tau\!-\!\ket{x},i\!+\!1}.
   \end{array}}
}}
\end{equation}

\noindent Also in this case there are no infinite transitions from an
initial position $(\upsilon,\tau,0)$, because the probability that
certain colours do not occur in P\'olya draws becomes zero as the size
of draws goes to infinity.

\section{Hypergeometric and P\'olya distributions via (negative) 
   binomials}\label{BinomialSec}

In Section~\ref{MulnomHypgeomPolyaSec} we have have introduced the
(ordinary, non-negative) hypergeometric and P\'olya distributions
$\hypergeometric[K](\upsilon)$ and $\polya[K](\upsilon)$, for an urn
$\upsilon$. It is known that these distributions can also be obtained
via conditioning, namely of parallel binomials in the hypergeometric
case, and of parallel negative binomials in the P\'olya case (see
\textit{e.g.}~\cite{SibuyaYS64}, for the bivariate, and
also~\cite{Janardan74} for the multivariate case). These conditionings
build on Propositions~\ref{HypgeomCountProp}, \ref{PolyaCountProp} and
fit very well in the current account, and are therefore included here,
in fully multivariate form. In order to do so we need to recall the
basics of probabilistic conditioning, see
\textit{e.g.}~\cite{Jacobs19b,Jacobs19c,Jacobs19d,Jacobs21b,Jacobs21a}
for more information.

Let $\omega\in\Dst(X)$ be distribution and $p\colon X \rightarrow
[0,1]$ be a (fuzzy) predicate. We write $\omega\models p \coloneqq
\sum_{x} \omega(x)\cdot p(x)$ for the validity (expected value) of $p$
in $\omega$. If this validity is non-zero, we can define the updated
distribution $\omega|_{p} \in\Dst(X)$ as the normalised product:
\[ \begin{array}{rcl}
\omega|_{p}(x)
& \coloneqq &
\displaystyle\frac{\omega(x)\cdot p(x)}{\omega\models p}.
\end{array} \]

\noindent See \textit{e.g.}~\cite{Jacobs19b,Jacobs19c,Jacobs21a} for
more information.

The two propositions below describe the two conditioning results for
hypergeometric and P\'olya distributions. They both use the following
sum predicate $\mathsl{sum}_{K} \colon \NNO^{\ell} \rightarrow [0,1]$,
for $K\in\NNO$.
\begin{ceqn}
\begin{equation}
\label{SumPredEqn}
\begin{array}{rcl}
\mathsl{sum}_{K}(n_{1}, \ldots, n_{\ell})
& \coloneqq &
\begin{cases}
1 & \mbox{if } n_{1} + \cdots + n_{\ell} = K
\\
0 & \mbox{otherwise.}
\end{cases}
\end{array}
\end{equation}
\end{ceqn}

\begin{proposition}
\label{HypgeomViaBinomProp}
Conditioning parallel binomials with this sum
predicate~\eqref{SumPredEqn} yields the hypergeometric distribution:
for $K \leq \sum_{i}k_{i}$,
\[ \begin{array}{rcl}
\big(\binomial[k_{1}](r) \otimes \cdots \otimes 
   \binomial[k_{\ell}](r)\big)\big|_{\mathsl{sum}_K}
& = &
\hypergeometric[K]\big(\sum_{i} k_{i}\ket{i}\big).
\end{array} \]

\noindent This works for any number $r\in [0,1]$.
\end{proposition}

\begin{myproof}
We first compute the validity:
\allowdisplaybreaks{\begin{alignat}{2}
\lefteqn{\binomial[k_{1}](r) \otimes \cdots \otimes \binomial[k_{\ell}](r)
   \models \mathsl{sum}_K}
\\
&\notag = 
\displaystyle\sum_{n_{i} \leq k_{i}} \,
   \big(\binomial[k_{1}](r) \otimes \cdots \otimes \binomial[k_{\ell}](r)\big)
      (n_{1}, \ldots, n_{\ell}) \cdot \mathsl{sum}_{K}(n_{1}, \ldots, n_{\ell})
\\
& \notag= 
\displaystyle\sum_{n_{i} \leq k_{i}, \, {\scriptscriptstyle\sum_{i}}n_{i} \,=\, K} \,
   \binomial[k_{1}](r)(n_{1}) \cdot \ldots \cdot \binomial[k_{\ell}](r)(n_{\ell})
\\
&\notag = 
\displaystyle\sum_{n_{i} \leq k_{i}, \, {\scriptscriptstyle\sum_{i}}n_{i} \,=\, K} \,
   \binom{k_1}{n_1}\cdot r^{n_1}\cdot (1\!-\!r)^{k_{1}-n_{1}} \cdot \ldots \cdot
   \binom{k_\ell}{n_\ell}\cdot r^{n_\ell}\cdot (1\!-\!r)^{k_{\ell}-n_{\ell}}
\\
&\notag = 
\displaystyle\sum_{n_{i} \leq k_{i}, \, {\scriptscriptstyle\sum_{i}}n_{i} \,=\, K} \,
   \binom{\sum_{i}k_{i}\ket{i}}{\sum_{i}n_{i}\ket{i}}\cdot r^{K}\cdot 
   (1\!-\!r)^{(\sum_{i}k_{i})-K}
\\
&\notag = 
\displaystyle\binom{\sum_{i}k_{i}}{K}\cdot r^{K}\cdot 
   (1\!-\!r)^{(\sum_{i}k_{i})-K} \qquad 
   \mbox{by Proposition~\ref{HypgeomCountProp}}
\\
&\notag = 
\binomial[{\scriptstyle\sum_{i}}k_{i}](r)(K).
\end{alignat}}

\noindent Now we can move on to the conditioning itself. We see that
the probability $r\in[0,1]$ drops out of the calculation.
\[\begin{array}{rcl}
\lefteqn{\big(\binomial[k_{1}](r) \otimes \cdots \otimes 
      \binomial[k_{\ell}](r)\big)\big|_{\mathsl{sum}_K}}
\\[0.4em]
& = &
\displaystyle\sum_{n_{i} \leq k_{i}} \,
\frac{\big(\binomial[k_{1}](r) \otimes \cdots \otimes 
   \binomial[k_{\ell}](r)\big)(\vec{n}) \cdot \mathsl{sum}_{K}(\vec{n})}
   {\binomial[k_{1}](r) \otimes \cdots \otimes \binomial[k_{\ell}](r)
   \models \mathsl{sum}_K} \, \bigket{n_{1}, \ldots, n_{\ell}}
\\[+1.2em]
& = &
\displaystyle\sum_{n_{i} \leq k_{i}, \, {\scriptscriptstyle\sum_{i}}n_{i} \,=\, K} 
   \frac{\binom{k_1}{n_1}\!\cdot\! r^{n_1}\!\cdot\! (1\!-\!r)^{k_{1}-n_{1}}\!\cdot 
   \ldots \cdot\! \binom{k_\ell}{n_\ell} 
   \!\cdot\! r^{n_\ell} \!\cdot\! (1\!-\!r)^{k_{\ell}-n_{\ell}}}
   {\binom{\sum_{i}k_{i}}{K}\cdot r^{K}\cdot (1\!-\!r)^{(\sum_{i}k_{i})-K}}
   \bigket{n_{1}, \ldots, n_{\ell}}
\\[+1.2em]
& = &
\displaystyle\sum_{n_{i} \leq k_{i}, \, {\scriptscriptstyle\sum_{i}}n_{i} \,=\, K} 
   \frac{\prod_{i}\binom{k_i}{n_i} \cdot r^{K}\cdot (1\!-\!r)^{(\sum_{i}k_{i})-K}}
   {\binom{\sum_{i}k_{i}}{K}\cdot r^{K}\cdot (1\!-\!r)^{(\sum_{i}k_{i})-K}}
   \bigket{n_{1}, \ldots, n_{\ell}}
\\[+1.2em]
& = &
\displaystyle\sum_{n_{i} \leq k_{i}, \, {\scriptscriptstyle\sum_{i}}n_{i} \,=\, K} 
   \frac{\prod_{i}\binom{k_i}{n_i}}{\binom{\sum_{i}k_{i}}{K}}
   \bigket{n_{1}, \ldots, n_{\ell}}
\\[+1.8em]
& \smash{\stackrel{\eqref{HypgeomEqn}}{=}} &
\hypergeometric[K]\Big(\sum_{i} k_{i}\ket{i}\Big).
\end{array}\]

\noindent In the last line we implicitly identify the sequence $n_{1},
\ldots, n_{\ell}$ with the multiset $n_{1}\ket{1} + \cdots +
n_{\ell}\ket{\ell}$. \QED
\end{myproof}

There is a similar result for P\'olya distributions, using negative
binomials. It requires some care since it involves a shift of
arguments, since negative distributions (on $\NNO$) start only after a
certain number of steps.

\begin{proposition}
\label{PolyaViaNegbinomProp}
The multivariate P\'olya distribution can be described as conditioning
of negative binomials: for $K \geq \sum_{i}k_{i}$,
\[ \begin{array}{rcl}
\lefteqn{\big(\negbinomial[k_{1}](r) \otimes \cdots \otimes 
   \negbinomial[k_{\ell}](r)\big)\big|_{\mathsl{sum}_K}\big(k_{1}\!+\!n_{1}, 
   \ldots, k_{\ell}\!+\!n_{\ell}\big)}
\\[+0.4em]
& = &
\begin{cases}
\polya[K\!-\!\sum_{i}k_{i}]\big(\sum_{i} k_{i}\ket{i}\big)
   \big(\sum_{i}n_{i}\ket{i}\big) 
   & \mbox{if } \sum_{i}k_{i}\!+\!n_{i} = K
\\
0 & \mbox{otherwise.}
\end{cases}
\end{array} \]

\noindent The number $r\in[0,1]$ is arbitrary, and the predicate
$\mathsl{sum}_K$ is from~\eqref{SumPredEqn}.
\end{proposition}

\begin{myproof}
We start with the validity:
\allowdisplaybreaks{\begin{alignat}{2}
\lefteqn{\big(\negbinomial[k_{1}](r) \otimes \cdots \otimes 
   \negbinomial[k_{\ell}](r)\big)\models\mathsl{sum}_K}
\\
&\notag = 
\displaystyle\sum_{n_{1},\ldots,n_{\ell}} \,
   \big(\negbinomial[k_{1}](r)\otimes\cdots\otimes\negbinomial[k_{\ell}](r)\big)
      (n_{1}, \ldots, n_{\ell}) \cdot \mathsl{sum}_{K}(n_{1}, \ldots, n_{\ell})
\\
&\notag =
\displaystyle\sum_{n_{i}, \, {\scriptscriptstyle\sum_{i}} k_{i}+n_{i} \,=\, K} \,
   \binomial[k_{1}](r)(k_{1}\!+\!n_{1}) \cdot \ldots \cdot 
   \binomial[k_{\ell}](r)(k_{\ell}\!+\!n_{\ell})
\\
&\notag \smash{\stackrel{\eqref{NegBinomEqn}}{=}} 
\displaystyle\sum_{n_{i}, \, {\scriptscriptstyle\sum_{i}} k_{i}+n_{i} \,=\,K} \, \textstyle
   {\displaystyle\prod}_{i}\, \displaystyle 
   \bibinom{k_i}{n_i}\cdot r^{k_i}\cdot (1\!-\!r)^{n_{i}}
\\
&\notag = 
\displaystyle\sum_{n_{i}, \, {\scriptscriptstyle\sum_{i}}n_{i} \,=\, 
     K-{\scriptscriptstyle\sum_{i}}k_{i}} \, \textstyle
   {\displaystyle\prod}_{i}\, \displaystyle 
   \bibinom{\sum_{i}k_{i}\ket{i}}{\sum_{i}n_{i}\ket{i}}
   \cdot r^{\sum_{i}k_i}\cdot (1\!-\!r)^{\sum_{i}n_{i}}
\\
&\notag = 
\displaystyle \bibinom{\sum_{i}k_{i}}{K-\sum_{i}k_{i}}
   \cdot r^{\sum_{i}k_i}\cdot (1\!-\!r)^{K-\sum_{i}k_{i}}
   \qquad \mbox{by Proposition~\ref{PolyaCountProp}}
\\
&\notag =
\negbinomial[{\scriptstyle\sum_{i}}k_{i}](r)(K).
\end{alignat}}

\noindent Conditioning yields, when $\sum_{i}k_{i}\!+\!n_{i} = K$,
\[ \hspace*{-3em}\begin{array}[b]{rcl}
\lefteqn{\big(\negbinomial[k_{1}](r) \otimes \cdots \otimes 
   \negbinomial[k_{\ell}](r)\big)\big|_{\mathsl{sum}_K}\big(k_{1}\!+\!n_{1}, 
   \ldots, k_{\ell}\!+\!n_{\ell}\big)}
\\[0.4em]
& = &
\displaystyle
\frac{\big(\binomial[k_{1}](r) \otimes \cdots \otimes 
   \binomial[k_{\ell}](r)\big)((k_{1}\!+\!n_{1}, \ldots, k_{\ell}\!+\!n_{\ell})}
   {\negbinomial[k_{1}](r) \otimes \cdots \otimes \negbinomial[k_{\ell}](r)
   \models \mathsl{sum}_K}
\\[+1.0em]
& = &
\displaystyle
\frac{\prod_{i}\,\big(\!\binom{k_i}{n_i}\!\big)\cdot r^{k_i}\cdot (1\!-\!r)^{n_{i}}}
    {\big(\!\binom{\sum_{i}k_{i}}{K-\sum_{i}k_{i}}\!\big)
   \cdot r^{\sum_{i}k_i}\cdot (1\!-\!r)^{K-\sum_{i}k_{i}}}
\\[+2.2em]
& = &
\displaystyle
\frac{\prod_{i}\,\big(\!\binom{k_i}{n_i}\!\big)}
    {\big(\!\binom{\sum_{i}k_{i}}{K-\sum_{i}k_{i}}\!\big)}
\\[+1.8em]
& \smash{\stackrel{\eqref{PolyaEqn}}{=}} &
\polya[K\!-\!\sum_{i}k_{i}]\big(\sum_{i} k_{i}\ket{i}\big)
   \big(\sum_{i}n_{i}\ket{i}\big).
\end{array} \eqno{\QEDbox} \]
\end{myproof}

\section{Number-theoretic corollaries}\label{CorollarySec}

In this final section we extract several number-theoretic equations
from the fact that first-full and negative probabilities form
distributions and thus add up to one. These equations are obtained
from the bivariate case, with two tubes. Recall, that the bivariate
situation studied in the literature involves one tube only. As far as
we know, the equations given below are new (or at least, not very
familiar). For some of them --- like
Corollaries~\ref{BinFirstFullCor}~\eqref{BinFirstFullCorMulnom}
and~\ref{BinNegativeCor}~\eqref{BinNegativeCorMulnom} --- the author
has direct proofs, but not for the others.

We first look at what follows from the bivariate first-full
distributions.

\begin{corollary}
\label{BinFirstFullCor}
Fix numbers $n > 0$ and $m > 0$.
\begin{enumerate}
\item \label{BinFirstFullCorMulnom}For $r,s\in [0,1]$ with $r + s = 1$
  one has:
\[ \begin{array}{rcl}
r^{n} \cdot \displaystyle\sum_{j<m}\, \bibinom{n}{j}\cdot s^{j}
\;+\;
s^{m} \cdot \displaystyle\sum_{i<n}\, \bibinom{m}{i}\cdot r^{i}
& = &
1.
\end{array} \]

\item \label{BinFirstFullCorHypgeom} For $N\geq n$ and $M \geq m$ one has:
\[ \begin{array}{rcl}
\displaystyle\sum_{j<m} \bibinom{n}{j}\cdot\binom{N\!-\!n+M\!-\!j}{N\!-\!n}
 \;+\;
   \sum_{i<n} \bibinom{m}{i}\cdot\binom{N\!-\!i+M\!-\!m}{M\!-\!m}
& = &
\displaystyle\binom{N+M}{N}. 
\end{array} \]

\item \label{BinFirstFullCorPolya} For $N>0$ and $M>0$,
\[ \begin{array}{rcl}
\displaystyle n\cdot\bibinom{N}{n}\cdot\sum_{j< m} 
   \frac{\big(\!\binom{M}{j}\!\big)}
      {\big(\!\binom{n+j}{N+M}\!\big)} \,+\,
   m\cdot\bibinom{M}{m}\cdot\sum_{i<n} 
   \frac{\big(\!\binom{N}{i}\!\big)}
      {\big(\!\binom{i+m}{N+M}\!\big)}
& = &
N\!+\!M.
\end{array} \]
\end{enumerate}
\end{corollary}

\begin{myproof}
Take a binary space $X = \{a,b\}$ with tubes $\tau = n\ket{a} +
m\ket{b}$.
\begin{enumerate}
\item The numbers $r,s$ form a state $\omega = r\ket{a} + s\ket{b}$.
  Since $\mnff(\omega,\tau)$ is a distribution on $X$ we get
  $\mnff(\omega,\tau)(a) + \mnff(\omega,\tau)(b) = 1$. According to
  Definition~\ref{FirstFullDstDef}~\eqref{FirstFullDstDefMulnom} this
  means:

\allowdisplaybreaks{\begin{alignat}{2}
1
&\notag = 
\displaystyle
   \sum_{j<m}\, \mulnom(\omega)((n\!-\!1)\ket{a} + j\ket{b})\cdot\omega(a) 
\\
&\notag  \hspace*{4em} +\; \displaystyle
   \sum_{i<n}\, \mulnom(\omega)(i\ket{a} + (m\!-\!1)\ket{b})\cdot\omega(b) 
\\
&\notag =
\displaystyle
   \sum_{j<m}\, \binom{n\!-\!1\!+\!j}{n\!-\!1}\cdot r^{n-1}\cdot s^{j}\cdot r
   \;+\;
   \sum_{i<n}\, \binom{i\!+\!m\!-\!1}{i}\cdot r^{i}\cdot s^{m-1}\cdot s
\\
&\notag = 
\displaystyle
   r^{n} \cdot \sum_{j<m}\, \bibinom{n}{j}\cdot s^{j}
   \;+\;
   s^{m} \cdot \sum_{i<n}\, \bibinom{m}{i}\cdot r^{i}.
\end{alignat}}

\item For urn $\upsilon = N\ket{a} + M\ket{b}$ we have a
  hypergeometric first-fill distribution $\hgff(\upsilon,\tau)$, so
  that by
  Definition~\ref{FirstFullDstDef}~\eqref{FirstFullDstDefHypgeom}:
\[ \hspace*{-0.2em}\begin{array}{rcl}
1
& = &
\hgff(\upsilon,\tau)(a) + \hgff(\upsilon,\tau)(b)
\\[+0.3em]
& = &
\displaystyle
   \sum_{j<m}\, \hypgeom(\upsilon)((n\!-\!1)\ket{a} + j\ket{b})\cdot
      \flrn\big(\upsilon-(n\!-\!1)\ket{a} - j\ket{b}\big)(a)
\\[-0.5em]
& & \hspace*{2em} +\, \displaystyle
   \sum_{i<n}\, \hypgeom(\upsilon)(i\ket{a} + (m\!-\!1)\ket{b})\cdot
      \flrn\big(\upsilon-i\ket{a} - (m\!-\!1)\ket{b}\big)(b)
\\[+1.4em]
& = &
\displaystyle
   \sum_{j<m}\, \frac{\binom{N}{n-1}\cdot\binom{M}{j}}{\binom{N+M}{n-1+j}}\cdot
     \frac{N\!-\!n\!+\!1}{N\!+\!M\!-\!n\!+\!1\!-\!j}
\\[+1.2em]
& & \hspace*{1em} +\, \displaystyle
   \sum_{i<n}\, \frac{\binom{N}{i}\cdot\binom{M}{m-1}}{\binom{N+M}{i+m-1}}\cdot
     \frac{M\!-\!m\!+\!1}{N\!+\!M\!-\!i\!-\!m\!+\!1}
\\[+1.4em]
& = &
\displaystyle
   \sum_{j<m} \frac{N!\cdot M!\cdot (n\!-\!1\!+\!j)!\cdot 
    (N\!+\!M\!-\!n\!+\!1\!-\!j)!}
   {(n\!-\!1)!\cdot(N\!-\!n\!+\!1)!\cdot j!\cdot (M\!-\!j)! \cdot (N\!+\!M)!}
      \cdot\frac{N\!-\!n\!+\!1}{N\!+\!M\!-\!n\!+\!1\!-\!j}
\\[+1.2em]
& & +\, \displaystyle
   \sum_{i<n} \frac{N!\cdot M!\cdot (i\!+\!m\!-\!1)!\cdot 
    (N\!+\!M\!-\!i\!-\!m\!+\!1)!}
   {i!\cdot (N\!-\!i)! \cdot (m\!-\!1)!\cdot(M\!-\!m\!+\!1)!\cdot (N\!+\!M)!}
      \cdot\frac{M\!-\!m\!+\!1}{N\!+\!M\!-\!i\!-\!m\!+\!1}
\\[+1.0em]
& = &
\displaystyle \frac{1}{\binom{N+M}{N}}\cdot\left(
   \sum_{j<m} \bibinom{n}{j}\cdot\binom{N\!-\!n+M\!-\!j}{N\!-\!n}
   \;+\;
   \sum_{i<n} \bibinom{m}{i}\cdot\binom{N\!-\!i+M\!-\!m}{M\!-\!m}\right).
\end{array} \]

\item By rewriting the equation $\plff(\upsilon,\tau)(a) +
  \plff(\upsilon,\tau)(b) = 1$ in a similar manner. \QED

\auxproof{
Similarly, by Definition~\ref{FirstFullDstDef}~\eqref{FirstFullDstDefPolya}:
\[ \begin{array}{rcl}
1
& = &
\plff(\upsilon,\tau)(a) + \plff(\upsilon,\tau)(b)
\\[+0.3em]
& = &
\displaystyle
   \sum_{j<m}\, \pol(\upsilon)((n\!-\!1)\ket{a} + j\ket{b})\cdot
      \flrn\big(\upsilon+(n\!-\!1)\ket{a} + j\ket{b}\big)(a)
\\[-0.5em]
& & \hspace*{2em} +\, \displaystyle
   \sum_{i<n}\, \pol(\upsilon)(i\ket{a} + (m\!-\!1)\ket{b})\cdot
      \flrn\big(\upsilon+i\ket{a} + (m\!-\!1)\ket{b}\big)(b)
\\[+1.4em]
& = &
\displaystyle
   \sum_{j<m}\, \frac{\big(\!\binom{N}{n-1}\!\big)\cdot
      \big(\!\binom{M}{j}\!\big)}{\big(\!\binom{N+M}{n-1+j}\!\big)}\cdot
     \frac{N\!+\!n\!-\!1}{N\!+\!M\!+\!n\!-\!1\!+\!j}
\\[+1.2em]
& & \hspace*{1em} +\, \displaystyle
   \sum_{i<n}\, \frac{\big(\!\binom{N}{i}\!\big)\cdot
      \big(\!\binom{M}{m-1}\!\big)}{\big(\!\binom{N+M}{i+m-1}\!\big)}\cdot
     \frac{M\!+\!m\!-\!1}{N\!+\!M\!+\!i\!+\!m\!-\!1}
\\[+1.4em]
& = &
\displaystyle
   \sum_{j<m} \frac{n\cdot\big(\!\binom{N}{n}\!\big)\cdot
   \big(\!\binom{M}{j}\!\big) \cdot (n\!-\!1\!+\!j)!\cdot 
    (N\!+\!M\!-\!1)!}
   {(N\!+\!M\!+\!n\!-\!2\!+\!j)!}
      \cdot\frac{1}{N\!+\!M\!+\!n\!-\!1\!+\!j}
\\[+1.2em]
& & +\, \displaystyle
   \sum_{i<n} \frac{\big(\!\binom{N}{i}\!\big)\cdot
      m\cdot\big(\!\binom{M}{m}\!\big)\cdot (i\!+\!m\!-\!1)!\cdot 
    (N\!+\!M\!-\!1)!}
   {(N\!+\!M\!+\!i\!+\!m\!-\!2)!}
      \cdot\frac{1}{N\!+\!M\!+\!i\!+\!m\!-\!1}
\\[+1.4em]
& = &
\displaystyle\frac{1}{N\!+\!M}\cdot\left(
   n\cdot\bibinom{N}{n}\cdot\sum_{j< m} 
   \frac{\big(\!\binom{M}{j}\!\big)}
      {\big(\!\binom{n+j}{N+M}\!\big)} \,+\,
   m\cdot\bibinom{M}{m}\cdot\sum_{i< n} 
   \frac{\big(\!\binom{N}{i}\!\big)}
      {\big(\!\binom{i+m}{N+M}\!\big)}\right).
\end{array} \]
}
\end{enumerate}
\end{myproof}

Next we look at the consequences of having (bivariate) negative distributions.

\begin{corollary}
\label{BinNegativeCor}
Let arbitrary numbers $n>0$ and $m>0$ be given.
\begin{enumerate}
\item \label{BinNegativeCorMulnom} For probabilities $r,s\in(0,1)$
  with $r+s=1$ one has:
\[ \begin{array}{rcl}
\displaystyle\sum_{i\geq 0}\, \bibinom{n}{m+i}\cdot s^{i} \,+\,
   \bibinom{m}{n+i}\cdot r^{i}
& = &
\displaystyle\frac{1}{r^{n}\cdot s^{m}}.
\end{array} \]

\item \label{BinNegativeCorHypgeom} For $N\geq n$ and $M\geq m$,
\[ \begin{array}{l}
\displaystyle \sum_{j\leq M-m} \bibinom{n}{m\!+\!j}\cdot
   \binom{N\!-\!n+M\!-\!m\!-\!j}{N\!-\!n}
\\
\qquad + \; \displaystyle
   \sum_{i\leq N-n} \bibinom{m}{n\!+\!i}\cdot
   \binom{N\!-\!n\!-\!i+M\!-\!m}{M\!-\!m}
\hspace*{\arraycolsep}=\hspace*{\arraycolsep}
\displaystyle\binom{N+M}{N}.
\end{array} \]

\item \label{BinNegativeCorPolya} For $N > 0$ and $M > 0$ we have:
\[ \begin{array}{rcl}
\displaystyle n\cdot\bibinom{N}{n}\cdot\sum_{j\geq m} 
   \frac{\big(\!\binom{M}{j}\!\big)}
      {\big(\!\binom{n+j}{N+M}\!\big)} \,+\,
   m\cdot\bibinom{M}{m}\cdot\sum_{i\geq n} 
   \frac{\big(\!\binom{N}{i}\!\big)}
      {\big(\!\binom{i+m}{N+M}\!\big)}
& = &
N\!+\!M.
\end{array} \]
\end{enumerate}
\end{corollary}

\begin{myproof}
We only do the first and third item and leave the second one to the
interested reader. Like in the proof of
Corollary~\ref{BinFirstFullCor} we fix a space $X = \{a,b\}$ with
tubes $\tau = n\ket{a} + m\ket{b}$.  For state $\omega = r\ket{a} +
s\ket{b}$ we use that the negative multinomial
$\negmultinomial(\omega, \tau)$ is a distribution on $\NNO$ and unpack
its description from
Definition~\ref{NegativeDstDef}~\eqref{NegativeDstDefMulnom}.
\[ \begin{array}[b]{rcl}
1
& = &
\displaystyle\sum_{i\geq 0}\, \negmultinomial(\omega, \tau)(i)
\\[+1.2em]
& = &
\displaystyle\sum_{i\geq 0}\, 
   \mulnom(\omega)((n\!-\!1)\ket{a}+(m\!+\!i)\ket{b})\cdot \omega(a) 
   \;+\;
   \mulnom(\omega)((n\!+\!i)\ket{a}+(m\!-\!1)\ket{b})\cdot \omega(b)
\\[+0.4em]
& = &
\displaystyle\sum_{i\geq 0}\, 
   \binom{n\!+\!i\!+\!m\!-\!1}{n\!-\!1} \cdot r^{n-1}\cdot s^{m+i} \cdot r
   \,+\,
   \binom{n\!+\!i\!+\!m\!-\!1}{n\!+\!i} \cdot r^{n+i}\cdot s^{m-1} \cdot s
\\[+1.2em]
& = &
\displaystyle r^{n}\cdot s^{m} \cdot 
   \sum_{i\geq 0}\, \bibinom{n}{m+i}\cdot s^{i} \,+\, \bibinom{m}{n+i}\cdot r^{i}.
\end{array} \]

\auxproof{
For item~\eqref{BinNegativeCorHypgeom} We use urn $\upsilon = N\ket{a}
+ M\ket{b}$ in:
\[ \begin{array}{rcl}
1
& = &
\displaystyle\sum_{j\leq M-m} 
   \hypgeom(\upsilon)((n\!-\!1)\ket{a} + (m\!+\!j)\ket{b})
\\[-0.5em]
& & \hspace*{6em} \cdot\, 
   \flrn(\upsilon - (n\!-\!1)\ket{a} - (m\!+\!j)\ket{b})(a)
\\[+0.5em]
& & \hspace*{3em} +\, \displaystyle 
   \sum_{i\leq N-n} 
   \hypgeom(\upsilon)((n\!+\!i)\ket{a} + (m\!-\!1)\ket{b})
\\[-0.5em]
& & \hspace*{9em} \cdot\, 
   \flrn(\upsilon - (n\!+\!i)\ket{a} - (m\!-\!1)\ket{b})(b)
\\[+0.5em]
& = &
\displaystyle\sum_{j\leq M-m} 
   \frac{\binom{N}{n-1}\cdot\binom{M}{m+j}}{\binom{N+M}{n-1+m+j}} \cdot
   \frac{N\!-\!n\!+\!1}{N\!+\!M\!-\!n\!+\!1\!-\!m\!-\!j}
\\[+1.3em]
& & \hspace*{3em} +\, \displaystyle 
   \sum_{i\leq N-n} 
   \frac{\binom{N}{n+i}\cdot\binom{M}{m-1}}{\binom{N+M}{n+i+m-1}} \cdot
   \frac{M\!-\!m\!+\!1}{N\!+\!M\!-\!n\!-\!i\!-\!m\!+\!1}
\\[+1.4em]
& = &
\displaystyle\sum_{j\leq M-m} 
   \frac{N!\cdot M!\cdot (n\!-\!1\!+\!m\!+\!j)!\cdot 
   (N\!+\!M\!-\!n\!+\!1\!-\!m\!-\!j)!}
      {(n\!-\!1)!\cdot (N\!-\!n\!+\!1)! \cdot (m\!+\!j)!\cdot 
   (M\!-\!m\!-\!j)!\cdot (N\!+\!M)!} \cdot
   \frac{N\!-\!n\!+\!1}{N\!+\!M\!-\!n\!+\!1\!-\!m\!-\!j}
\\[+1.3em]
& & +\, \displaystyle 
   \sum_{j\leq M-m} 
   \frac{N!\cdot M!\cdot (n\!+\!i\!+\!m\!-\!1)!\cdot 
   (N\!+\!M\!-\!n\!-\!i\!-\!m\!+\!1)!}
      {(n\!+\!i)!\cdot (N\!-\!n\!-\!i)! \cdot (m\!-\!1)!\cdot 
   (M\!-\!m\!+\!1)!\cdot (N\!+\!M)!} \cdot
   \frac{M\!-\!m\!+\!1}{N\!+\!M\!-\!n\!-\!i\!-\!m\!+\!1}
\\[+1.4em]
& = &
\displaystyle\frac{1}{\binom{N+M}{N}}\cdot
   \left(\sum_{i\leq N-n} \bibinom{n}{m\!+\!j}\cdot
   \binom{N\!-\!n+M\!-\!m\!-\!j}{N\!-\!n} \,+\,
   \sum_{i\leq N-n} \bibinom{m}{n\!+\!i}\cdot
   \binom{N\!-\!n\!-\!i+M\!-\!m}{M\!-\!m}\right).
\end{array} \]
}

\noindent For item~\eqref{BinNegativeCorPolya} we use:
\[ \begin{array}[b]{rcl}
1
& = &
\displaystyle\sum_{i\geq 0} 
   \pol(\upsilon)((n\!-\!1)\ket{a} + (m\!+\!j)\ket{b}) \cdot
   \flrn(\upsilon + (n\!-\!1)\ket{a} + (m\!+\!j)\ket{b})(a)
\\[-0.5em]
& & \hspace*{3em} +\; \displaystyle 
   \pol(\upsilon)((n\!+\!i)\ket{a} + (m\!-\!1)\ket{b}) \cdot
   \flrn(\upsilon + (n\!+\!i)\ket{a} + (m\!-\!1)\ket{b})(b)
\\[+0.6em]
& = &
\displaystyle\sum_{i\geq 0} 
   \frac{\big(\!\binom{N}{n-1}\!\big)\cdot\big(\!\binom{M}{m+i}\!\big)}
      {\big(\!\binom{N+M}{n-1+m+i}\!\big)} \cdot
   \frac{N\!+\!n\!-\!1}{N\!+\!M\!+\!n\!-\!1\!+\!m\!+\!i}
\\[+1.3em]
& & \hspace*{3em} +\, \displaystyle 
   \frac{\big(\!\binom{N}{n+i}\!\big)\cdot\big(\!\binom{M}{m-1}\!\big)}
      {\big(\!\binom{N+M}{n+i+m-1}\!\big)} \cdot
   \frac{M\!+\!m\!-\!1}{N\!+\!M\!+\!n\!+\!i\!+\!m\!-\!1}
\\[+1.4em]
& = &
\displaystyle\sum_{i\geq 0} 
   \frac{n\cdot\big(\!\binom{N}{n}\!\big)\cdot\big(\!\binom{M}{m+i}\!\big)\cdot 
   (n\!-\!1\!+\!m\!+\!i)!\cdot (N\!+\!M\!-\!1)!}
      {(N\!+\!M\!+\!n\!-\!2\!+\!m\!+\!i)!} \cdot
   \frac{1}{N\!+\!M\!+\!n\!-\!1\!+\!m\!+\!i}
\\[+1.3em]
& & \hspace*{1em}+\, \displaystyle 
   \frac{\big(\!\binom{N}{n+i}\!\big)\cdot m\cdot\big(\!\binom{M}{m}\!\big)\cdot
      (n\!+\!i\!+\!m\!-\!1)!\cdot (N\!+\!M\!-\!1)!}
      {(N\!+\!M\!+\!n\!+\!i\!+\!m\!-\!2)!} \cdot
   \frac{1}{N\!+\!M\!+\!n\!+\!i\!+\!m\!-\!1}
\\[+1.4em]
& = &
\displaystyle\frac{1}{N\!+\!M}\cdot\left(
   n\cdot\bibinom{N}{n}\cdot\sum_{j\geq m} 
   \frac{\big(\!\binom{M}{j}\!\big)}
      {\big(\!\binom{n+j}{N+M}\!\big)} \,+\,
   m\cdot\bibinom{M}{m}\cdot\sum_{i\geq n} 
   \frac{\big(\!\binom{N}{i}\!\big)}
      {\big(\!\binom{i+m}{N+M}\!\big)}\right).
\end{array} \eqno{\QEDbox} \]
\end{myproof}

\section{Conclusions}

This paper extends the familiar urn model to an urn \& tubes model. It
raises several research questions, with possible applications in risk
modeling.  The extension is first used to introduce first-full
distributions, which have a historical basis in the `problem of
points' of Pascal and Fermat. Next, the urn \& tubes models is used
for negative distributions. The contribution of this paper lies in
systematisation, via a clear model, formalised via multisets (for
urns, draws, tubes), covering the three main drawing modes
(multinomial, hypergeometric, P\'olya).

This paper concentrates on the conceptual basis, formalisation, and
illustration of first-full and negative distributions. There is more
to say, for instance about associated statistical properties like mean
and (co)variance. They exist in the literature for the single tube
case. Extension to general tubes is a challenge that is left open
here.

The urn \& tubes model may be generalised, for instance to multiple
urns, where there is a choice from which urn one wishes to draw a
ball. When the contents of the urns are known, one can consider
different strategies for such choices, in different drawing modes, via
Markov decision processes (see
\textit{e.g.}~\cite{BaierK08,Puterman94}). When the contents are
unknown, the setting may be used for reinforcement
learning~\cite{KaelblingLM96,SuttonB18}: jointly learning these
contents and developing a strategy.

\medskip

\appendix

Section~\ref{FirstFullDstSec} and Subsection~\ref{NegativeDstSubsec}
use a compositional approach for showing that first-full and negative
probabilities add up to one, and thus form proper distributions.  The
distributions appear after iteratively self-composing a Markov model
with output, in the form of a coalgebra $c\colon Y \rightarrow
\Dst(Y+X)$, see~\eqref{MMODiag}. There is a little bit of category
theory underlying this composition, which we make explicit in this
appedix.

It is well known that the mapping $A \mapsto \Dst(A)$, sending a set $A$
to the set of (discrete, finite) probability distributions on $A$ is
a monad, on the category $\Sets$ of sets and functions. The unit $\eta$
and multiplication $\mu$ of this monad $\Dst$ are:
\[ \xymatrix@R-2.5pc{
A\ar[r]^-{\eta} & \Dst(A)
& &
\Dst\big(\Dst(A)\big)\ar[r]^-{\mu} & \Dst(A)
\\
a\ar@{|->}[r] & 1\bigket{a}
& &
\sum_{i} r_{i}\bigket{\omega_i}\ar@{|->}[r] &
\displaystyle\sum_{a\in A} \Big(\textstyle\sum_{i} 
   r_{i} \cdot \omega_{i}(a)\Big) \bigket{a}.
} \]

\noindent We write $A+B$ for the coproduct (disjoint union) of two
sets $A,B$, with coprojections $A \stackrel{\kappa_1}{\longrightarrow}
A+B \stackrel{\kappa_2}{\longleftarrow} B$, and cotuple $[f,g] \colon
A+B \rightarrow C$, for $f\colon A \rightarrow C$, $g\colon B
\rightarrow C$. For a fixed set $X$, the mapping $A \mapsto A+X$ is
also a monad with unit $\kappa_{1} \colon A \rightarrow A+X$ and
multiplication $[\idmap, \kappa_{2}] \colon (A + X) + X \rightarrow
A+X$.

These two monads $\Dst$ and $(-)+X$ are connected via a distributive
law, of the form:
\[ \xymatrix{
\Dst(A) + X\ar[r]^-{\lambda} & \Dst(A+X)
\qquad\mbox{namely}\qquad
{\begin{array}{rcl}
\lambda
& = &
\big[\Dst(\kappa_{1}), \, \eta \after \kappa_{2}\big].
\end{array}}
} \]

\noindent A general categorical result, see
\textit{e.g.}~\cite{BarrW85}, now says that the composite
$\Dst\big((-)+X\big)$ is then also a monad. In particular, if we have
maps $c\colon A \rightarrow \Dst\big(B+X\big)$ and $d\colon B
\rightarrow \Dst\big(C+X\big)$ we can form a composition $d \klafter c
\colon A \rightarrow \Dst\big(C+X\big)$ as:
\[ \xymatrix@C-1pc@R-1pc{
d \klafter c \;=\;\Big(A\ar[r]^-{c} &
   \Dst\big(B\!+\!X\big)\ar[rr]^-{\Dst(d+\idmap)} & &
   \Dst\big(\Dst\big(C\!+\!X\big) \!+\! X\big)\ar[d]_-{\Dst(\lambda)}
\\
& & & \Dst\big(\Dst\big((C\!+\!X) \!+\! X\big)\big)
   \ar[rrr]^-{\Dst(\Dst([\idmap, \kappa_{2}]))} & & &
   \Dst\big(\Dst\big(C\!+\!X\big)\big)\ar[r]^-{\mu} & \Dst\big(C\!+\!X\big)\Big)
} \]

Now assume that we have a Markov model with output (MMO) $c\colon Y
\rightarrow \Dst\big(Y+X\big)$. We can form self-composites $c^{n}
\colon Y \rightarrow \Dst\big(Y+X\big)$, for $n\in\NNO$, in the
following manner:
\[ \begin{array}{rclcrcl}
c^{0}
& = &
\eta \after \kappa_{1}
& \qquad\mbox{and}\qquad &
c^{n+1}
& = &
c^{n} \klafter c.
\end{array} \]

\noindent By elaborating the details we get the self-composition
formulas~\eqref{IterationEqn} used for MMOs.


\end{document}